\documentclass[11pt]{article}
\usepackage{amsfonts,amsmath,amssymb}
\usepackage{graphicx}
\usepackage{hyperref}
\usepackage{color}

\newcommand{\status}{}

\newcommand{\file}{$\ti{\ }$/wisk/points/revis/netpoints\underline{\ }rev.tex
\quad}
\renewcommand{\file}{}

\newcommand{\plaat}[2]{{#2}}

\newcommand{\detail}[1]{\par\noi{\bf[Proof detail\ }{#1}
\hfill{\bf ]}\par\noi\hspace{-4pt}}
\renewcommand{\detail}[1]{}

\newcommand{\dis}{\displaystyle}
\newcommand{\txt}{\textstyle}

\newcommand{\noi}{\noindent}

\newcommand{\halmos}{\rule{1ex}{1.4ex}}
\def \qed {\nopagebreak{\hspace*{\fill}$\halmos$\medskip}}
\newcommand{\med}{\medskip}

\newtheorem{theorem}{Theorem}[section]
\newtheorem{proposition}[theorem]{Proposition}
\newtheorem{corollary}[theorem]{Corollary}
\newtheorem{conjecture}[theorem]{Conjecture}
\newtheorem{lemma}[theorem]{Lemma}

\newtheorem{remark}[theorem]{Remark}
\newtheorem{defi}[theorem]{Definition}

\newcommand{\bd}{\begin{defi}}
\newcommand{\ed}{\end{defi}}
\newcommand{\bt}{\begin{theorem}}
\newcommand{\et}{\end{theorem}}
\newcommand{\bl}{\begin{lemma}}
\newcommand{\el}{\end{lemma}}
\newcommand{\bp}{\begin{proposition}}
\newcommand{\ep}{\end{proposition}}
\newcommand{\bcor}{\begin{corollary}}
\newcommand{\ecor}{\end{corollary}}
\newcommand{\br}{\begin{remark}\rm}
\newcommand{\er}{\end{remark}}
\newcommand{\bcon}{\begin{conjecture}}
\newcommand{\econ}{\end{conjecture}}

\renewcommand{\theequation}{\thesection .\arabic{equation}}

\newcommand{\be}{\begin{equation}}
\newcommand{\ee}{\end{equation}}
\newcommand{\ben}{\begin{equation*}}
\newcommand{\een}{\end{equation*}}
\newcommand{\bea}{\begin{eqnarray}}
\newcommand{\eea}{\end{eqnarray}}
\newcommand{\bean}{\begin{eqnarray*}}
\newcommand{\eean}{\end{eqnarray*}}
\newcommand{\ba}{\begin{array}}
\newcommand{\ea}{\end{array}}
\newcommand{\bc}{\be\begin{array}{r@{\,}c@{\,}l}}
\newcommand{\ec}{\end{array}\ee}

\newcommand{\de}{\delta}
\newcommand{\De}{\Delta}
\newcommand{\eps}{\varepsilon}

\newcommand{\sig}{\sigma}

\newcommand{\si}{\ensuremath{\sigma}}

\newcommand{\Bi}{{\cal B}}
\newcommand{\Ci}{{\cal C}}
\newcommand{\Di}{{\cal D}}
\newcommand{\Ei}{{\cal E}}
\newcommand{\Fi}{{\cal F}}

\newcommand{\Hi}{{\cal H}}
\newcommand{\Ii}{{\cal I}}
\newcommand{\Ki}{{\cal K}}
\newcommand{\Li}{{\cal L}}

\newcommand{\Ni}{{\cal N}}

\newcommand{\Ri}{{\cal R}}

\newcommand{\Ti}{{\cal T}}
\newcommand{\Ui}{{\cal U}}

\newcommand{\Wi}{{\cal W}}

\newcommand{\R}{{\mathbb R}}
\newcommand{\N}{{\mathbb N}}
\newcommand{\Z}{{\mathbb Z}}

\newcommand{\Q}{{\mathbb Q}}

\renewcommand{\P}{{\mathbb P}}
\newcommand{\E}{{\mathbb E}}

\newcommand{\volgt}{\ensuremath{\Rightarrow}}
\newcommand{\up}{\uparrow}
\newcommand{\down}{\downarrow}
\newcommand{\sub}{\subset}
\newcommand{\beh}{\backslash}

\newcommand{\asto}[1]{\underset{{#1}\to\infty}{\longrightarrow}}
\newcommand{\Asto}[1]{\underset{{#1}\to\infty}{\Longrightarrow}}

\newcommand{\Aston}[1]{\underset{{#1}\to 0}{\Longrightarrow}}

\newcommand{\ti}{\tilde}

\newcommand{\ffrac}[2]{{\textstyle\frac{{#1}}{{#2}}}}
\newcommand{\dif}[1]{\ffrac{\partial}{\partial{#1}}}
\newcommand{\diff}[1]{\ffrac{\partial^2}{{\partial{#1}}^2}}
\newcommand{\difif}[2]{\ffrac{\partial^2}{\partial{#1}\partial{#2}}}

\newcommand{\di}{\mathrm{d}}
\newcommand{\half}{{[0,\infty)}}

\begin{document}

\makeatletter\@addtoreset{equation}{section}
\makeatother\def\theequation{\thesection.\arabic{equation}}

\renewcommand{\labelenumi}{{(\roman{enumi})}}

\newcommand{\Rc}{R^2_{\rm c}}
\newcommand{\Wl}{\Wi^{\rm l}}
\newcommand{\Wr}{\Wi^{\rm r}}
\newcommand{\El}{\Ei^{\rm l}}
\newcommand{\Er}{\Ei^{\rm r}}

\title{\vspace*{-1cm}Special points of the Brownian net}
\author{Emmanuel~Schertzer
\and
Rongfeng~Sun
\and
Jan~M.~Swart
}

\date{{\small\file} December 19, 2008}

\maketitle\vspace{-.8cm}
\status

\begin{abstract}\noi
The Brownian net, which has recently been introduced by Sun and
Swart~\cite{SS06}, and independently by Newman, Ravishankar and
Schertzer~\cite{NRSc08}, generalizes the Brownian web by allowing
branching. In this paper, we study the structure of the Brownian net
in more detail. In particular, we give an almost sure classification
of each point in $\R^2$ according to the configuration of the Brownian
net paths entering and leaving the point. Along the way, we establish
various other structural properties of the Brownian net.
\end{abstract}

\vspace{.3cm}

\noi
\vspace{12pt}

\noi
{\it Acknowledgement.} R.~Sun is supported by a postdoc position in the DFG
Forschergruppe 718 ``Analysis and Stochastics in Complex Physical
  Systems''. J.M.~Swart is sponsored by GA\v CR grant 201/07/0237. E.~Schertzer
and J.M.~Swart thank the Berlin-Leipzig Forschergruppe for hospitality and
support during a visit to the TU Berlin.

{\setlength{\parskip}{-2pt}\tableofcontents}

\section{Introduction and results}

\subsection{Introduction}

The Brownian web, $\Wi$, is essentially a collection of one-dimensional
coalescing Brownian motions starting from every point in space and time
$\R^2$. It originated from the work of Arratia \cite{A79,A81} on the scaling
limit of the voter model, and arises naturally as the diffusive scaling limit
of the system of one-dimensional coalescing random walks dual to the voter
model; see also \cite{FINR04} and \cite{NRSu05}. In the language of stochastic
flows, the coalescing flow associated with the Brownian web is known as the
{\em Arratia flow}. A detailed analysis of the Brownian web was carried out by
T\'oth and Werner in \cite{TW98}. More recently, Fontes, Isopi, Newman and
Ravishankar \cite{FINR04} introduced a by now standard framework in which the
Brownian web is regarded as a random compact set of paths, which (in a
suitable topology) becomes a random variable taking values in a Polish
space. It is in this framework that the object initially proposed by Arratia
in \cite{A81} takes on the name {\em the Brownian web}.

Recently, Sun and Swart \cite{SS06} introduced a generalization of the
Brownian web, called {\em the Brownian net}, $\Ni$, in which paths not only
coalesce, but also branch. From a somewhat different starting point, Newman,
Ravishankar and Schertzer independently arrived at the same object. Their
alternative construction of the the Brownian net will be published in
\cite{NRSc08}. The motivation in \cite{SS06} comes from the study of the
diffusive scaling limit of one-dimensional branching-coalescing random walks
with weak branching, while the motivation in \cite{NRSc08} comes from the
study of one-dimensional stochastic Ising and Potts models with boundary
nucleation. The different constructions of the Brownian net given in
\cite{SS06} and \cite{NRSc08} complement each other and give different
insights into the structure of the Brownian~net.

In the Brownian web, at a typical, deterministic point in $\R^2$, there is
just a single path leaving the point and no path entering the point. There
are, however, random, {\em special} points, where more than one path leaves,
or where paths enter. A full classification of these special points is given
in \cite{TW98}, see also \cite{FINR06}. The special points of the Brownian web
play an important role in the construction of the so-called {\em marked
  Brownian web} \cite{FINR06}, and also in the construction of the Brownian
net in \cite{NRSc08}. A proper understanding of the Brownian net thus calls
for a similar classification of special points of the Brownian net, which is
the main goal of this paper. Along the way, we will establish various
properties for the Brownian net $\Ni$, and the left-right Brownian web
$(\Wl,\Wr)$, which is the key intermediate object in the construction of the
Brownian net in \cite{SS06}.

Several models have been studied recently which have close connections to the
Brownian net. One such model is the so-called {\em dynamical Browian web},
which is a Browian web evolving in time in such a way that at random times,
paths switch among outgoing trajectories at points with one incoming, and two
outgoing paths. Such a model is similar in spirit to {\em dynamical
  percolation}, see e.g.\ \cite{H98}. In \cite{HW07}, Howitt and Warren
characterized the two-dimensional distributions of the dynamical Brownian
web. This leads to two coupled Brownian webs which are similar in spirit to
the left-right Brownian web $(\Wl,\Wr)$ in \cite{SS06}. Indeed, there is a
close connection between these objects. In \cite{NRSc08}, the dynamical
Brownian web and the Brownian net are constructed in the same framework, and
questions of exceptional times (of the former) are investigated. A discrete
space-time version of the dynamical Brownian web was studied in \cite{FNRS07}.

A second model closely related to the Brownian net is a class of stochastic
flows of kernels introduced by Howitt and Warren \cite{HW06}.  These
stochastic flows are families of random transition kernels, describing a
Brownian motion evolving in a random space-time environment. It turns out that
these stochastic flows can be constructed through a random switching between
outgoing paths in a `reference' Brownian web, which plays the role of the
random environment, similar to the construction of the dynamical Brownian web
and the Brownian net in \cite{NRSc08}. A subclass of the stochastic flows of
kernels of Howitt and Warren turns out to be supported on the Brownian net.
This is the subject of the ongoing work \cite{SSS08}. Results established in
the present paper, as well as \cite{NRSc08}, will provide important tools to
analyze the Howitt-Warren flows.

Finally, there are close connections between the Brownian net and low
temperature scaling limits of one-dimensional stochastic Potts models. In
\cite{NRS09}, these scaling limits will be constructed with the help of a
graphical representation based on a marking of paths in the Brownian
net. Their construction uses in an essential way one of the results in the
present paper (the local finiteness of relevant separation points proved in
Proposition~\ref{P:finrel} below).

In the rest of the introduction, we recall the characterization of the
Brownian web and its dual from \cite{FINR04,FINR06}, the characterization of
the left-right Brownian web and the Brownian net from \cite{SS06}, the
classification of special points of the Brownian web from \cite{TW98, FINR06},
and lastly we formulate our main results on the classification of special
points for the left-right Brownian web and the Brownian net according to the
configuration of paths entering and leaving a point.

\subsection{The Brownian web, left-right Brownian web, and Brownian net}
\label{S:wwn}

\begin{figure}
\begin{center}
\setlength{\unitlength}{.7cm}
\begin{picture}(10,8)(-5,-4)
\linethickness{.4pt}
\qbezier(0,-3)(0,0)(0,3)
\qbezier(-3,0)(0,0)(3,0)
\linethickness{.4pt}
\qbezier(0,-3)(-6,0)(0,3)
\qbezier(0,-3)(-5.3,0)(0,3)
\qbezier(0,-3)(-4.5,0)(0,3)
\qbezier(0,-3)(-3,0)(0,3)
\qbezier(0,-3)(3,0)(0,3)
\qbezier(0,-3)(4.5,0)(0,3)
\qbezier(0,-3)(5.3,0)(0,3)
\qbezier(0,-3)(6,0)(0,3)
\qbezier(-1.9,-1.8)(0,-1.8)(1.9,-1.8)
\qbezier(-2.2,-1.5)(0,-1.5)(2.2,-1.5)
\qbezier(-2.65,-1)(0,-1)(2.65,-1)
\qbezier(-2.65,1)(0,1)(2.65,1)
\qbezier(-2.2,1.5)(0,1.5)(2.2,1.5)
\qbezier(-1.9,1.8)(0,1.8)(1.9,1.8)

\put(0,-3){\circle*{.25}}
\put(0,3){\circle*{.25}}
\put(0,0){\circle*{.25}}
\put(2.2,1.5){\circle*{.25}}
\put(-2.65,-1){\circle*{.25}}

\put(-.6,-3.5){$(\ast,-\infty)$}
\put(-.6,3.5){$(\ast,+\infty)$}
\put(-1.3,-.5){$(0,0)$}
\put(2.3,1.6){$(+\infty,2)$}
\put(-5.3,-1.1){$(-\infty,-1)$}

\end{picture}
\caption[The compactification $\Rc$ of $\R^2$.]{The compactification
$\Rc$ of $\R^2$.}\label{fig:comp}
\end{center}
\end{figure}
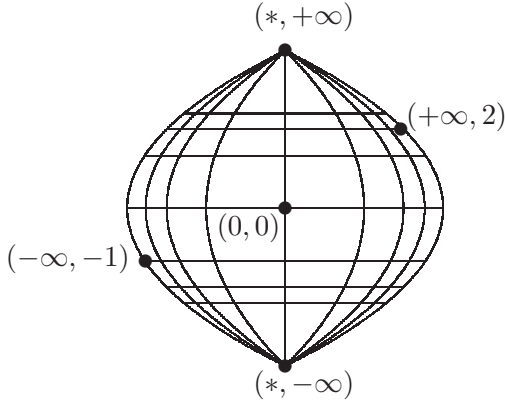

Let us first recall from \cite{FINR04} the space of {\em compact sets of
  paths} in which the Brownian web and the Brownian net take their values. Let
$\Rc$ denote the completion of the space-time plane $\R^2$ w.r.t.\ the metric
\be\label{rho}
\rho\big((x_1, t_1), (x_2,t_2)\big) = \left|\tanh(t_1)-\tanh(t_2)\right|
\ \vee\ \left|\frac{\tanh(x_1)}{1+|t_1|}-\frac{\tanh(x_2)}{1+|t_2|}\right|.
\ee
As a topological space, $\Rc$ can be identified with the continuous image of
$[-\infty, \infty]^2$ under a map that identifies the line
$[-\infty,\infty]\times\{\infty\}$ with a single point $(*, \infty)$, and the
line $[-\infty,\infty]\times\{-\infty\}$ with the point $(*,-\infty)$, see
Figure~\ref{fig:comp}.

A path $\pi$ in $\Rc$, whose starting time we denote by
$\sigma_\pi\in [-\infty,\infty]$, is a mapping $\pi :
[\sigma_\pi,\infty] \to [-\infty, \infty] \cup\{*\}$ such that
$\pi(\infty)=*$, $\pi(\sigma_\pi)=*$ if $\sigma_\pi=-\infty$, and $t
\to (\pi(t), t)$ is a continuous map from $[\sigma_\pi,\infty]$ to
$(\Rc, \rho)$. We then define $\Pi$ to be the space of all
paths in $\Rc$ with all possible starting times in $[-\infty,\infty]$.
Endowed with the metric
\be\label{PId}
\!
d(\pi_1, \pi_2) = \Big|\!\tanh(\sigma_{\pi_1})-\tanh(\sigma_{\pi_2})\Big|
\ \vee \sup_{t\geq \sigma_{\pi_1} \wedge \sigma_{\pi_2}}\!
\left|\frac{\tanh(\pi_1(t\vee \sigma_{\pi_1}))}{1+|t|}
-\frac{\tanh(\pi_2(t\vee \sigma_{\pi_2}))}{1+|t|}\right|, \!\!\!
\ee
$(\Pi, d)$ is a complete separable metric space. Note that
convergence in the metric $d$ can be desrcibed as locally uniform convergence
of paths plus convergence of starting times. (The metric
$d$ differs slightly from the original choice in \cite{FINR04}, which
is slightly less natural as explained in the appendix of \cite{SS06}.)

We can now define $\Hi$, the space of compact subsets of $(\Pi, d)$,
equipped with the Hausdorff metric
\be\label{dH}
d_{\Hi}(K_1, K_2) = \sup_{\pi_1\in K_1} \inf_{\pi_2\in K_2}\!\!
d(\pi_1, \pi_2)\ \vee  \sup_{\pi_2\in K_2}
\inf_{\pi_1\in K_1} d(\pi_1, \pi_2).
\ee
The space $(\Hi, d_\Hi)$ is also a complete separable metric
space. Let $\Bi_\Hi$ be the Borel sigma-algebra associated with
$d_\Hi$. The Brownian web $\Wi$ and the Brownian net $\Ni$ will be
$(\Hi, \Bi_\Hi)$-valued random variables.

In the rest of the paper, for $K\in\Hi$ and $A\subset\Rc$, let
$K(A)$ denote the set of paths in $K$ with starting points in
$A$. When $A=\{z\}$ for $z\in\Rc$, we also write $K(z)$ instead
of $K(\{z\})$.

We recall from \cite{FINR04} the following characterization of the
Brownian web.
\bt\label{T:webchar}{\bf[Characterization of the Brownian web]}\\
There exists a $(\Hi, \Bi_\Hi)$-valued random variable $\Wi$, called
the standard Brownian web, whose distribution is uniquely determined
by the following properties:
\begin{itemize}
\item[{\rm (a)}] For each deterministic $z\in\R^2$, almost surely
there is a unique path $\pi_z\in \Wi(z)$.

\item[{\rm (b)}] For any finite deterministic set of points $z_1, \ldots, z_k
  \in\R^2$, the collection $(\pi_{z_1}, \ldots, \pi_{z_k})$
  is distributed as coalescing Brownian motions.

\item[{\rm (c)}] For any deterministic countable dense subset
$\Di \subset\R^2$, almost surely, $\Wi$ is the closure of
$\{\pi_z : z\in\Di\}$ in $(\Pi, d)$.
\end{itemize}
\et

To each Brownian web $\Wi$, there is associated a {\em dual Brownian web}
$\hat\Wi$, which is a random set of paths running backward in time
\cite{A81,TW98,FINR06}. The pair $(\Wi,\hat\Wi)$ is called the {\em double
  Brownian web}. By definition, a {\em backward path} $\hat\pi$, with starting
time denoted by $\hat\sig_{\hat\pi}$, is a function
$\hat\pi:[-\infty,\hat\sig_{\hat\pi}]\to[-\infty,\infty]\cup\{\ast\}$, such
that $t\mapsto(t,\hat\pi(t))$ is a continuous map from
$[-\infty,\hat\sig_{\hat\pi}]$ to $\Rc$. We let $(\hat\Pi,\hat d)$ denote the
space of backward paths, which is in a natural way isomorphic to $(\Pi,d)$
through time reversal, and we let $(\hat\Hi,d_{\hat\Hi})$ denote the space of
compact subsets of $(\hat\Pi,\hat d)$, equipped with the Hausdorff metric. We
say that a dual path $\hat\pi$ {\em crosses} a (forward) path $\pi$ if there
exist $\sigma_\pi\leq s<t\leq\hat\sigma_{\hat\pi}$ such that
$(\pi(s)-\hat\pi(s))(\pi(t)-\hat\pi(t))<0$. The next theorem follows from
\cite[Theorem 3.7]{FINR06}; a slightly different construction of the dual
Brownian web can be found in \cite[Theorem 1.9]{SS06}.

\bt\label{T:dwebchar}
{\bf [Characterization of the dual Brownian web]}\\
Let $\Wi$ be the standard Brownian web. Then there exists a $\hat\Hi$-valued
random variable $\hat\Wi$, defined on the same probability space as $\Wi$,
called the dual Brownian web, which is almost surely uniquely determined by
the following properties:
\begin{itemize}
\item[{\rm (a)}] For any deterministic $z\in\R^2$, almost surely $\hat\Wi(z)$
  consists of a single path $\hat\pi_z$, which is the unique path in
  $\hat\Pi(z)$ that does not cross any path in $\Wi$.
\item[{\rm (b)}] For any deterministic countable dense subset $\Di
  \subset\R^2$, almost surely, $\hat\Wi$ is the closure of $\{\hat\pi_z :
  z\in\Di\}$ in $(\hat\Pi,\hat d)$.
\end{itemize}
\et
It is known that, modulo a time reversal, $\hat\Wi$ is equally distributed
with $\Wi$. Moreover, paths in $\Wi$ and $\hat\Wi$ interact via Skorohod
reflection. (This follows from the results in \cite{STW00}, together with the
standard discrete aproximation of the double Brownian web.)

We now recall the {\em left-right Brownian web} $(\Wl,\Wr)$, which is
the key intermediate object in the construction of the Brownian net in
\cite{SS06}. Following \cite{SS06}, we call $(l_1, \ldots, l_m; r_1,
\ldots, r_n)$ a collection of {\em left-right coalescing Brownian
motions}, if $(l_1, \ldots, l_m)$ is distributed as coalescing
Brownian motions each with drift $-1$, $(r_1, \ldots, r_n)$ is
distributed as coalescing Brownian motions each with drift $+1$, paths
in $(l_1, \ldots, l_m; r_1, \ldots, r_n)$ evolve independently when
they are apart, and the interaction between $l_i$ and $r_j$ when they
meet is described by the two-dimensional stochastic
  differential equation
\bc\label{lrsde}
\dis\di L_t\ &=&\dis 1_{\{L_t\neq R_t\}}\di B^{\rm l}_t
+1_{\{L_t=R_t\}}\di B^{\rm s}_t-\di t,\\[5pt]
\dis\di R_t\ &=&\dis 1_{\{L_t\neq R_t\}}\di B^{\rm r}_t
+1_{\{L_t=R_t\}}\di B^{\rm s}_t+\di t,
\ec
where $B^{\rm l}_t, B^{\rm r}_t, B^{\rm s}_t$ are independent standard
Brownian motions, and $(L,R)$ are subject to the constraint that
\be\label{lrorder}
L_t\leq R_t\mbox{ for all }t\geq\tau_{L,R},
\ee
where, for any two paths $\pi,\pi'\in\Pi$, we let
\be\label{taudef}
\tau_{\pi,\pi'}:=\inf\{t>\sig_\pi\vee\sig_{\pi'}:\pi(t)=\pi'(t)\}
\ee
denote the first meeting time of $\pi$ and $\pi'$, which may be $\infty$. It
can be shown that subject to the condition (\ref{lrorder}), solutions to the
SDE (\ref{lrsde}) are unique in distribution \cite[Proposition 2.1]{SS06}.
(See also \cite[Proposition 14]{HW07} for a martingale problem
characterization of sticky Brownian motions with drift.)  The interaction
between left-most and right-most paths is a form of sticky reflection; in
particular, if $\{t:L_t=R_t\}$ is nonempty then it is a nowhere dense set
with positive Lebesgue measure \cite[Proposition 3.1]{SS06}.

We cite the following characterization of the left-right
Brownian web from \cite[Theorem~1.5]{SS06}.
\bt\label{T:lrwebchar}{\bf[Characterization of the left-right Brownian web]}\\
There exists a $(\Hi^2, \Bi_{\Hi^2})$-valued random variable $(\Wl,
\Wr)$, called the standard left-right Brownian web, whose distribution
is uniquely determined by the following properties:

\begin{itemize}
\item[{\rm (a)}] For each deterministic $z\in\R^2$, almost surely there are
  unique paths $l_z\in\Wl$ and $r_z\in\Wr$.
\item[{\rm (b)}] For any finite deterministic set of points $z_1,\ldots,z_m,
  z'_1,\ldots,z'_n\in\R^2$, the collection
  $(l_{z_1},\ldots,l_{z_n};r_{z'_1},\ldots,r_{z'_n})$ is distributed as
  left-right coalescing Brownian motions.
\item[{\rm (c)}] For any deterministic countable dense subset $\Di\sub\R^2$,
  almost surely $\Wl$ is the closure of $\{l_z:z\in\Di\}$ and $\Wr$ is the
  closure of $\{r_z:z\in\Di\}$ in the space $(\Pi,d)$.
\end{itemize}
\et
Comparing Theorems~\ref{T:webchar} and \ref{T:lrwebchar}, we see that $\Wl$
and $\Wr$ are distributed as Brownian webs tilted with drift $-1$ and $+1$,
respectively. Therefore, by Theorem~\ref{T:dwebchar}, the Brownian webs $\Wl$
and $\Wr$ a.s.\ uniquely determine dual webs $\hat\Wl$ and $\hat\Wr$,
respectively. It turns out that $(\hat\Wl,\hat\Wr)$ is equally distributed
with $(\Wl, \Wr)$ modulo a rotation by $180^{\rm o}$.

Based on the left-right Brownian web, \cite{SS06} gave three
equivalent characterizations of the Brownian net, which are called
respectively the {\em hopping, wedge, and mesh characterizations}. We
first recall what is meant by hopping, and what are wedges and
meshes.\medskip

\noi
{\bf Hopping:} Given two paths $\pi_1, \pi_2\in\Pi$, any $t >
\sigma_{\pi_1} \vee \sigma_{\pi_2}$ (note the strict inequality) is called
an {\em intersection time} of $\pi_1$ and $\pi_2$ if $\pi_1(t)=\pi_2(t)$. By
hopping from $\pi_1$ to $\pi_2$, we mean the construction of a new
path by concatenating together the piece of $\pi_1$ before and the
piece of $\pi_2$ after an intersection time. Given the left-right
Brownian web $(\Wl, \Wr)$, let $H(\Wl\cup\Wr)$ denote the set of paths
constructed by hopping a finite number of times among paths in
$\Wl\cup\Wr$.\med

\begin{figure}[tp] 
\centering
\plaat{5cm}{\includegraphics[width=10cm]{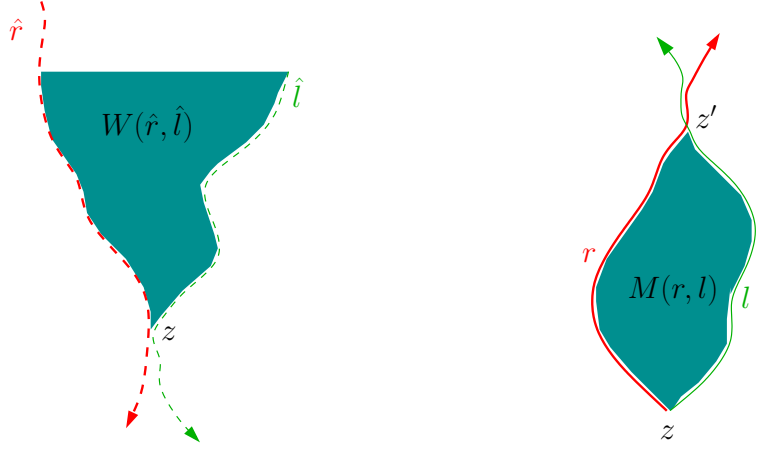}}
\caption{A wedge $W(\hat r,\hat l)$ with
bottom point $z$ and a mesh $M(r,l)$ with bottom and top points $z$ and $z'$.}
\label{fig:meshwedge}
\end{figure}

\noi
{\bf Wedges:} Let $(\hat\Wl, \hat\Wr)$ be the dual left-right Brownian web
almost surely determined by $(\Wl,\Wr)$. Recall that $\hat\sigma_{\hat \pi}$
denotes the starting time of a backward path $\hat\pi$. Any pair $\hat
l\in\hat\Wl$, $\hat r\in \hat\Wr$ with $\hat r(\hat\sigma_{\hat l}\wedge
\hat\sigma_{\hat r}) < \hat l(\hat\sigma_{\hat l}\wedge \hat\sigma_{\hat r})$
defines an open set (see Figure~\ref{fig:meshwedge})
\be\label{wedge}
W(\hat r, \hat l) = \{(x,u)\in\R^2 :\hat\tau_{\hat r,\hat l}<u<
\hat\sigma_{\hat l}\wedge \hat\sigma_{\hat r},\ \hat r(u) < x < \hat l(u) \},
\ee
where, in analogy with (\ref{taudef}), $\hat\tau_{\hat r,\hat
  l}:=\sup\{t<\hat\sigma_{\hat l}\wedge \hat\sigma_{\hat r}: \hat r(t) =\hat
l(t)\}$ denotes the first (backward) hitting time of $\hat r$ and $\hat
l$. Such an open set is called a {\em wedge} of $(\hat\Wl, \hat\Wr)$. If
$\hat\tau_{\hat r,\hat l}>-\infty$, we call $\tau_{\hat r,\hat l}$ the bottom
time, and $(\hat l(\hat\tau_{\hat r,\hat l}),\tau_{\hat r,\hat l})$ the bottom
point of the wedge $W(\hat r,\hat l)$.
\vspace{10pt}

\noindent
{\bf Meshes:} By definition, a {\em mesh} of $(\Wl,\Wr)$ is an open set of
the form (see Figure~\ref{fig:meshwedge})
\be\label{mesh}
M=M(r,l)=\{(x,t)\in\R^2:\sigma_l<t<\tau_{l,r},\ r(t)<x<l(t)\},
\ee
where $l\in\Wl$, $r\in\Wr$ are paths such that $\sigma_l=\sigma_r$,
$l(\sig_l)=r(\sig_r)$ and $r(s)<l(s)$ on $(\sigma_l,\sigma_l+\epsilon)$ for
some $\epsilon>0$.  We call $(l(\sigma_l),\sigma_l)$ the bottom point,
$\sigma_l$ the bottom time, $(l(\tau_{l,r}),\tau_{l,r})$ the top point,
$\tau_{l,r}$ the top time, $r$ the left (!) boundary, and
$l$ the right boundary of $M$.\med

Given an open set $A\subset \R^2$ and a path $\pi\in\Pi$, we say $\pi$
{\em enters} $A$ if there exist $\sigma_\pi <s<t$ such that $\pi(s)
\notin A$ and $\pi(t)\in A$. We say $\pi$ {\em enters $A$ from
outside} if there exists $\sigma_\pi <s<t$ such that $\pi(s) \notin
\bar A$ and $\pi(t)\in A$. We now recall the following
characterization of the Brownian net from \cite[Theorems 1.3, 1.7, 1.10]{SS06}.

\bt\label{T:netchar}{\bf[Characterization of the Brownian net]}\\
There exists a $(\Hi, \Bi_\Hi)$-valued random variable $\Ni$, the
standard Brownian net, whose distribution is uniquely determined by
property {\rm (a)} and any of the three equivalent properties {\rm
(b1)--(b3)} below:
\begin{itemize}
\item[{\rm (a)}] There exist $\Wl, \Wr\subset \Ni$ such that $(\Wl, \Wr)$ is
distributed as the left-right Brownian web.
\item[{\rm (b1)}] Almost surely, $\Ni$ is the closure of $H(\Wl \cup\Wr)$.
\item[{\rm (b2)}] Almost surely, $\Ni$ is the set of paths in $\Pi$ which do
not enter any wedge of $(\hat\Wl, \hat\Wr)$
  from outside.
\item[{\rm (b3)}] Almost surely, $\Ni$ is the set of paths in $\Pi$ which do
not enter any mesh of $(\Wl, \Wr)$.
\end{itemize}
\et
{\bf Remark.} Properties (b1)--(b3) in fact imply that the left-right Brownian
web $(\Wl, \Wr)$ contained in a Brownian net $\Ni$ is
almost surely uniquely determined by the latter, and for
each deterministic $z\in\R^2$, the path in $\Wl$, resp.\ $\Wr$, starting from
$z$ is just the left-most, resp.\ right-most, path among all paths in $\Ni$
starting from $z$. Since $(\Wl, \Wr)$ uniquely determines a dual left-right
Brownian web $(\hat\Wl,\hat\Wr)$, there exists a dual Brownian net $\hat\Ni$
uniquely determined by and equally distributed with $\Ni$ (modulo time
reversal).\med

The construction of the Brownian net from the left-right Brownian web can be
regarded as an outside-in approach because $\Wl$ and $\Wr$ are the
``outermost'' paths among all paths in $\Ni$. On the other hand, the
construction of the Brownian net in \cite{NRSc08} can be regarded as an
inside-out approach, since they start from a standard Brownian web, which may
be viewed as a collection of ``innermost'' paths, to which new paths are added
by allowing branching to the left and right.  More precisely, they allow
hopping at a set of marked points, which is a Poisson subset of the set of all
special points with one incoming and two outgoing paths. One may call this
construction the {\em marking construction} of the Brownian net. In this
paper, we will only use the characterizations provided by
Theorem~\ref{T:netchar}.

\subsection{Classification of special points}\label{S:class}

To classify each point of $\R^2$ according to the local configuration
of paths in the Brownian web or net, we first formulate a notion of
equivalence among paths entering, resp.\ leaving, a point, which
provides a unified framework. We say that a path $\pi\in\Pi$ {\em
enters} a point $z=(x,t)\in\R^2$ if $\sig_\pi<t$ and $\pi(t)=x$. We
say that $\pi$ {\em leaves} $z$ if $\sig_\pi\leq t$ and $\pi(t)=x$.

\bd\label{D:pathequiv}{\bf[Equivalence of paths entering and
leaving a point]}\\
We say $\pi_1,\pi_2\in\Pi$ are {\em equivalent paths entering}
$z=(x,t)\in\R^2$, denoted by $\pi_1\sim^{z}_{\rm in}\pi_2$, if $\pi_1$
and $\pi_2$ enter $z$ and $\pi_1(t-\eps_n)=\pi_2(t-\eps_n)$ for a
sequence $\eps_n\down 0$. We say $\pi_1,\pi_2$ are {\em equivalent
paths leaving} $z$, denoted by $\pi_1\sim^z_{\rm out}\pi_2$, if
$\pi_1$ and $\pi_2$ leave $z$ and $\pi_1(t+\eps_n)=\pi_2(t+\eps_n)$
for a sequence $\eps_n\down 0$.
\ed

Note that, on $\Pi$, $\sim^z_{\rm in}$ and $\sim^z_{\rm out}$ are not
equivalence relations. However, almost surely, they define equivalence
relations on the set of all paths in the Brownian web $\Wi$ entering
or leaving $z$. Due to coalescence, almost surely for all $\pi_1, \pi_2\in \Wi$
and
$z=(x,t)\in\R^2$,
\be\label{wequiv}
\begin{aligned}
\pi_1 \sim^z_{\rm in} \pi_2 \quad &\text{iff}\quad \pi_1=\pi_2
\text{ on } [t-\eps, \infty) \text{ for
  some } \eps>0,  \\
\pi_1 \sim^z_{\rm out} \pi_2 \quad &\text{iff}\quad \pi_1=\pi_2
\text{ on } [t,\infty).
\end{aligned}
\ee
Let $m_{\rm in}(z)$, resp.\ $m_{\rm out}(z)$, denote the number of equivalence
classes of paths in $\Wi$ entering, resp.\ leaving, $z$, and let $\hat
m_{\rm in}(z)$ and $\hat m_{\rm out}(z)$ be defined similarly for the
dual Brownian web $\hat\Wi$. For the Brownian web, points $z\in\R^2$
are classified according to the value of $(m_{\rm in}(z),m_{\rm
out}(z))$. Points of type (1,2) are further divided into
types $(1,2)_{\rm l}$ and $(1,2)_{\rm r}$, where the subscript
$\rm l$ (resp.\ $\rm r$) indicates that the left (resp.\ right) of the
two outgoing paths is the continuation of the (up to equivalence)
unique incoming path. Points in the dual Brownian web $\hat\Wi$ are
labelled according to their type in the Brownian web obtained by
rotating the graph of $\hat\Wi$ in $\R^2$ by $180^{\rm o}$. We cite
the following result from \cite[Proposition~2.4]{TW98} or
\cite[Theorems~3.11--3.14]{FINR06}.

\begin{figure}[tp] 
\centering
\plaat{3.5cm}
{\includegraphics[width=15.6cm]{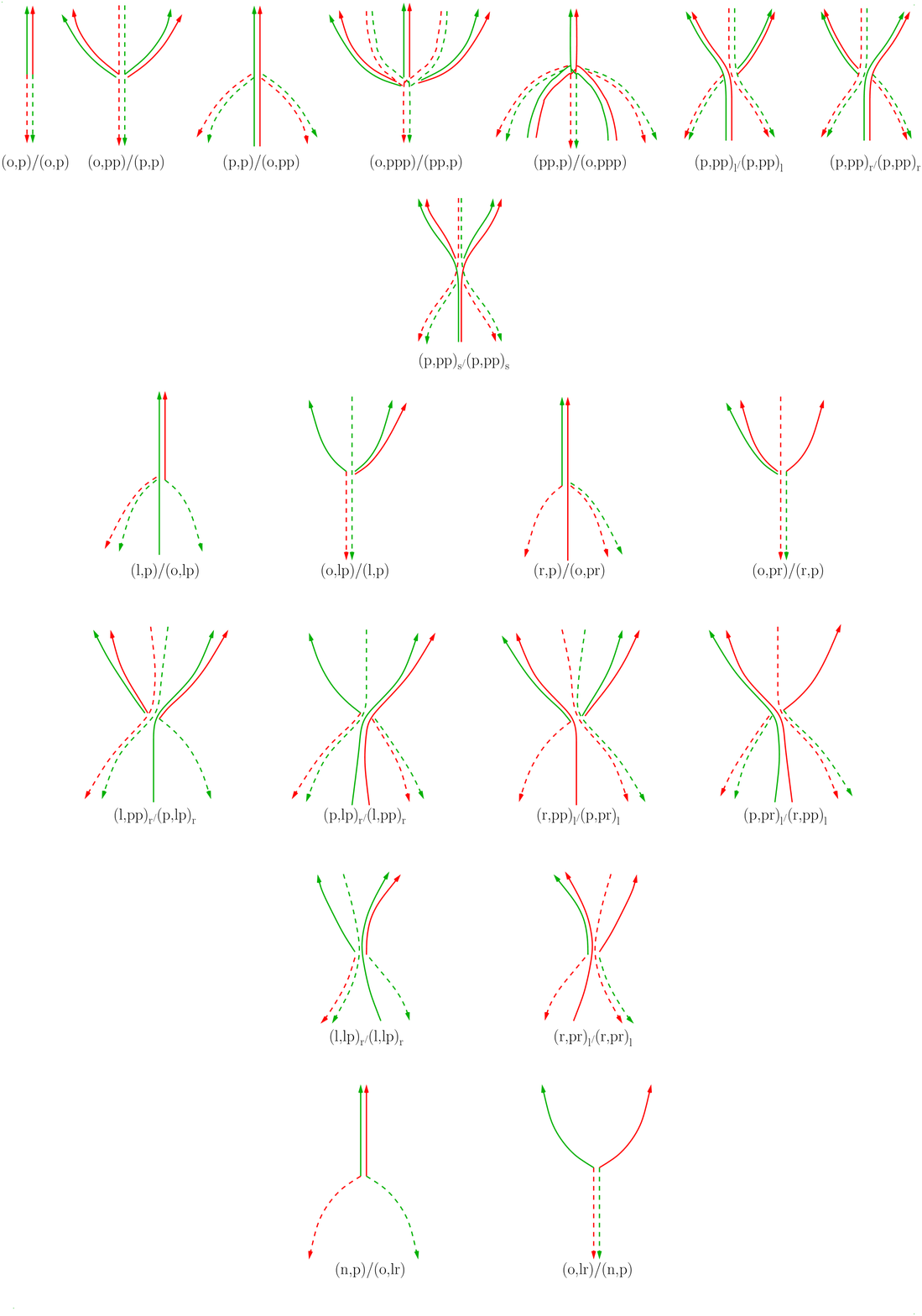}}
\caption{Classification of special points of the Brownian net.}
\label{fig:LRwebmult}
\end{figure}

\bt\label{T:classweb}{\bf[Classification of special points of the
Brownian web]}\\
Let $\Wi$ be the Brownian web and $\hat\Wi$ its dual. Then almost surely, each
$z\in\R^2$ satisfies $m_{\rm out}(z)= \hat m_{\rm in}(z)+1$ and $\hat m_{\rm
  out}(z)=m_{\rm in}(z)+1$, and $z$ is of one of the following seven types in
$\Wi/\hat\Wi$: $(0,1)/(0,1)$, $(0,2)/(1,1)$, $(1,1)/(0,2)$, $(0,3)/(2,1)$,
$(2,1)/(0,3)$, $(1,2)_{\rm l}/(1,2)_{\rm l}$, and $(1,2)_{\rm r}/(1,2)_{\rm
  r}$. For each deterministic $t\in\R$, almost surely, each point in
$\R\times\{t\}$ is of either type $(0,1)/(0,1)$, $(0,2)/(1,1)$, or
$(1,1)/(0,2)$. A deterministic point in $\R^2$ is almost surely of type
$(0,1)/(0,1)$.
\et
We do not give a picture to demonstrate Theorem~\ref{T:classweb}. However, if
in the first row in Figure~\ref{fig:LRwebmult}, one replaces each pair
consisting of one left-most (green) and right-most (red) path by a single
path, and likewise for pairs of dual (dashed) paths, then one obtains a
schematic depiction of the 7 types of points in $\Wi/\hat\Wi$.

We now turn to the problem of classifying the special points of the Brownian
net. We start by observing that also in the left-right Brownian web, a.s.\ for
each $z\in\R^2$, the relations $\sim^z_{\rm in}$ and $\sim^z_{\rm out}$ define
equivalence relations on the set of paths in $\Wl \cup \Wr$ entering,
resp.\ leaving $z$. This follows from (\ref{wequiv}) for the Brownian web
and the fact that a.s.\ for each $l\in\Wl$, $r\in\Wr$ and
$\sig_l\vee\sig_r<s<t$ such that $l(s)=r(s)$, one has $l(t)\leq r(t)$ (see
Prop.~3.6~(a) of \cite{SS06}). Moreover, by the same facts, the equivalence
classes of paths in $\Wl\cup\Wr$ entering, resp.\ leaving, $z$ are naturally
ordered from left to right.

Our classification of points in the Brownian net $\Ni$ is mainly based
on the equivalence classes of incoming and outgoing paths in
$\Wl\cup\Wr$. To denote the type of a point, we first list the
incoming equivalence classes of paths from left to right, and then,
separated by a comma, the outgoing equivalence classes of paths from
left to right. If an equivalence class contains only paths in $\Wl$
resp.\ $\Wr$ we will label it by l, resp.\ r, while if it contains
both paths in $\Wl$ and in $\Wr$ we will label it by p, standing for
pair. For points with (up to equivalence) one incoming and two
outgoing paths, a subscript $\rm l$ resp.\ $\rm r$ means that all
incoming paths belong to the left one, resp.\ right one, of the two
outgoing equivalence classes; a subscript $\rm s$ indicates that
incoming paths in $\Wl$ belong to the left outgoing equivalence class,
while incoming paths in $\Wr$ belong to the right outgoing equivalence
class. If at a point there are no incoming paths in $\Wl\cup\Wr$, then
we denote this by o or n, where o indicates that there are no incoming
paths in the net $\Ni$, while n indicates that there are incoming
paths in~$\Ni$ (but none in $\Wl\cup\Wr$).

Thus, for example, a point is of type $(\rm p,\rm lp)_{\rm r}$ if at
this point there is one equivalence class of incoming paths in
$\Wl\cup\Wr$ and there are two outgoing equivalence classes. The
incoming equivalence class is of type p while the outgoing equivalence
classes are of type l and p, from left to right. All incoming paths in
$\Wl\cup\Wr$ continue as paths in the outgoing equivalence class of type p.

Points in the dual Brownian net $\hat\Ni$, which is defined in terms of the
dual left-right Brownian web $(\hat\Wl,\hat\Wr)$, are labelled according to
their type in the Brownian net and left-right web obtained by rotating the
graphs of $\hat\Ni$ and $(\hat\Wl,\hat\Wr)$ in $\R^2$ by $180^{\rm o}$.  With
the notation introduced above, we can now state our first main result on the
classification of points in $\R^2$ for the Brownian net.

\bt\label{T:classnet}{\bf[Classification of special points of the
Brownian net]}\\
Let $(\Wl,\Wr)$ be the standard left-right Brownian web, let
$(\hat\Wl,\hat\Wr)$ be its dual, and let $\Ni$ and $\hat\Ni$ be the
associated Brownian net and its dual. Then almost surely, each point in $\R^2$
is of one of the following $20$ types in
$\Ni/\hat\Ni$ (see Figure~\ref{fig:LRwebmult}):
\begin{itemize}
\item[{\rm (1)}] {\rm (o,p)/(o,p), (o,pp)/(p,p), (p,p)/(o,pp),
(o,ppp)/(pp,p), (pp,p)/(o,ppp), ${\rm (p,pp)}_{\rm l}/{\rm
(p,pp)}_{\rm l}$, ${\rm (p,pp)}_{\rm r}/{\rm (p,pp)}_{\rm r}$ };
\item[{\rm (2)}] ${\rm (p,pp)}_{\rm s}/{\rm (p,pp)}_{\rm s}$;
\item[{\rm (3)}] {\rm (l,p)/(o,lp), (o,lp)/(l,p), (r,p)/(o,pr), (o,pr)/(r,p)};
\item[{\rm (4)}] ${\rm (l,pp)}_{\rm r}/{\rm (p,lp)}_{\rm r}$,
${\rm (p,lp)}_{\rm r}/{\rm (l,pp)}_{\rm r}$,
${\rm (r,pp)}_{\rm l}/{\rm (p,pr)}_{\rm l}$,
${\rm (p,pr)}_{\rm l}/{\rm (r,pp)}_{\rm l}$;
\item[{\rm (5)}] ${\rm (l,lp)}_{\rm r}/{\rm (l,lp)}_{\rm r}$,
${\rm (r,pr)}_{\rm l}/{\rm (r,pr)}_{\rm l}$;
\item[{\rm (6)}] {\rm (n,p)/(o,lr), (o,lr)/(n,p)};
\end{itemize}
and all of these types occur. For each deterministic time $t\in\R$,
almost surely, each point in $\R\times\{t\}$ is of either type
$\rm(o,p)/(o,p)$, $\rm(o,pp)/(p,p)$, or $\rm(p,p)/(o,pp)$, and all of
these types occur. A deterministic point $(x,t)\in\R^2$ is almost
surely of type $\rm(o,p)/(o,p)$.
\et
{\bf Remark 1} Our classification is mainly based on the configuration
of equivalence classes of paths in the left-right Brownian web
$\Wl\cup\Wr$ entering and leaving a point. For points of types
$\rm(o,p)$ and $\rm(n,p)$, however, we also use information about
paths which belong to the Brownian net $\Ni$ but not to
$\Wl\cup\Wr$. By distinguishing notationally between these types, we
achieve that the type of a point in $\hat\Ni$ is uniquely determined
by its type in $\Ni$. Moreover, by counting $\rm n$ as an incoming
path (and counting equivalence classes of types $\rm l,r,p$ as one
incoming resp.\ outgoing path each), we achieve that $m_{\rm
out}(z)=\hat m_{\rm in}(z)+1$ and $\hat m_{\rm out}(z)=m_{\rm
in}(z)+1$, in analogy with the Brownian web.\med

\noi
{\bf Remark 2} Modulo symmetry between left and right, and between
forward and dual paths, there are only 9 types of points in
Theorem~\ref{T:classnet}: four types from group (1), and one type each
{f}rom groups (2)--(6). Group (1) corresponds to the 7 types of points
of the Brownian web, except that each equivalence class of paths is
now of type p. We call points of type $\rm(pp,p)$ {\em meeting
points}, and points of type $\rm(p,pp)_s$ (from group~(2)) {\em
separation points}. Points in groups (3)--(6) are `cluster points'
that arise as the limit of a nested sequence of excursions between
paths in the left-right Brownian web, or its dual (see
Proposition~\ref{P:imgclas} and Figure~\ref{fig:clusterpts} below).

\subsection{Structure of special points}

When contemplating the special points of the Brownian net as depicted in
Figure~\ref{fig:LRwebmult}, one is struck by the fact that points from
groups~(3)--(6) (the `cluster points') look just like certain points from
group~(1), except that some paths are missing. For example, points of type
$\rm(l,p)/(o,lp)$ look like points of type $\rm(p,p)/(o,pp)$, points of type
$\rm(l,pp)_r/(p,lp)_r$ look like points of $\rm(p,pp)_r/(p,pp)_r$, and points
of type $\rm(n,p)/(o,lp)$ look like points of type $\rm(p,p)/(o,pp)$. In most
cases, when one member of a (dual) left-right pair seems to be missing, the
other member of the pair is still present, but for points of group~(6), a
whole (dual) incoming pair seems to have disappeared. In the present section,
we will show how to make these `missing' paths visible in the form of {\em
  reflected left-most} (or {\em right-most}) {\em paths}. These reflected
left-most (resp.\ right-most) paths are not elements of $\Wl$ (resp.\ $\Wr$),
but they turn out to be countable concatenations of paths in $\Wl$
(resp.\ $\Wr$).

Except for fulfilling the aesthetic role of bringing to light what we feel is
missing, these reflected paths also serve the more practical aim of putting
limitations on how general paths in the Brownian net $\Ni$ (and not just the
left-right Brownian web $\Wl\cup\Wr$) can enter and leave points of various
types. Recall that our classicifation theorem (Theorem~\ref{T:classnet})
primarily describes how left-most and right-most paths enter and leave points
in the plane. In Theorems~\ref{T:locspec} and \ref{T:netpt} below, we first
describe the structure of a special set of reflected left-most and right-most
paths near the special points, and then we use this to describe the local
structure of the Brownian net at these points. Our results will show that,
with the exception of outgoing paths at points of type $\rm(o,lr)$, all
Brownian net paths must enter or leave a point squeezed between a pair
consisting of one (reflected) left-most and (reflected) right-most path.

We start with some definitions. By definition, we say that a set $\Ki$ of
paths, all starting at the same time $t$, has a {\em maximum} (resp.\ {\em
  minimum}) if there exists a path $\pi\in\Ki$ such that $\pi'\leq\pi$
(resp.\ $\pi\leq\pi'$) on $[t,\infty)$ for all $\pi'\in\Ki$. We denote the
(necessarily unique) maximum (resp.\ minimum) of $\Ki$ by ${\rm max}(\Ki)$
(resp.\ ${\rm min}(\Ki)$).
\bl\label{L:extl}{\bf[Reflected left-most and right-most paths]}\\
Almost surely, for each $\hat\pi\in\hat\Ni$ and $z=(x,t)\in\Rc$ such
that $t\leq\hat\sig_{\hat\pi}$ and $\hat\pi(t)\leq x$ (resp.\
$x\leq\hat\pi(t)$), there exists a unique path $l_{z,\hat\pi}\in\Ni$
(resp.\ $r_{z,\hat\pi}\in\Ni$), defined by
\bc\label{extl}
\dis l_{z,\hat\pi}&:=&\dis\min\big\{\pi\in\Ni(z):
\hat\pi\leq\pi\mbox{ on }[t,\hat\sig_{\hat\pi}]\big\},\\[5pt]
\dis r_{z,\hat\pi}&:=&\dis\max\big\{\pi\in\Ni(z):
\pi\leq\hat\pi\mbox{ on }[t,\hat\sig_{\hat\pi}]\big\}.
\ec
The set
\be\label{Ridef}
\Fi(l_{z,\hat\pi}):=\big\{s\in(t,\hat\sig_{\hat\pi}):
\hat\pi(s)=l_{z,\hat\pi}(s),\ (\hat\pi(s),s)
\mbox{ is of type }{\rm(p,pp)_{\rm s}}\big\}
\ee
is a locally finite subset of $(t,\hat\sig_{\hat\pi})$. Let
$\Fi':=\Fi(l_{z,\hat\pi})\cup\{t,\hat\sig_{\hat\pi},\infty\}$ if
$\hat\sig_{\hat\pi}$ is a cluster point of $\Fi(l_{z,\hat\pi})$ and
$\Fi':=\Fi\cup\{t,\infty\}$ otherwise. Then, for each $s,u\in\Fi'$ such that
$s<u$ and $(s,u)\cap\Fi'=\emptyset$, there exists an $l\in\Wl$ such that
$l=l_{z,\hat\pi}$ on $[s,u]$. If $u<\hat\sig_{\hat\pi}$, then $l\leq\hat\pi$
on $[u,\hat\sig_{\hat\pi}]$ and $\inf\{s'>s:l(s')<\hat\pi(s')\}=u$. Analogous
statements hold for $r_{z,\hat\pi}$.
\el

We call the path $l_{z,\hat\pi}$ (resp.\ $r_{z,\hat\pi}$) in (\ref{extl}) the
{\em reflected left-most} (resp.\ {\em right-most}) {\em path relative to
  $\hat\pi$}. See Figure~\ref{fig:seqcross} below for a picture of a reflected
right-most path relative to a dual left-most path.  Lemma~\ref{L:extl} says
that such reflected left-most (resp.\ right-most) paths are well-defined, and
are concatenations of countably many left-most (resp.\ right-most)
paths. Indeed, it can be shown that in a certain well-defined way, a reflected
left-most path $l_{z,\hat\pi}$ always `turns left' at separation points,
except those in the countable set $\Fi(l_{z,\hat\pi})$, where turning left
would make it cross $\hat\pi$.  We call the set $\Fi(l_{z,\hat\pi})$ in
(\ref{Ridef}) the set of {\em reflection times} of $l_{z,\hat\pi}$, and define
$\Fi(r_{z,\hat\pi})$ analogously. In analogy with (\ref{extl}), we also define
{\em reflected dual paths} $\hat l_{z,\pi}$ and $\hat r_{z,\pi}$ relative to a
forward path $\pi\in\Ni$.

For a given point $z$ in the plane, we now define special
classes of reflected left-most and right-most paths, which extend the classes
of paths in $\Wl$ and $\Wr$ entering and leaving $z$.
\bd\label{D:extpath}{\bf[Extended left-most and right-most paths]}\\
For each $z=(x,t)\in\R^2$, we define
\be
\Wl_{\rm in}(z):=\{l\in\Wl:l\mbox{ enters }z\}\quad
\mbox{and}\quad\Wl_{\rm out}(z):=\{l\in\Wl:l\mbox{ leaves }z\}.
\ee
Similar definitions apply to $\Wr,\hat\Wl,\hat\Wr,\Ni$, and $\hat\Ni$.
We define
\be
\El_{\rm in}(z):=\big\{l_{z',\hat r}\in\Ni_{\rm in}(z):
\hat r\in\hat\Wr_{\rm out}(z),\ z'=(x',t'),\ t'<t,\ \hat r(t')\leq x'\big\},
\ee
and define $\Er_{\rm in}(z),\hat\El_{\rm in}(z),\hat\Er_{\rm in}(z)$
analogously, by symmetry. Finally, we define
\be\label{El}
\El_{\rm out}(z):=\big\{l_{z,\hat r}\in\Ni_{\rm out}(z):
\hat r\in\hat\Er_{\rm in}(z'),\ z'=(x',t),\ x'\leq x\big\},
\ee
and we define $\Er_{\rm out}(z),\hat\El_{\rm out}(z),\hat\Er_{\rm out}(z)$
analogously. We call the elements of $\El_{\rm in}(z)$ and $\El_{\rm out}(z)$
{\em extended left-most paths} entering, resp.\ leaving $z$.
See Figure \ref{fig:locref} for an illustration.
\ed

As we will see in Theorem~\ref{T:locspec} below, the extended left- and
right-most paths we have just defined are exactly the `missing' paths in
groups~(3)--(6) from Theorem~\ref{T:classnet}. Finding the right definition of
extended paths that serves this aim turned out to be rather subtle. Note that
paths in $\El_{\rm in}(z)$ and $\Er_{\rm in}(z)$ reflect off paths in
$\hat\Wr_{\rm out}(z)$ and $\hat\Wl_{\rm out}(z)$, respectively, but paths in
$\El_{\rm out}(z)$ and $\Er_{\rm out}(z)$ could be reflected off reflected
paths. If, instead of (\ref{El}), we would have defined $\El_{\rm out}(z)$ as
the set of all reflected left-most paths that reflect off paths in
$\hat\Wr_{\rm in}(z)$ or $\hat\Ni_{\rm in}(z)$, then one can check that for
points of type $\rm(o,lr)$, in the first case we would not have found all
`missing' paths, while in the second case we would obtain too many paths.
Likewise, just calling any countable concatenation of left-most paths an
extended left-most path would, in view of our aims, yield too many paths.

To formulate our results rigorously, we need one more definition.
\bd\label{D:strongequiv}{\bf[Strong equivalence of paths]}\\
We say that two paths $\pi_1,\pi_2\in\Pi$ entering a point $z=(x,t)\in\R^2$
are {\em strongly equivalent}, denoted by $\pi_1=^{z}_{\rm in}\pi_2$, if
$\pi_1=\pi_2$ on $[t-\eps,t]$ for some $\eps>0$. We say that two paths
$\pi_1,\pi_2\in\Pi$ leaving $z$ are {\em strongly equivalent}, denoted by
$\pi_1=^{z}_{\rm out}\pi_2$, if $\pi_1=\pi_2$ on $[t,t+\eps]$ for some
$\eps>0$. We say that two classes $\Ci,\ti\Ci$ of paths entering
(resp.\ leaving) $z$ are equal {\em up to strong equivalence} if for each
$\pi\in\Ci$ there exists a $\ti\pi\in\ti\Ci$ such that $\pi=^z_{\rm in}\ti\pi$
(resp.\ $\pi=^z_{\rm out}\ti\pi$), and vice versa.
\ed
Note that by (\ref{wequiv}), a.s.\ for each $l_1,l_2\in\Wl$ and $z\in\R^2$,
$l_1\sim^z_{\rm in}l_2$ implies $l_1=^z_{\rm in}l_2$ and $l_1\sim^z_{\rm
  out}l_2$ implies $l_1=^z_{\rm out}l_2$. Part~(a) of the next theorem shows
that our extended left-most and right-most paths have the same property.
Parts~(b) and (c) say that extended left- and right-most always pair up, and
are in fact the `missing' paths from groups~(3)--(6) of
Theorem~\ref{T:classnet}.

\bt\label{T:locspec}{\bf[Extended left-most and right-most paths]}\\
Almost surely, for each $z=(x,t)\in\R^2$:
\begin{itemize}
\item[\rm(a)] If $l_1,l_2\in\El_{\rm in}(z)$ satisfy $l_1\sim^z_{\rm in}l_2$,
  then $l_1=^z_{\rm in}l_2$. If $l_1,l_2\in\El_{\rm out}(z)$ satisfy
  $l_1\sim^z_{\rm out}l_2$, then $l_1=^z_{\rm out}l_2$.

\item[\rm(b)] If $z$ is of type $\rm(\cdot,p)/(o,\cdot)$,
  $\rm(\cdot,pp)/(p,\cdot)$, or $\rm(\cdot,ppp)/(pp,\cdot)$ in $\Ni/\hat\Ni$,
  then, up to strong equivalence, one has $\El_{\rm out}(z)=\Wl_{\rm out}(z)$,
  $\Er_{\rm out}(z)=\Wr_{\rm out}(z)$, $\hat\El_{\rm in}(z)=\hat\Wl_{\rm
    in}(z)$, and $\hat\Er_{\rm in}(z)=\hat\Wr_{\rm in}(z)$.

\item[\rm(c)] If $z$ is of type $\rm(\cdot,pp)/(p,\cdot)$,
  $\rm(\cdot,lp)/(l,\cdot)$, $\rm(\cdot,pr)/(r,\cdot)$, or
  $\rm(\cdot,lr)/(n,\cdot)$ in $\Ni/\hat\Ni$, then up to strong equivalence,
  one has $\El_{\rm out}(z)=\{l,l'\}$, $\Er_{\rm out}(z)=\{r,r'\}$,
  $\hat\El_{\rm in}(z)=\{\hat l\}$, and $\hat\Er_{\rm in}(z)=\{\hat r\}$,
  where $l$ is the left-most element of $\Wl(z)$, $r$ is the right-most
  element of $\Wr(z)$, and $l',r',\hat l,\hat r$ are paths satisfying
  $l\sim^z_{\rm out}r'$, $\hat r\sim^z_{\rm in}\hat l$, and $l'\sim^z_{\rm
    out}r$. (See Figure~\ref{fig:locref}.)
\end{itemize}
\et
Note that parts~(b) and (c) cover all types of points from
Theorem~\ref{T:classnet}. Indeed, if we are only interested in the
structure of $\Ni_{\rm out}(z)$ and $\hat\Ni_{\rm in}(z)$, then a.s.\ every
point $z\in\R^2$ is of type $\rm(\cdot,p)/(o,\cdot)$, $\rm(\cdot,pp)/(p,\cdot)$,
$\rm(\cdot,ppp)/(pp,\cdot)$, $\rm(\cdot,lp)/(l,\cdot)$,
$\rm(\cdot,pr)/(r,\cdot)$, or $\rm(\cdot,lr)/(n,\cdot)$.

\begin{figure}[tp] 
\centering
\plaat{5cm}
{\includegraphics[width=8.5cm]{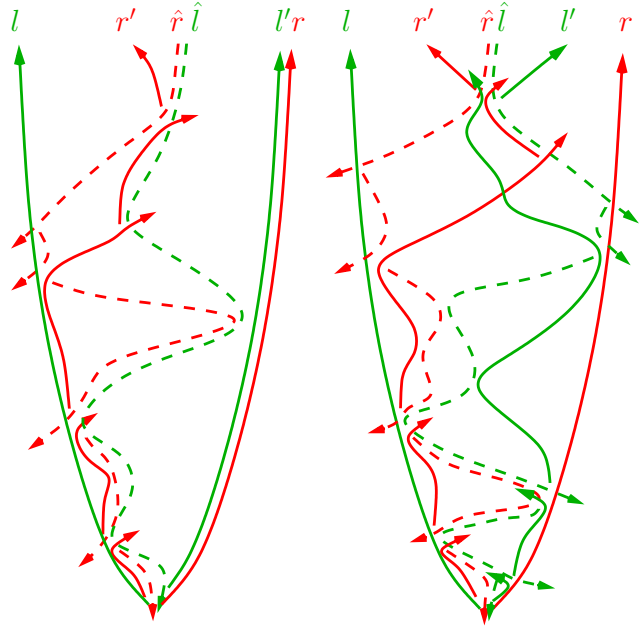}}
\caption{Local structure of outgoing extended left-most and right-most paths
  and incoming extended dual left-most and right-most paths in points of types
  ${\rm(\,\cdot\,,lp)}$ and ${\rm(o,lr)}$. In the picture on the left, $r'$
  and $\hat r$ are reflected paths that are not elements of $\Wr$ and
  $\hat\Wr$, respectively. In the picture on the right, only $l$ and $r$ are
  `true' left-most and right-most paths, while all other paths are `missing'
  paths that are not visible in the schematic picture for points of type
  ${\rm(o,lr)}$ in Figure~\ref{fig:LRwebmult}.}
\label{fig:locref}
\end{figure}

We finally turn our attention to the way general paths in $\Ni$ (and not just
our extended paths) enter and leave the special points from
Theorem~\ref{T:classnet}. The next theorem shows that with the exception of
outgoing paths at points of type $\rm(o,lr)$, all Brownian net paths must
enter and leave points squeezed between a pair consisting of one extended
left-most and one extended right-most path.

\bt\label{T:netpt}{\bf[Structure of special points]}\\
Almost surely, for each $z=(x,t)\in\R^2$:
\begin{itemize}
\item[\rm(a)] The relation $\sim^z_{\rm in}$ is an equivalence relation on
  $\Ni_{\rm in}(z)$. Each equivalence class $\Ci$ of paths in $\Ni_{\rm
  in}(z)$ contains an $l\in\El_{\rm in}(z)$ and $r\in\Er_{\rm in}(z)$, which
  are unique up to strong equivalence, and each $\pi\in\Ci$ satisfies
  $l\leq\pi\leq r$ on $[t-\eps,t]$ for some $\eps>0$.

\item[\rm(b)] If $z$ is not of type ${\rm(o,lr)}$, then the relation
  $\sim^z_{\rm out}$ is an equivalence relation on $\Ni_{\rm out}(z)$. Each
  equivalence class $\Ci$ of paths in $\Ni_{\rm out}(z)$ contains an
  $l\in\El_{\rm out}(z)$ and $r\in\Er_{\rm out}(z)$, which are unique up to
  strong equivalence, and each $\pi\in\Ci$ satisfies $l\leq\pi\leq r$ on
  $[t,t+\eps]$ for some $\eps>0$.

\item[\rm(c)] If $z$ is of type ${\rm(o,lr)}$, then there exist $\pi\in\Ni$
  such that $l\sim^z_{\rm out}\pi\sim^z_{\rm out}r$ while $l\not\sim^z_{\rm
    out}r$, where $l$ and $r$ are the unique outgoing paths at $z$ in $\Wl$
  and $\Wr$, respectively.

\item[\rm(d)] At points of types with the subscript ${\rm l}$ (resp.\ ${\rm
  r}$), all incoming paths in $\Ni$ continue in the left (resp.\ right)
  outgoing equivalence class. Except for this restriction, any concatenation
  of a path in $\Ni_{\rm in}(z)$ up to time $t$ with a path in $\Ni_{\rm
    out}(z)$ after time $t$ is again a path in $\Ni$.
\end{itemize}
\et

\subsection{Outline and open problems}\label{S:outline}

In this section, we outline the main structure of our proofs and
mention some open problems.

In Section~\ref{S:inter} we study separation points, i.e., points of
type $\rm(p,pp)_s$ from Theorem~\ref{T:classnet}. In a sense, these
are the most important points in the Brownian net, since at these
points paths in the Brownian net have a choice whether to turn left or
right. Also, these are exactly the marked points in \cite{NRSc08}, and
their marking construction shows that the Brownian net is a.s.\
determined by its set of separation points and an embedded Brownian
web.

In order to prepare for our study of separation points, in
Section~\ref{S:inter1}, we investigate the interaction between forward
right-most and dual left-most paths. It turns out that the former are Skorohod
reflected off the latter, albeit they may cross the latter from left to right
at some random time. In Section~\ref{S:inter3}, the results from
Section~\ref{S:inter1} are used to prove that crossing points between forward
right-most and dual left-most paths are separation points between left-most
and right-most paths, and that these points are of type $\rm(p,pp)_s$.

In Section~\ref{S:inter4}, we study `relevant' separation points. By
definition, a point $z=(x,t)\in\R^2$ is called an {\em $(s,u)$-relevant
  separation point} for some $-\infty\leq s<u\leq\infty$, if $s<t<u$, there
exists a $\pi\in\Ni$ such that $\sig_\pi=s$ and $\pi(t)=x$, and there exist
incoming $l\in\Wl$ and $r\in\Wr$ at $z$ such that $l<r$ on $(t,u)$. At
$(s,u)$-relevant separation points, a Brownian net path going from time $s$ to
time $u$ has a choice whether to turn left or right that may be relevant for
where it ends up at time $u$. The main result of Section~\ref{S:inter4} says
that for deterministic $-\infty<s<t<\infty$, the set of $(s,u)$-relevant
separation points is a.s.\ a locally finite subset of $\R\times(s,u)$. This
fact has several useful consequences. As a first application, in
Section~\ref{S:inter5}, we prove Lemma~\ref{L:extl}.

In Section~\ref{S:image}, we study the image set of
the Brownian net started at a fixed time $T$. Let
$\Ni_T:=\{\pi\in\Ni:\sig_\pi=T\}$ denote the set of paths in the
Brownian net starting at a given time $T\in[-\infty,\infty]$, and let
$N_T$ be the {\em image set} of $\Ni_T$ in $\Rc$, i.e.,
\be\label{NT}
N_T:=\big\{(\pi(t),t):\pi\in\Ni_T,\ t\in[T,\infty]\big\}.
\ee
By \cite[Prop.~1.13]{SS06}, a.s.\ for each $T\in[-\infty,\infty]$,
\be\label{improp}
\Ni_T=\big\{\pi\in\Pi:\sig_\pi=T,\ \pi\sub N_T\big\}.
\ee
In view of (\ref{improp}), much can be learned about the Brownian net
by studying the closed set $N_T$.

In Section~\ref{S:image1}, it is shown that the connected components
of the complement of $N_T$ relative to $\{(x,t)\in\Rc:t\geq T\}$ are
meshes of a special sort, called {\it maximal $T$-meshes}. In
Section~\ref{S:image2}, it is shown that $N_T$ has a local
reversibility property that allows one, for example, to deduce
properties of meeting points from properties of separation points.
Using these facts, in Section~\ref{S:image4}, we give a
preliminary classification of points in $\R^2$ based only on the
structure of incoming paths in $\Ni$. In Section~\ref{S:image5}, we
use the fractal structure of $N_T$ to prove the existence, announced
in \cite{SS06}, of random times $t>T$ when $\{x\in\R:(x,t)\in N_T\}$
is not locally finite.

It turns out that to determine the type of a point $z\in\R^2$ in
$\Ni$, according to the classification of Theorem~\ref{T:classnet},
except for one trivial ambiguity, it suffices to know the structure of
the incoming paths at $z$ in both $\Ni$ and $\hat\Ni$ according to the
preliminary classification from Section~\ref{S:image4}. Therefore, in
order to prove Theorem~\ref{T:classnet}, we need to know which
combinations of types in $\Ni$ and $\hat\Ni$ are possible according to
the latter classification. In particular, proving the existence of
points in groups (4) and (5) from Theorem~\ref{T:classnet} depends on
showing that there are points where the incoming paths in both $\Ni$
and $\hat\Ni$ form a nested sequence of
excursions. Section~\ref{S:excur} contains some excursion theoretic
arguments that prepare for this. In particular, in
Section~\ref{S:excur1}, we prove that there are many excursions
between a given left-most path $l$ and dual left-most path $\hat l$ on the
left of $l$ that are entered by some dual right-most path. In
Section~\ref{S:excur2}, we prove that there are many points where
$\hat l$ hits $l$ while at the same time some right-most path makes an
excursion away from $l$.

In Section~\ref{S:struct}, we finally prove our main
results. Section~\ref{S:struct1} contains the proof of
Theorem~\ref{T:classnet}, while Section~\ref{S:struct2} contains the proofs of
Theorems~\ref{T:locspec} and \ref{T:netpt}.

We conclude the paper with two appendices (Appendices~\ref{A:extra1}
and \ref{A:extra2}) containing some facts and proofs that are not used
in the main argument but may be of independent interest.

Our investigations leave open a few questions that we believe are
important for understanding the full structure of the Brownian net,
and which we hope to settle in future work. The first question we
would like to mention concerns the image set in (\ref{NT}). Fix
$-\infty\leq s<u\leq\infty$, and say that a point $z\in\R^2$ is {\em
$n$-connected} to $s$ and $u$ if $z\in\{(x,t)\in N_s:s<t<u\}$ and one
needs to remove at least $n$ points from $N_s$ to disconnect $z$ from
$\{(x,t)\in N_s:t=s\mbox{ or }t=u\}$. Note that the notion of connectedness
in $N_s$ is graph theoretic and does not respect the time direction
inherent in the Brownian net. Here is a conjecture:
\bcon{\bf[3-connected points]}\label{C:three}\\
Almost surely for each $-\infty\leq s<u\leq\infty$, the set of
3-connected points to $s$ and $u$ is a locally finite subset of
$\R\times(s,u)$, and all 3-connected points are either meeting or
separation points.
\econ
It is easy to see that all points in $\{(x,t)\in N_s:s<t<u\}$ are 2-connected
to $s$ and $u$. The results in our present paper imply that meeting and
separation points are not 4-connected. All $(s,u)$-relevant separation points
are 3-connected, but not all 3-connected separation points are
$(s,u)$-relevant. If Conjecture~\ref{C:three} is correct, then the structure
of all Brownian net paths going from time $s$ to $u$ can be described by a
purely 3-connected, locally finite graph, whose vertices are the 3-connected
points.

Other open problems concern the reflected left-most and right-most
paths described in Theorem~\ref{T:locspec}. Here is another conjecture:
\bcon{\bf[Reflection points]}\label{C:refl}\\
Almost surely for all $z=(x,t)\in\R^2$, if $l,l',\hat l,r,r',\hat r$ are as in
Theorem~\ref{T:locspec}~(c), then there exists an $\eps>0$ such that $\hat
r(u)=\hat l(u)$ for all $u\in\Fi(\hat r)\cap[t,t+\eps]$.
\econ
Conjecture~\ref{C:refl} says that near $t$, all reflection
points of $\hat r$ off $l$ lie on $\hat l$ (and hence, by symmetry, a similar
statement holds for the  reflection points of $\hat l$ off $r$).
If this is true, then the picture in Figure~\ref{fig:locref}
simplifies a lot. In particular, there exists an $\eps>0$ such that on
$[t,t+\eps]$, the paths are eventually ordered as $l\leq r'\leq\hat
r\leq\hat l\leq l'\leq r$. Note from Figure~\ref{fig:locref} that at
present, for points of type $\rm(o,lr)$, we cannot even rule out that
$l'(u_n)<r'(u_n)$ for a sequence of times $u_n\down t$.

It seems that Conjectures~\ref{C:three} and \ref{C:refl} cannot be proved with
the methods developed in the present paper and \cite{SS06}. Instead, we hope
to tackle these problems by calling in the marking construction of the
Brownian net developed in \cite{NRSc08}.

Another open problem is to determine the Hausdorff dimension of the sets of
points of each of the types from Theorem~\ref{T:classnet}. For the Brownian
web, the Hausdorff dimensions of all types of points are known, see
\cite[Theorem~3.12]{FINR06}. We believe that points from group~(1) of
Theorem~\ref{T:classnet} have the same Hausdorff dimension as the
corresponding points in the Brownian web. Separation points (group~(2)) are
countable. About the Hausdorff dimensions of points from groups~(3)--(6), we
know nothing.

\subsection{List of global notations}

For ease of reference, we collect here some notations that will
be used throughout the rest of the paper. In the proofs, some new notation
might be derived from the global notations listed here, such as by adding
superscripts or subscripts, in which case the objects they encode will be
closely related to what the corresponding global notation stands for.\med

\begin{tabular}{r@{\ \ }l}
\multicolumn{2}{l}{\bf General notation:}\\[3pt]
$\Li(\,\cdot\,)$\;:
& law of a random variable.\\
$\Di$\;: & a deterministic countable dense subset of $\R^2$.\\
$S_{\eps}$\;: & the diffusive scaling map, applied to subsets of $\R^2$,
paths, and sets of\\
& paths, defined as $S_\eps(x,t):=(\eps x,\eps^2t)$ for $(x,t)\in\R^2$.\\
$\tau$\;: & a stopping time;\\
& in Section~\ref{S:excur}: a time
  of increase of the reflection process (see (\ref{tauincr})).\\
$\tau_{\pi,\pi'}$\;: & the first meeting time of $\pi$ and
  $\pi'$, defined in (\ref{taudef}).\\
$(L_t,R_t)$\;: & a pair of diffusions solving the
left-right SDE (\ref{lrsde}).\\
$\sim^z_{\rm in}, \sim^z_{\rm out}$\;: & equivalence of paths entering,
resp.\ leaving $z\in\R^2$, see Definition \ref{D:pathequiv}.\\
$=^z_{\rm in}, =^z_{\rm out}$\;: & strong equivalence of paths entering,
resp.\ leaving $z\in\R^2$,\\
& see Definition \ref{D:strongequiv}.
\end{tabular}
\med

\begin{tabular}{r@{\ \ }l}
\multicolumn{2}{l}{\bf Paths, Space of Paths:}\\[3pt]
$(\Rc,\rho)$\;: & the compactification of $\R^2$ with the metric $\rho$, see
(\ref{rho}).\\
$z=(x,t)$\;: & point in $\Rc$, with
  position $x$ and time $t$.\\
$s,t,u,S,T,U$\;: & times.\\
$(\Pi,d)$\;: & the space of continuous paths in $(\Rc,\rho)$ with metric $d$,
see (\ref{PId}).\\
$\Pi(A),\Pi(z)$\;:
& the set of paths in $\Pi$ starting from a set $A\sub\Rc$
resp.\ a point $z\in\Rc$.\\
& The same notatation applies to any subset of $\Pi$ such as
$\Wi,\Wl,\Ni$.\\
$\pi$\;: & a path in $\Pi$.\\
$\sigma_\pi$\;: & the starting time of the path $\pi$\\
$\pi(t)$\;: & the position of $\pi$ at time $t\geq\sigma_\pi$.\\
$(\hat\Pi,\hat d)$\;: & the space of continuous backward paths in $(\Rc,\rho)$
  with metric $\hat d$,\\
& see Theorem \ref{T:dwebchar}.\\
$\hat\Pi(A),\hat\Pi(z)$\;:
& the set of backward paths in $\hat\Pi$ starting from
$A\sub\Rc$ resp.\ $z\in\Rc$.\\
$\hat\pi$\;: & a path in $\hat\Pi$.\\
$\hat\sigma_{\hat\pi}$\;: & the starting time of the backward path $\hat\pi$.\\
$\Ii(\pi,\hat\pi)$\;: 
& the set of intersection times of $\pi$ and $\hat\pi$,
see (\ref{Iidef}).\\
$(\Hi, d_\Hi)$\;: & the space of compact subsets of $(\Pi,d)$ with Hausdorff
  metric $d_\Hi$,\\ 
& see (\ref{dH}).
\end{tabular}

\begin{tabular}{r@{\ \ }l}
\multicolumn{2}{l}{\bf Brownian webs:}\\[3pt]
$(\Wi,\hat\Wi)$\;: & a double Brownian web
      consisting of a Brownian web and its dual.\\
$\pi_z,\hat\pi_z$\;: & the a.s.\ unique
paths in $\Wi$ resp.\ $\hat\Wi$ starting from\\
& a deterministic point $z\in\R^2$.\\
$(\Wl,\Wr)$\;: & the left-right Brownian web, see Theorem \ref{T:lrwebchar}.\\
$\Wi_{\rm in}(z),\Wi_{\rm out}(z)$\;: & the set of paths in $\Wi$ entering,
resp.\ leaving $z$.\\
$\El_{\rm in}(z),\El_{\rm out}(z)$\;: &
the sets of extended left-most paths entering,
resp.\ leaving $z$,\\
& see Definition~\ref{D:extpath}.\\
$l,r,\hat l,\hat r$\;: & a path in $\Wl$, resp. $\Wr$, $\hat\Wl$, $\hat\Wr$,
often called a left-most, resp.\\
& right-most, dual left-most, dual right-most path.\\
$l_z,r_z,\hat l_z,\hat r_z$\;: & the a.s.\ unique
paths $\Wl(z)$, resp. $\Wr(z)$, $\hat\Wl(z)$, $\hat\Wr(z)$\\
& starting from a deterministic point $z\in\R^2$.\\
$l_{z,\hat\pi}, r_{z,\hat\pi}$\;: & reflected left-most and right-most paths
starting from $z$,\\
& and reflected off $\hat\pi\in\hat\Ni$, see (\ref{extl}).\\
$W(\hat r, \hat l)$\;: & a wedge defined by the paths $\hat r\in\hat\Wr$ and
$\hat l\in\hat\Wl$, see (\ref{wedge}).\\
$M(r, l)$\;: & a mesh defined by the paths $r\in\Wr$ and $l\in\Wl$, see
(\ref{mesh}).
\end{tabular}

\begin{tabular}{r@{\ \ }l}
\multicolumn{2}{l}{\bf Brownian net:}\\[3pt]
$(\Ni,\hat\Ni)$\;: & the Brownian net and the dual Brownian net,
see Theorem \ref{T:netchar}.\\
$\Ni_T$\;: & the set of paths in $\Ni$ starting at time
$T\in[-\infty,\infty]$.\\
$N_T$\;: & the image set of $\Ni_T$ in $\R^2$, see (\ref{NT}).\\
$\xi^K_t$\;: &
$\{\pi(t):\pi\in\Ni(K),\ \sig_\pi\leq t\}$, the set of
  positions through which there passes\\
& a path in the Brownian net starting from the set
$K\sub\Rc$.\\
$\xi^{(T)}_t$\;: &
$\{\pi(t):\pi\in\Ni,\ \sig_\pi=T\leq t\}$,
the same as $\xi^K_t$, but now for Brownian\\
& net paths starting from a given time
$T\in[-\infty,\infty]$.\\
$\hat\xi^K_t$, $\hat\xi^{(T)}_t$\;: & 
the same as $\xi^K_t$, $\xi^{(T)}_t$, but defined in
terms of the dual Brownian net $\hat\Ni$.\\
*-mesh\;: & see Definition \ref{D:starmesh}
\end{tabular}

\section{Separation points}\label{S:inter}

\subsection{Interaction between forward and dual paths}\label{S:inter1}

We know from \cite{STW00} that paths in $\Wl$ and $\hat\Wl$, resp.\
$\Wr$ and $\hat\Wr$, interact by Skorohod reflection. More precisely,
conditioned on $\hat l\in\hat\Wl$ with deterministic starting point
$\hat z=(\hat x,\hat t)$, the path $l\in\Wl$ with deterministic
starting point $z=(x,t)$, where $t<\hat t$, is distributed as a Brownian
motion with drift $-1$ Skorohod reflected off $\hat l$. It turns out
that conditioned on $\hat l\in\hat\Wl$, the interaction of $r\in
\Wr$ starting at $z$ with $\hat l$ is also Skorohod reflection, albeit
$r$ may cross $\hat l$ at a random time.

\bl\label{L:forback}{\bf[Interaction between dual and forward paths]}\\
Let $\hat l\in\hat \Wl$, resp.\ $r\in\Wr$, start from deterministic points
$\hat z=(\hat x,\hat t)$, resp.\ $z=(x, t)$,
with $t<\hat t$. Then almost surely,
\begin{itemize}
\item[{\rm (a)}] Conditioned on $(\hat l(s))_{t\leq s\leq \hat t}$
with $\hat l(t)<x$, the path $r$ is given (in distribution) by the unique
solution to the Skorohod equation
\be\label{roffhatl}
\begin{aligned}
\di r(s)&=\di B(s)+\di s+\di\De(s),\qquad\qquad & t\leq s&\leq \hat t,\\
\di r(s)&=\di B(s)+\di s,                      & \hat t\leq s&,
\end{aligned}
\ee
where $B$ is a standard Brownian motion, $\De$ is a nondecreasing
process increasing only when $r(s)=\hat l(s)$ (i.e, $\int_t^{\hat t}
1_{\{r(s)\neq\hat l(s)\}}d\De(s)=0$), and $r$ is subject to the
constraint that $\hat l(s)\leq r(s)$ for all $t\leq s\leq \hat t$.

\item[{\rm (b)}] Conditioned on $(\hat l(s))_{t\leq s\leq \hat t}$ with
$x<\hat l(t)$, the path $r$ is given (in distribution) by the unique
solution to the Skorohod equation
\be\label{roffhatl2}
\begin{aligned}
\di r(s)&=\di B(s)+\di s-\di\De(s),\qquad\qquad
&t\leq s&\leq \hat t,\ \De(s)<T,\\
\di r(s)&=\di B(s)+\di s +\di\De(s),
&t\leq s&\leq \hat t,\ T\leq\De(s),\\
\di r(s)&=\di B(s)+\di s,
&\hat t\leq s&,
\end{aligned}
\ee
where $B$ is a standard Brownian motion, $\De$ is a nondecreasing
process increasing only when $r(s)=\hat l(s)$, $T$ is an independent
mean $1/2$ exponential random variable, and $r$ is subject to the
constraints that $r(s)\leq\hat l(s)$ resp.\ $\hat l(s)\leq r(s)$ for
all $s\in[t,\hat t]$ such that $\De(s)<T$ resp.\ $T\leq\De(s)$.
\end{itemize}
The interaction between paths in $\Wl$ and $\hat\Wr$ is similar by
symmetry. If $z=(\hat l(t),t)$, where $t$ is a deterministic time, then
exactly two paths $r_1,r_2\in\Wr$ start from $z$, where one solves
(\ref{roffhatl}) and the other solves (\ref{roffhatl2}). Conditional on $\hat
l$, the paths $r_1$ and $r_2$ evolve independently up to the first time they
meet, at which they coalesce.
\el
{\bf Remark} Lemma~\ref{L:forback} gives an almost sure construction
of a pair of paths $(\hat l,r)$ starting from deterministic points
$\hat z$ and $z$. By the same argument as in~\cite{STW00}, we can
extend Lemma~\ref{L:forback} to give an almost sure construction of a
finite collection of paths in $(\hat\Wl, \Wr)$ with deterministic
starting points, and the order in which the paths are constructed is
irrelevant. (The technical issues involved in consistently defining
multiple coalescing-reflecting paths in our case are the same as those
in~\cite{STW00}.) Therefore, by Kolmogorov's extension theorem,
Lemma~\ref{L:forback} may be used to construct $(\hat\Wr,\Wl)$
restricted to a countable dense set of starting points in $\R^2$.
Taking closures and duals, this gives an alternative construction of
the double left-right Brownian web $(\Wl,\Wr,\hat\Wl,\hat\Wr)$, and
hence of the Brownian net.\med

\noi
{\bf Proof of Lemma~\ref{L:forback}} We will prove the following
claim. Let $\hat z=(\hat x,\hat t)$ and $z=(x, t)$ be deterministic
points in $\R^2$ with $t<\hat t$ and let $\hat l\in\hat \Wl$ and
$r\in\Wr$ be the a.s.\ unique dual left-most and forward right-most
paths starting from $\hat z$ and $z$, respectively. We will show that
it is possible to construct a standard Brownian motion $B$ and a mean
$1/2$ exponential random variable $T$ such that $\hat l,B$, and $T$ are
independent, and such that $r$ is the a.s.\ unique solution to the
equation
\be\label{rint}
r(s)=W(s)+\int_t^s(1_{\{\tau\leq u\}}-1_{\{u<\tau\}})\,\di\De(u),
\quad \quad s\geq t,
\ee
where $W(s):=x+B(s-t)+(s-t)$, $s\geq t$, is a Brownian motion with
drift 1 started at time $(x,t)$, $\De$ is a nondecreasing process
increasing only at times $s\in[t,\hat t]$ when $r(s)=\hat l(s)$,
\be\label{chi}
\tau:=\left\{\ba{ll}
\inf\{s\in[t,\hat t]:\De(s)\geq T\}\quad&\mbox{if }x<\hat l(t),\\
t\quad&\mbox{if }\hat l(t)\leq x,
\ea\right.
\ee
and $r$ is subject to the constraints that $r(s)\leq\hat l(s)$ resp.\
$\hat l(s)\leq r(s)$ for all $s\in[t,\hat t]$ such that $s<\tau$
resp.\ $\tau\leq s$. Note that (\ref{rint}) is a statement about the
joint law of $(\hat l,r)$, which implies the statements about the
conditional law of $r$ given $\hat l$ in Lemma~\ref{L:forback}.

Since $\P[\hat l(t)=x]=0$, to show that (\ref{rint})--(\ref{chi}) has an a.s.\
unique solution, we may distinguish the cases $\hat l(t)<x$ and
$x<\hat l(t)$. If $\hat l(t)<x$, then (\ref{rint}) is a usual Skorohod
equation with reflection (see Section 3.6.C of \cite{KS91}) which is
known to have the unique solution $r(s)=W(s)+\De(s)$, where
\be
\De(s)=\sup_{t\leq u\leq s\wedge\hat t}\big(\hat l(u)-W(u)\big)\vee 0
\qquad(s\geq t).
\ee
If $x<\hat l(t)$, then by the same arguments $r(s)=W(s)-\De(s)$
$(t\leq s\leq\tau)$, where
\be
\De(s)=\sup_{t\leq u\leq s\wedge\hat t}\big(W(u)-\hat l(u)\big)\vee 0
\qquad(t\leq s\leq\tau),
\ee
and
\be
\tau:=\inf\big\{s\geq t:\sup_{t\leq u\leq s\wedge\hat t}
\big(W(u)-\hat l(u)\big)\vee 0=T\big\},
\ee
which may be infinite. Note that for a.e.\ path $\hat l$, the time $\tau$
is a stopping time for $W$. If $\tau<\infty$, then for $s\geq\tau$ our
equation is again a Skorohod equation, with reflection in the other
direction, so $r(s)=W^\tau(s)+\De^\tau(s)$ where
$W^\tau(s):=r(\tau)+W(s)-W(\tau)$ $(s\geq\tau)$, and
\be
\De^\tau(s)=\sup_{\tau\leq u\leq s\wedge\hat t}
\big(\hat l(u)-W^\tau(u)\big)\vee 0
\qquad(s\geq\tau).
\ee

To prove that there exist $W$ and $T$ such that $r$ solves
(\ref{rint}), we follow the approach in~\cite{STW00} and use
discrete approximation. First, we recall from \cite{SS06} the discrete
system of branching-coalescing random walks on $\Z^2_{\rm
even}=\{(x,t) :x,t\in\Z, x+t \text{ is even}\}$ starting from every
site of $\Z^2_{\rm even}$. Here, for $(x,t)\in\Z^2_{\rm even}$, the
walker that is at time $t$ in $x$ jumps at time $t+1$ with probability
$\frac{1-\epsilon}{2}$ to $x-1$, with the same probability to $x+1$,
and with the remaining probability $\epsilon$ branches in two walkers
situated at $x-1$ and $x+1$. Random walks that land on the same
position coalesce immediately. Following \cite{SS06}, let
$\Ui_\epsilon$ denote the set of branching-coalescing random walk
paths on $\Z^2_{\rm even}$ (linearly interpolated between integer
times), and let $\Ui^{\rm l}_\epsilon$, resp.\ $\Ui^{\rm r}_\epsilon$,
denote the set of left-most, resp.\ right-most, paths in $\Ui_\epsilon$
starting from each $z\in\Z^2_{\rm even}$. There exists a natural dual
system of branching-coalescing random walks on $\Z^2_{\rm odd} =
\Z^2\backslash \Z^2_{\rm even}$ running backward in time, where
$(x,t)\in\Z^2_{\rm even}$ is a branching point in the forward system
if and only if $(x,t+1)$ is a branching point in the backward system,
and otherwise the random walk jumping from $(x,t)$ in the forward
system and the random walk jumping from $(x,t+1)$ in the backward
system are coupled so that they do not cross. Denote the dual
collection of branching-coalescing random walk paths by
$\hat\Ui_\eps$, and let $\hat\Ui^{\rm l}_\eps$, resp.\ $\hat\Ui^{\rm
r}_\eps$, denote the dual set of left-most, resp.\ right-most, paths in
$\hat\Ui_\eps$. Let $S_\eps:\R^2\to\R^2$ be the diffusive scaling map
$S_\eps (x,t)=(\eps x,\eps^2 t)$, and define $S_\eps$ applied to a
subset of $\R^2$, a path, or a set of paths analogously.

Let $\eps_n$ be a sequence satisfying $\eps_n \downarrow 0$. Choose
$z_n=(x_n,t_n)\in\Z^2_{\rm even}$ and $\hat z_n=(\hat x_n,\hat
t_n)\in\Z^2_{\rm odd}$ such that $S_{\epsilon_n}(z_n,\hat z_n)
\to(z,\hat z)$ as $n\to\infty$. If we denote by $\hat l_n$ the unique
path in $\hat\Ui_{\eps_n}^{\rm l}(\hat z_n)$ and by $r_n$ the unique
path in $\Ui^{\rm r}_{\eps_n}(z_n)$, then by Theorem~5.2 of
\cite{SS06}, we have
\be
\Li\big(S_{\eps_n}(\hat l_n, r_n)\Big)\Aston{\eps_n}\Li\big(\hat l,r\big),
\ee
where $\Li$ denotes law and $\Rightarrow$ denotes weak convergence.

The conditional law of $r_n$ given $\hat l_n$ has the following
description. Let
\be
I_n:=\{(y,s)\in\Z^2_{\rm even}:t\leq s<\hat t,\ \hat l(s+1)=y\}
\ee
be the set of points where a forward path can meet $\hat l$. Let
$I^{\rm l}_n:=\{(y,s)\in I_n:\hat l(s+1)<\hat l(s)\}$ and $I^{\rm
r}_n:=\{(y,s)\in I_n:\hat l(s)<\hat l(s+1)\}$ be the sets of points
where a forward path can meet $\hat l$ from the left and right,
respectively. Conditional on $\hat l_n$, the process $(r_n(s))_{s\geq
t_n}$ is a Markov process such that if $(r_n(s),s)\not\in I_n$, then
$r_n(s+1)=r_n(s)+1$ with probability $(1+\eps_n)/2$, and
$r_n(s+1)=r_n(s)-1$ with the remaining probability. If $(r_n(s),s)\in
I^{\rm l}_n$, then $r_n(s+1)=r_n(s)+1$ with probability
$2\eps_n/(1+\eps_n)$ and $r_n(s+1)=r_n(s)-1$ with the remaining
probability. Here $2\eps_n/(1+\eps_n)$ is the conditional probability
that $(r_n(s),s)$ is a branching point of the forward random walks
given that $\hat l_n(s+1)<\hat l_n(s)$. Finally, if $(r_n(s),s)\in I^{\rm
r}_n$, then $r_n(s+1)=r_n(s)+1$ with probability 1.

In view of this, we can construct $\hat r_n$ as follows.
Independently of $\hat l_n$, we choose a random walk $(W_n(s))_{s\geq
t_n}$ starting at $(x_n,t_n)$ that at integer times jumps from $y$ to
$y-1$ with probability $(1-\eps_n)/2$ and to $y+1$ with
probability $(1+\eps_n)/2$. Moreover, we introduce i.i.d.\
Bernoulli random variables $(X_n(i))_{i\geq 0}$ with
$\P[X_n(i)=1]=4\eps_n/(1+\eps_n)^2$. We then inductively construct
processes $\De_n$ and $r_n$, starting at time $s$ in $\De_n(s)=0$ and
$r_n(s)=x_n$, by putting
\be\label{Dedisc}
\De_n(s+1):=\left\{\ba{ll}
\De_n(s)+1\quad&\mbox{if }(r_n(s),s)\in I^{\rm l}_n,\ W_n(s+1)>W_n(s),
\mbox{ and }X_{\De_n(s)}=0,\\
\De_n(s)+1\quad&\mbox{if }(r_n(s),s)\in I^{\rm r}_n\mbox{ and }W_n(s+1)<W_n(s),\\
\De_n(s)\quad&\mbox{otherwise,}
\ea\right.\ee
and
\be
r_n(s+1):=\left\{\ba{ll}
r_n(s)+(W_n(s+1)-W_n(s))\quad&\mbox{if $\De_n(s+1)=\De_n(s)$,}\\
r_n(s)-(W_n(s+1)-W_n(s))\quad&\mbox{if $\De_n(s+1)>\De_n(s)$.}
\ea\right.\ee
This says that $r_n$ evolves as $W_n$, but is reflected off $\hat l_n$
with probability $1-4\eps_n/(1+\eps_n)^2$ if it attempts to cross from
left to right, and with probability $1$ if it attempts to cross from
right to left. Note that if $(r_n(s),s)\in I^{\rm l}_n$, then $r_n$
attempts to cross with probability $(1+\eps_n)/2$, hence the
probability that it crosses is
$4\eps_n/(1+\eps_n)^2\cdot(1+\eps_n)/2=2\eps_n/(1+\eps_n)$, as
required.

Extend $W_n(s),\De_n(s)$, and $r_n(s)$ to all real $s\geq t$ by linear
interpolation, and set $T_n:=\inf\{i\geq
0:X_i=1\}$. Then
\be\label{rint2}
r_n(s)=W_n(s)+2\!\int_t^s\!(1_{\{\tau_n\leq u\}}
-1_{\{u<\tau_n\}})\,\di\De_n(u)\qquad(s\geq t),
\ee
where
\be\label{chi2}
\tau_n:=\left\{\ba{ll}
\inf\{s\in[t_n,\hat t_n]:\De_n(s)\geq T_n\}\quad&\mbox{if }x_n<\hat l_n(t),\\
t\quad&\mbox{if }\hat l_n(t)<x_n.
\ea\right.
\ee
The process $r_n$ satisfies $r_n(s)\leq\hat l_n(s)-1$ resp.\ $\hat
l_n(s)+1\leq r_n(s)$ for all $s\in[t_n,\hat t_n]$ such that
$s\leq\tau_n$ resp.\ $\tau_n+1\leq s$, and $\De_n$ increases only for
$s\in[t_n,\hat t_n]$ such that $|r_n(s)-\hat l_n(s)|=1$.

By a slight abuse of notation, set $T_{\eps_n}:=\eps_nT_n$,
$\tau_{\eps_n}:=\eps_n^2\tau_n$, and let $\hat
l_{\eps_n},r_{\eps_n},W_{\eps_n}$, and $\De_{\eps_n}$ be the
counterparts of $l_n,r_n,W_n$, and $\De_n$, diffusively rescaled with
$S_{\eps_n}$. Then
\be\label{basco}
\Li(\hat l_{\eps_n},r_{\eps_n})
\Asto{n}\Li(\hat l,r)\quad\mbox{and}\quad
\Li(\hat l_{\eps_n},W_{\eps_n},2T_{\eps_n})
\Asto{n}\Li(\hat l,W,T),
\ee
where $\hat l,r$ are the dual and forward path we are interested in,
$W$ is a Brownian motion with drift 1, started at $(x_n,t_n)$, $T$ is
an exponentially distributed random variable with mean $1/2$, and
$\hat l,W,T$ are independent.

It follows from the convergence in (\ref{basco}) that the laws of the
processes $\ti\De_{\eps_n}(s):=\ffrac{1}{2}(r_n(s)-W_n(s))$ are
tight. By (\ref{rint2})--(\ref{chi2}), for $n$ large enough, one has
$\ti\De_{\eps_n}=\De_{\eps_n}$ on the event that $\hat l(t)<x$, while on the
complementary event, $\ti\De_{\eps_n}$ is $-\De_{\eps_n}$ reflected at
the level $-T_{\eps_n}$. Using this, it is not hard to see that the
processes $\De_{\eps_n}$ are tight. Therefore, going to a subsequence if
necessary, by Skorohod's representation theorem (see e.g.~Theorem~6.7
in \cite{Bi99}), we can find a coupling such that
\be
(\hat l_{\eps_n},r_{\eps_n},W_{\eps_n},2\De_{\eps_n},2T_{\eps_n})
\asto{n}(\hat l,r,W,\De,T)\quad{\rm a.s.},
\ee
where the paths converge locally uniformly on compacta. Since the
$\De_{\eps_n}$ are nondecreasing and independent of $T_n$, and $T$ has a
continuous distribution, one has $\tau_n\to\tau$ a.s. Taking the limit
in (\ref{rint2})--(\ref{chi2}), it is not hard to see that $r$ solves
the equations (\ref{rint})--(\ref{chi}).

If $z$ is of the form $z=(\hat l(t),t)$ for some deterministic $t$,
then the proof is similar, except that now we consider two
approximating sequences of right-most paths, one started at $(\hat
l_n(t_n)-1,t_n)$ and the other at $(\hat l_n(t_n)+1,t_n)$.\qed

\begin{figure}[tp] 
\centering
\plaat{0.5cm}
{\includegraphics[width=2.5cm]{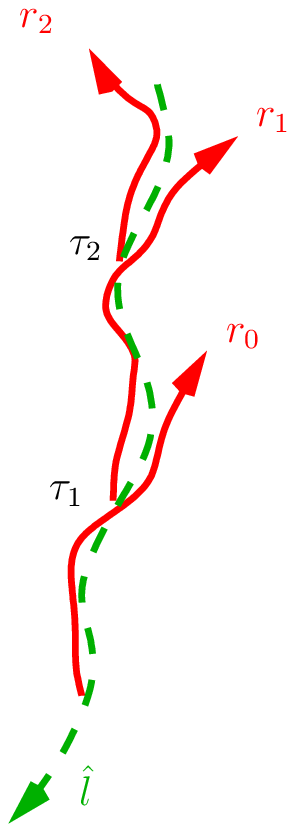}}
\caption{A sequence of paths in $\Wr$ crossing $\hat l\in\hat\Wl$.}
\label{fig:seqcross}
\end{figure}

The next lemma is very similar to Lemma~\ref{L:forback} (see
Figure~\ref{fig:seqcross}).
\bl\label{L:nocro}{\bf[Sequence of paths crossing a dual path]}\\
Let $\hat z=(\hat x,\hat t)\in\R^2$ and $t<\hat t$ be a deterministic point
and time. Let $\hat B,B_i$ $(i\geq 0)$ be independent, standard
Brownian motions, and let $(T_i)_{i\geq 0}$ be independent mean $1/2$
exponential random variables.

Set $\hat l(\hat t-s):=\hat x+\hat B_s+s$ $(s\geq 0)$. Set $\tau_0:=t$ and
define inductively paths $(r_i(s))_{s\geq\tau_i}$ $(i=0,\ldots,M)$ starting at
$r_i(\tau_i)=\hat l(\tau_i)$ by the unique solutions to the Skorohod equation
\be\label{Skorind}
\begin{aligned}
\di r_i(s)&=\di\ti B_i(s)+\di s-\di\De_i(s),\qquad\qquad
&\tau_i\leq s&\leq\tau_{i+1},\\
\di r_i(s)&=\di\ti B_i(s)+\di s +\di\De_i(s),
&\tau_{i+1}\leq s&\leq\hat t,\\
\di r_i(s)&=\di\ti B_i(s)+\di s,
&\hat t\leq s&,
\end{aligned}
\ee
where $\De_i$ is a nondecreasing process increasing only when $r_i(s)=\hat
l(s)$, the process $r_i$ is subject to the constraints that $r_i(s)\leq\hat
l(s)$ on $[\tau_i,\tau_{i+1}]$ and $\hat l(s)\leq r_i(s)$ on $[\tau_{i+1},\hat
  t]$, and $\tau_{i+1}:=\inf\{s\in(\tau_i,\hat t):\De_i(s)>T_i\}$
$(i=0,\ldots,M-1)$. The induction terminates at $M:=\inf\{i\geq 0:\De_i(\hat
t)\leq T_i\}$, and we set $\tau_{M+1}:=\hat t$. The Brownian motions $\ti B_i$
in (\ref{Skorind}) are inductively defined as $\ti B_0:=B_0$, and, for
$i=1,\ldots,M$,
\be\label{tiBi}
\ti B_i(s):=\left\{\ba{l@{\qquad}l}
B_i(s)&\tau_i\leq s\leq\sig_i,\\
\ti B_{i-1}(s)&\sig_i\leq s,
\ea\right.\ee
where $\sig_i:=\inf\{s\geq\tau_i:r_i(s)=r_{i-1}(s)\}$.

Then $M<\infty$ a.s.\ and we can couple $\{\hat l,r_0,\ldots,r_M\}$ to a
left-right Brownian web and its dual in such a way that $\hat l\in\hat\Wl$ and
$r_i\in\Wr$ for $i=0,\ldots,M$.
\el
{\bf Proof} This follows by discrete approximation in the same way
as in the proof of Lemma~\ref{L:forback}. Note that (\ref{tiBi}) ensures that
$r_0,\ldots,r_M$ coalesce when they meet. To see that $M<\infty$ a.s.,
define $(r'(s))_{s\geq t}$ by $r':=r_i$ on $[\tau_i,\tau_{i+1}]$ and $r':=r_M$
on $[\hat t,\infty)$, and set $B(s):=\int_t^s\di B(u)$ and
$\De(s):=\int_t^s\di\De(u)$ where $\di B:=\di\ti B_i$ and $\di\De:=\di\De_i$ on
$[\tau_i,\tau_{i+1}]$ and $\di B:=\di\ti B_M$ on
$[\hat t,\infty)$. Then $B$ is a Brownian motion, $\De$ is a nondecreasing
process increasing only when $r'(s)=\hat l(s)$, and $r'$ solves
\be\label{pseudevol}
\begin{aligned}
\di r'(s)&=\di B(s)+\di s-\di\De(s),\qquad\qquad & t\leq s&\leq \hat t,\\
\di r'(s)&=\di B(s)+\di s,                      & \hat t\leq s&,
\end{aligned}
\ee
subject to the constraint that $r'(s)\leq\hat l(s)$ for all
$t\leq s\leq \hat t$. In particular, we have $\De(\hat t)
=\sup_{s\in[t,\hat t]}(B(s)-\hat l(s)+\hat l(t))<\infty$
and the times $\tau_1,\ldots,\tau_M$ are created by a Poisson point process on
$[t,\hat t]$ with intensity measure~$2\,\di\De$.\qed

\subsection{Structure of separation points}\label{S:inter3}

In this section we apply the results from the previous section to
study the structure of separation points. We start with some definitions.
First, we recall the definition of intersection points from \cite{SS06}, and
formally define meeting and separation points.
\bd\label{D:IMS}{\bf[Intersection, meeting, and separation points]}\\
We call $z=(x,t)\in\R^2$ an intersection point of $\pi_1,\pi_2\in\Pi$
if $\sig_{\pi_1},\sig_{\pi_2}<t$ and $\pi_1(t)=\pi_2(t)=x$. If
furthermore $\pi_1\neq\pi_2$ on $(t-\eps,t)$ (resp.\ $(t,t+\eps)$) for
some $\eps>0$, then we call $z$ a meeting (resp.\ separation) point of
$\pi_1$ and $\pi_2$. Intersection, meeting, and separation points of
dual paths are defined analogously.
\ed
Crossing points of two forward paths $\pi_1,\pi_2\in\Pi$ have been defined in
\cite{SS06}. Below, we define crossing points of a forward path $\pi$ and a
dual path $\hat\pi$.
\bd\label{D:cros}{\bf[Crossing points]}\\
We say that a path $\pi\in\Pi$ crosses a path $\hat\pi\in\hat\Pi$ from left to
right at time $t$ if there exist $\sig_\pi\leq
t_-<t<t_+\leq\hat\sig_{\hat\pi}$ such that $\pi(t_-)<\hat\pi(t_-)$,
$\hat\pi(t_+)<\pi(t_+)$, and $t=\inf\{s\in(t_-,t_+):\hat\pi(s)<\pi(s)\}
=\sup\{s\in(t_-,t_+):\pi(s)<\hat\pi(s)\}$. Crossing from right to left is
defined analogously. We call $z=(x,t)\in\R^2$ a crossing point of $\pi\in\Pi$
and $\hat\pi\in\hat\Pi$ if $\pi(t)=x=\hat\pi(t)$ and $\pi$ crosses $\hat\pi$
either from left to right or from right to left at time $t$.
 \ed
A disadvantage of the way we have defined crossing is that it is possible to
find paths $\pi\in\Pi$ and $\hat\pi\in\hat\Pi$ with
$\sig_\pi<\hat\sig_{\hat\pi}$, $\pi(\sig_\pi)<\hat\pi(\sig_\pi)$, and
$\hat\pi(\hat\sig_{\hat\pi})<\pi(\hat\sig_\pi)$, such that $\pi$ crosses
$\hat\pi$ from left to right at no time in $(\sig_\pi,\hat\sig_{\hat\pi})$.
The next lemma shows, however, that such pathologies do not happen for
left-most and dual right-most paths.
\bl\label{L:cros}{\bf[Crossing times]}\\
Almost surely, for each $r\in\Wr$ and $\hat l\in\hat\Wl$ such that
$\sig_r<\hat\sig_{\hat l}$, $r(\sig_r)<\hat l(\sig_r)$, and $\hat
l(\hat\sig_{\hat l})<r(\hat\sig_{\hat l})$, there exists a unique
$\sig_r<\tau<\hat\sig_{\hat l}$ such that $r$ crosses $\hat l$ from left to
right at time $\tau$. Moreover, one has $r\leq\hat l$ on $[\sig_r,\tau]$,
$\hat l\leq r$ on $[\tau,\hat\sig_{\hat l}]$, and there exist
$\eps_n,\eps'_n\down 0$ such that $r(\tau-\eps_n)=\hat l(\tau-\eps_n)$ and
$r(\tau+\eps'_n)=\hat l(\tau+\eps'_n)$. Analogous statements hold for
left-most paths crossing dual right-most paths from right to left.
\el
{\bf Proof} By \cite[Lemma~3.4~(b)]{SS06} it suffices to prove the statements
for paths with deterministic starting points. Set
\bc
\tau&:=&\dis\sup\{s\in(\sig_r,\hat\sig_{\hat l}):r(s)<\hat l(s)\},\\[5pt]
\tau'&:=&\dis\inf\{s\in(\sig_r,\hat\sig_{\hat l}):\hat l(s)<r(s)\}.
\ec
Lemma~\ref{L:forback} and the properties of Skorohod reflection imply
that $\tau=\tau'$, hence $r$ crosses $\hat l$ from left to right at the unique
crossing time $\tau$. To see that $r(\tau-\eps_n)=\hat l(\tau-\eps_n)$ and
for some $\eps_n\down 0$, we consider, in Lemma~\ref{L:forback}~(b),
the unique solution $(\ti r,\ti\De)$ to the Skorohod equation
\be\label{roffhat3}
\begin{aligned}
\di\ti r(s)&=\di B(s)+\di s-\di\ti\De(s),\qquad\qquad
& t\leq s&\leq \hat t,\\
\di\ti r(s)&=\di B(s)+\di s,   & \hat t\leq s&.
\end{aligned}
\ee
Then $(\ti r(s),\ti\De(s))=(r(s),\De(s))$ for all $s\leq\tau$ and $\ti\De$ is
independent of $T$. It follows that the set $\{s\in[t,\hat t]:\ti\De(s)=T\}$,
if it is nonempty, a.s.\ consists of one point, hence $\ti\De(s)<T$ for all
$s<\tau$, which implies our claim. By symmetry between forward and dual, and
between left and right paths, we also have $r(\tau+\eps'_n)=\hat
l(\tau+\eps'_n)$ for some $\eps'_n\down 0$.\qed

The next proposition says that the sets of crossing points of $(\hat\Wl,\Wr)$
and $(\Wl,\hat\Wr)$, of separation points of $(\Wl,\Wr)$ and
$(\hat\Wl,\hat\Wr)$, and of points of type $\rm(p,pp)_s/(p,pp)_s$ in
$\Ni/\hat\Ni$ as defined in Section~\ref{S:class}, all coincide.
\bp\label{P:sep}{\bf[Separation points]}\\
Almost surely for each $z\in\R^2$, the following statements are equivalent:
\begin{itemize}
\item[\rm(i)] $z$ is a crossing point of some $l\in\Wl$ and
$\hat r\in\hat\Wr$,
\item[\rm(ii)] $z$ is a crossing point of some $\hat l\in\hat\Wl$ and
$r\in\Wr$,
\item[\rm(iii)] $z$ is a separation point of some $l\in\Wl$ and $r\in\Wr$,
\item[\rm(iv)] $z$ is a separation point of some $\hat l\in\hat\Wl$ and
$\hat r\in\hat\Wr$,
\item[\rm(v)] $z$ is of type $\rm(p,pp)_s$ in $\Ni$.
\item[\rm(vi)] $z$ is of type $\rm(p,pp)_s$ in $\hat\Ni$.
\end{itemize}
Moreover, the set $\{z\in\R^2:z\mbox{ is of type $\rm(p,pp)_s$ in
  $\Ni$}\}$ is a.s.\ countable.
\ep
{\bf Proof} We will prove the implications
(i)$\volgt$(iii)$\volgt$(v)$\volgt$(ii). By symmetry between forward and dual
paths, this then implies that (ii)$\volgt$(iv)$\volgt$(vi)$\volgt$(i), hence
all conditions are equivalent.

To prove (i)$\volgt$(iii), let $z=(x,t)$ be a crossing point of $\hat
r\in\hat\Wr$ and $l\in\Wl$. By \cite[Lemma~3.4~(b)]{SS06}, without loss of
generality, we can assume that $\hat r$ and $l$ start from deterministic
points with $\sigma_l<t<\hat\sig_{\hat r}$. By Lemma~\ref{L:nocro} and the
fact that paths in $\hat\Wr$ cannot cross, there exists an $\hat
r'\in\hat\Wr(z)$ and $\eps>0$ such that $l\leq\hat r'$ on $[t-\eps,t]$. By
Lemma~\ref{L:cros}, we can find $s\in(t-\eps,t)$ such that $l(s)<\hat
r'(s)$. Now any path $r\in\Wr$ started at a point $(y,s)$ with $l(s)<x<\hat
r'(s)$ is confined between $l$ and $\hat r'$, hence passes through $z$. Since
$r$ cannot cross $\hat r$ we have $\hat r\leq r$ on $[t,\hat\sig_{\hat
    r}]$. Since $l$ and $r$ spend positive Lebesgue time together whenever
they meet by \cite[Prop.~3.6~(c)]{SS06}, while $r$ and $\hat r$ spend zero
Lebesgue time together by \cite[Prop.~3.2~(d)]{SS06}, $z$ must be a separation
point of $l$ and $r$.

To prove (iii)$\volgt$(v), let $z=(x,t)$ be a separation point of
$l\in\Wl$ and $r\in\Wr$ so that $l(s)<r(s)$ on $(t,t+\eps]$ for some
$\eps>0$. Without loss of generality, we can assume that $l$ and $r$
start from deterministic points with $\sigma_l,\sigma_r<t$. Choose
$\hat t\in(t,t+\eps]\cap\Ti$ where $\Ti\sub\R$ is some fixed,
deterministic countable dense set. By Lemma~\ref{L:nocro}, there exist
unique $\hat t>\tau_1>\cdots>\tau_M>\sig_l$, $\tau_0:=\hat t$,
$\tau_{M+1}:=-\infty$, with $0\leq M<\infty$, and $\hat
r_0\in\hat\Wr(l(\tau_0),\tau_0),\ldots,\hat
r_M\in\Wr(l(\tau_M),\tau_M)$, such that $l\leq\hat r_i$ on
$[\tau_{i+1}\vee\sig_l,\tau_i]$ for each $0\leq i\leq M$, and
$\tau_{i+1}$ is the crossing time of $\hat r_i$ and $l$ for each
$0\leq i<M$. (See Figure~\ref{fig:seqcross}, turned upside down.)
Note that all the paths $\hat r_i$ are confined to the left of $r$
because paths in $\Wr$ and $\hat\Wr$ cannot cross. Since $M<\infty$,
one of the paths, say $\hat r=\hat r_i$, must pass through the
separation point $z$. Since by Proposition 2.2~(b) of \cite{SS06}, $l$
and $r$ spend positive Lebesgue time together on $[t-\eps',t]$ for all
$\eps'>0$, and by Proposition 3.1~(d) of \cite{SS06}, $\hat r$ and $r$
spend zero Lebesgue time together, $\hat r$ must cross $l$ at $z$.
Similarly, there exists a $\hat l\in\hat\Wl$ starting at $(r(s'),s')$
for some $s'\in(t,t+\eps]$ such that $\hat l\leq r$ on $[t,s']$ and
$\hat l$ crosses $r$ in the point $z$. Again because paths in $\Wi$
and $\hat\Wi$ spend zero Lebesgue time together, $z$ must be a
separation point of $\hat r$ and $\hat l$.

For each point $(y,s)$ with $l(s)<y<\hat r(s)$ we can find a path
$r\in\Wr(y,s)$ that is confined between $l$ and $\hat r$, so using the
fact that $\Wr$ is closed we see that there exists a path
$r'\in\Wr(z)$ that is confined between $l$ and $\hat r$. By
Lemma~\ref{L:cros}, there exist $\eps_n\down 0$ such that
$l(t+\eps_n)=\hat r(t+\eps_n)$, so $l\sim^z_{\rm out}r'$. Similarly,
there exists $l'\in\Wl$ starting from $z$ which is confined between
$\hat l$ and $r$ and satisfies $l'\sim^z_{\rm out}r$. The point $z$
must therefore be of type $(1,2)_{\rm l}$ in $\Wl$ and of type
$(1,2)_{\rm r}$ in $\Wr$, and hence of type $\rm(p,pp)_s$ in
$\Ni$.

To prove (v)$\volgt$(ii), let $l\in\Wl$ and $r\in\Wr$ be the left-most and
right-most paths separating at $z$ and let $r'\in\Wr(z)$ and $l'\in\Wl(z)$ be
the right-most and left-most paths such that $l\sim^z_{\rm out}r'$ and
$l'\sim^z_{\rm out}r$. Then, by \cite[Prop.~3.2~(c) and Prop.~3.6~(d)]{SS06},
any $\hat l\in\hat\Wl$ started in a point $z'=(x',t')$ with $t<t'$ and
$r'(t')<x'<l'(t')$ is contained between $r'$ and $l'$, hence must pass through
$z$. Since $\hat l$ cannot cross $l$ and, by \cite[Prop.~3.2~(d)]{SS06},
spends zero Lebesgue time with $l$, while by \cite[Prop.~3.6~(c)]{SS06}, $l$
and $r$ spend positive Lebesgue time in $[t-\eps,t]$ for any $\eps>0$, the
path $\hat l$ must cross $r$ in $z$.

The fact that the set of separation points is countable, finally, follows from
the fact that by \cite[Lemma~3.4~(b)]{SS06}, each separation point between
some $l\in\Wl$ and $r\in\Wr$ is the separation point of some $l'\in\Wl(\Di)$
and $r'\in\Wr(\Di)$, where $\Di$ is a fixed, deterministic, countable dense
subset of $\R^2$.\qed

\subsection{Relevant separation points}\label{S:inter4}

In a sense, separation points are the most important points in the
Brownian net, since these are the points where paths have a choice
whether to turn left or right. In the present section, we prove that
for deterministic times $s<u$, there are a.s.\ only locally finitely
many separation points at which paths in $\Ni$ starting at time $s$
have to make a choice that is relevant for determining where they end
up at time $u$.

We start with a useful lemma.
\bl\label{L:NT}{\bf[Incoming net paths]}\\
Almost surely for each $-\infty\leq s<u<\infty$ and $-\infty<x_-\leq
x_+<\infty$:
\begin{itemize}
\item[\rm(a)] For each $\pi\in\Ni$ such that $\sig_\pi=s$ and
$x_-\leq\pi(u)\leq x_+$ there exist $\hat r\in\hat\Wr(x_-,u)$ and
$\hat l\in\hat\Wl(x_+,u)$ such that $\hat r\leq\pi\leq\hat l$ and
$\hat r<\hat l$ on $(s,u)$.

\item[\rm(b)] If there exist $\hat r\in\hat\Wr(x_-,u)$ and $\hat
l\in\hat\Wl(x_+,u)$ such that $\hat r<\hat l$ on $(s,u)$, then there
exists a $\pi\in\Ni$ such that $\sig_\pi=s$ and $\hat r\leq\pi\leq\hat l$
on $(s,u)$.
\end{itemize}
\el
{\bf Proof} Part~(b) follows from the steering argument used in
\cite[Lemma~4.7]{SS06}. To prove part~(a), choose $x^{(n)}_-\up x_-$,
$x^{(n)}_-\down x_-$, $\hat r_n\in\hat\Wr(x^{(n)}_-,u)$ and $\hat
l_n\in\hat\Wl(x^{(n)}_+,u)$. Since paths in the Brownian net do not
enter wedges from outside (see Theorem~\ref{T:netchar}~(b2)), one has
$\hat r_n\leq\pi\leq\hat l_n$ and $\hat r_n<\hat l_n$ on $(s,u)$. By
monotonicity, $\hat r_n\up\hat r$ and $\hat l_n\down\hat l$ for some
$\hat r\in\hat\Wr(x_-,u)$ and $\hat l\in\hat\Wl(x_+,u)$. The claim now
follows from the nature of convergence of paths in the Brownian web
(see \cite[Lemma~3.4~(a)]{SS06}).\qed

The next lemma introduces our objects of interest.
\bl\label{L:rel}{\bf[Relevant separation points]}\\
Almost surely, for each $-\infty\leq s<u\leq\infty$ and
$z=(x,t)\in\R^2$ with $s<t<u$, the following statements are
equivalent:
\begin{itemize}
\item[\rm(i)] There exists a $\pi\in\Ni$ starting at time $s$ such
that $\pi(t)=x$ and $z$ is the separation point of some $l\in\Wl$ and
$r\in\Wr$ with $l<r$ on $(t,u)$.
\item[\rm(ii)] There exists a $\hat\pi\in\Ni$ starting at time $u$ such
that $\hat\pi(t)=x$ and $z$ is the separation point of some $\hat
l\in\Wl$ and $\hat r\in\hat\Wr$ with $\hat r<\hat l$ on $(s,t)$.
\end{itemize}
\el
{\bf Proof} By symmetry, it suffices to prove (i)$\volgt$(ii). If $z$
satisfies (i), then by Lemma~\ref{L:NT}, there exists a
$\hat\pi\in\Ni$ starting at time $u$ such that $\hat\pi(t)=x$, and
there exist $\hat l\in\Wl(z)$, $\hat r\in\Wr(z)$ such that $\hat
r<\hat l$ on $(s,t)$. Since $z$ is the separation point of some
$l\in\Wl$ and $r\in\Wr$, by Proposition~\ref{P:sep}, $z$ is also the
separation point of some $\hat l'\in\Wl$ and $\hat r'\in\Wr$.  Again
by Proposition~\ref{P:sep}, $z$ is of type $\rm(p,pp)_s$, hence we
must have $\hat l'=\hat l$ and $\hat r'=\hat r$ on $[-\infty,t]$.\qed

\noi
If $z$ satisfies the equivalent conditions from Lemma~\ref{L:rel}, then, in
line with the definition given in Section~\ref{S:outline}, we say that $z$ is
an {\em $(s,u)$-relevant separation point}.  We will prove that for
deterministic $S<U$, the set of $(S,U)$-relevant separation points is
a.s.\ locally finite. Let
$\Phi(x):=\frac{1}{\sqrt{2\pi}}\int_{-\infty}^xe^{-y^2/2}\di y$ denote the
distribution function of the standard normal distribution and set
\be
\Psi(t):=\frac{e^{-t}}{\sqrt{\pi t}}+2\Phi(\sqrt{2t})\qquad(0<t\leq\infty).
\ee
For any $-\infty\leq S<U\leq\infty$, we write
\be\label{RSU}
R_{S,U}:=\big\{z\in\R^2:
z\mbox{ is an $(S,U)$-relevant separation point}\big\}.
\ee
Below, $|A|$ denotes the cardinality of a set $A$.
\bp\label{P:finrel}{\bf[Density of relevant separation points]}\\
For deterministic $-\infty\leq S\leq s<u\leq U\leq\infty$ and
$-\infty<a<b<\infty$,
\be\label{reldens}
\E\big[\big|R_{S,U}\cap([a,b]\times[s,u])\big|\big]
=2(b-a)\int_s^u\Psi(t-S)\Psi(U-t)\di t.
\ee
In particular, if $-\infty<S<U<\infty$, then $R_{S,U}$ is a.s.\ a
locally finite subset of $\R\times[S,U]$.
\ep
{\bf Proof} It suffices to prove the statement for deterministic $S<s<u<U$;
the general statement then follows by approximation. Set
\be
E_s=\{\pi(s):\pi\in\Ni,\ \sigma_\pi=S\}
\quad\mbox{and}\quad
F_u=\{\hat\pi(u):\hat\pi\in\hat\Ni,\ \hat\sigma_{\hat\pi}=U\}.
\ee
By \cite[Prop.~1.12]{SS06}, $E_s$ and $F_u$ are spatially homogeneous (in
law), locally finite point sets on $\R$ with densities $\Psi(s-S)$ and
$\Psi(U-u)$, respectively.

Since the restrictions of the Brownian net to $\R\times[S,s]$ and
$\R\times[s,u]$ are independent (which follows from the discrete
approximation in \cite[Thm.~1.1]{SS06}), at each point $x\in E_s$
there start a.s.\ unique paths $l_{(x,s)}\in\Wl(x,s)$ and
$r_{(x,s)}\in\Wr(x,s)$. Likewise, for each $y\in F_u$ there start
a.s.\ unique $\hat r_{(y,u)}\in\hat\Wr(y,u)$ and $\hat
l_{(y,u)}\in\hat\Wl(y,u)$. Let
\bc
\dis Q_{s,u}
&:=&\dis\big\{(x,y):
x\in E_s,\ y\in F_u,\ l_{(x,s)}(u)<y<r_{(x,s)}(u)\big\}\\[5pt]
&=&\dis\big\{(x,y):
x\in E_s,\ y\in F_u,\ \hat r_{(y,u)}(s)<x<\hat l_{(y,u)}(s)\big\},
\ec
where the equality follows from the fact that paths in $\Wl$ and
$\hat\Wl$ a.s.\ do not cross. Note that if $z\in R\times(s,u)$ is an
$(S,U)$-relevant separation point and $\pi,\hat\pi$ are as in
Lemma~\ref{L:rel}, then $(\pi(s),\hat\pi(u))\in Q_{s,u}$. Conversely,
if $(x,y)\in Q_{s,u}$, then the point $z=(w,\tau)$ defined by
$\tau:=\sup\{t\in(s,u):l_{(x,s)}(t)=r_{(x,s)}(t)\}$ and
$w:=l_{(x,s)}(\tau)$ is an $(S,U)$-relevant separation point.

Conditional on $E_s$, for each $x\in E_s$, the paths $l_{(x,s)}$ and
$r{(x,s)}$ are Brownian motions with drift $-1$ and $+1$,
respectively, hence $\E[r{(x,s)}(u)-l_{(x,s)}(u)]=2(u-s)$. Since the
restrictions of the Brownian net to $\R\times[S,s]$, $\R\times[s,u]$,
and $\R\times[u,U]$ are independent, and since the densities of $E_s$
and $F_u$ are $\Psi(s-S)$ and $\Psi(U-u)$, respectively, it follows
that for each $a<b$,
\be\label{Qdens}
\E\big[\big|(x,y)\in Q_{s,u}:x\in(a,b)\}|]=2(b-a)(u-s)\Psi(s-S)\Psi(U-u).
\ee
For $n\geq 1$, set $D_n:=\{S+k(U-S)/n:0\leq k\leq n-1\}$, and for
$t\in(S,U)$ write $\lfloor t\rfloor_n:=\sup\{t'\in D_n:t'\leq t\}$.
By our previous remarks and the equicontinuity of the Brownian net,
for each $z=(x,t)\in R_{S,U}$ there exist $(x_n,y_n)\in Q_{\lfloor
t\rfloor_n,\lfloor t\rfloor_n+1/n}$ such that $x_n\to x$ and $y_n\to
x$. It follows that for any $S<s<u<U$ and $a<b$,
\be
\big|R_{S,U}\cap((a,b)\times(s,u))\big|
\leq\liminf_{n\to\infty}
\big|\{(x,y)\in Q_{t,t+1/n}:t\in D_n\cap(s,u),\ x\in(a,b)\}|,
\ee
and therefore, by Fatou,
\bc
\dis\E\big[\big|R_{S,U}\cap((a,b)\times(s,u))\big|\big]
&\leq&\dis\lim_{n\to\infty}\sum_{t\in D_n\cap(s,u)}
\E\big[\big|\{(x,y)\in Q_{t,t+1/n}:x\in(a,b)\}|]\\[5pt]
&=&\dis 2(b-a)\int_s^u\Psi(t-S)\Psi(U-t)\di t,
\ec
where in the last step we have used (\ref{Qdens}) and Riemann sum
approximation. This proves the inequality $\leq$ in
(\ref{reldens}). In particular, our argument shows that the set
$R_{S,T}\cap((a,b)\times(S,T))$ is a.s.\ finite.

By our previous remarks, each point $(x,y)\in Q_{t,t+1/n}$ gives rise
to an $(S,U)$-relevant separation point $z\in\R\times(t,t+1/n)$. To
get the complementary inequality in (\ref{reldens}), we have to deal
with the difficulty that a given $z$ may correspond to more than one
$(x,y)\in Q_{t,t+1/n}$. For $\de>0$, set
\be
E^\de_t:=\big\{x\in E_t:E_t\cap(x-\de,x+\de)=\{x\}\big\}.
\ee
and define $F^\de_t$ similarly. By Lemma~\ref{L:tech} below, for each
$\eps>0$, we can find $\de>0$ such that
\bc
\dis\E[|E^\de_t\cap(a,b)|]&\geq&\dis(1-\eps)(b-a)\Psi(t-S),\\[5pt]
\dis\E[|F^\de_t\cap(a,b)|]&\geq&\dis(1-\eps)(b-a)\Psi(U-t)
\ec
for all $t\in[s,u]$. Arguing as before, we find that
\be\label{Ede}
\E\big[\big|\{(x,y)\in Q_{t,t'}:x\in(a,b)\cap E^\de_t\}|]
\geq2(1-\eps)(b-a)(t'-t)\Psi(t-S)\Psi(U-t')
\ee
for all $s\leq t\leq t'\leq u$. Similarly, by symmetry between forward
and dual paths,
\be\ba{l}\label{Fde}
\dis\E\big[\big|\{(x,y)\in Q_{t,t'}:x\in(a,b),\ y\in F^\de_{t'}\}|]\\[5pt]
\dis\quad=\E\big[\big|\{(x,y)\in Q_{t,t'}:y\in(a,b)\cap F^\de_{t'}\}|]
\geq2(1-\eps)(b-a)(t'-t)\Psi(t-S)\Psi(U-t').
\ec
Set
\be
Q^{\de,K}_{t,t'}:=\big\{(x,y)\in Q_{t,t'}:x\in E^\de_t,\ y\in F^\de_{t'},
\ x-K\leq l_{(x,t)}\leq r_{(x,t)}\leq x+K\mbox{ on }[t,t']\big\},
\ee
where $l_{(x,t)}$ and $r_{(x,t)}$ are the a.s.\ unique left-most and
right-most paths starting from $(x,t)$. Then, by (\ref{Ede}) and (\ref{Fde}),
for each $\eps>0$ we can choose $\de>0$ and $K>0$ such that
\be\label{Elow}
\E\big[\big|\{(x,y)\in Q^{\de,K}_{t,t'}:x\in(a,b)\}|]
\geq2(1-3\eps)(b-a)(t'-t)\Psi(t-S)\Psi(U-t').
\ee
By the equicontinuity of the net,
\be\ba{l}\label{lowbd}
\dis\big|R_{S,U}\cap((a,b)\times(s,u))\big|\\[5pt]
\dis\quad\geq\limsup_{n\to\infty}
\big|\{(x,y)\in Q^{\de,K}_{t,t+1/n}:t\in D_n\cap(s,u),\ x\in(a,b)\}\big|.
\ec
Since the random variables on the right-hand side of (\ref{lowbd}) are bounded
{f}rom above by the integrable random variable
\be
\big((b-a)/\de+1\big)\big|R_{S,U}\cap\big((a-K,b+K)\times(s,u)\big)\big|,
\ee
we can take expectations on both sides of (\ref{lowbd}) and let $\eps\to 0$,
to get the lower bound in (\ref{reldens}).

The final statement of the proposition follows by observing that the integral
on the right-hand side of (\ref{reldens}) is finite if
$-\infty<S=s<u=U<\infty$.\qed

To complete the proof of Proposition~\ref{P:finrel}, we need the following
lemma.
\bl\label{L:tech}{\bf[Uniform finiteness]} \\
For $0\leq t\leq\infty$, set $\xi^{(-t)}_0:=\{\pi(0):\pi\in\Ni,\ \sig_\pi=-t\}$
and $\xi^{(-t)}_{0,\de}:=\{x\in\xi^{(-t)}_0:\xi^{(-t)}_0
\cap(x-\de,x+\de)=\{x\}\}$. Then, for each compact set $K\sub(0,\infty]$ and
$-\infty<a<b<\infty$, one has
\be\label{tech}
\lim_{\de\down 0}\sup_{t\in K}
\E\big[|(\xi^{(-t)}_0\beh\xi^{(-t)}_{0,\de})\cap(a,b)|\big]=0.
\ee
\el
{\bf Proof} Set $F_t:=|\xi^{(-t)}_0\cap(a,b)|$,
$F^\de_t:=|\xi^{(-t)}_{0,\de}\cap(a,b)|$, and $f_\de(t):=\E[F_t-F^\de_t]$. By
\cite[Prop.~1.12]{SS06}, $\E[F_t]<\infty$ a.s.\ for all
$t\in(0,\infty]$. Since $F^\de_t\up F_t$ as $\de\down 0$, it follows that
$f_\de(t)=\E[F_t-F^\de_t]\down 0$ as $\de\down 0$, for each
$t\in(0,\infty]$. Since $F_s-F^\de_s\down F_t-F^\de_t$ as $s\up t$, the
$f_\de$ are continuous functions on $(0,\infty]$ decreasing to zero, hence
  $\lim_{\de\down 0}\sup_{t\in K}f_\de(t)=0$ for each compact
  $K\sub(0,\infty]$.\qed

The following simple consequence of Proposition~\ref{P:finrel} will often be
useful.
\bl\label{L:finrel}{\bf[Local finiteness of relevant separation points]}\\
Almost surely, for each $-\infty\leq s<u\leq\infty$, the set $R_{s,u}$ of all
$(s,u)$-relevant separation points is a locally finite subset of
$\R\times(s,u)$.
\el
{\bf Proof} Let $\Ti$ be a deterministic countable dense subset of $\R$. Then,
by Proposition~\ref{P:finrel}, $R_{s,u}$ is a locally finite subset of
$\R\times[s,u]$ for each $s<u$, $s,u\in\Ti$. For general $s<u$, we can choose
$s_n,u_n\in\Ti$ such that $s_n\down s$ and $u_n\up u$. Then $R_{s_n,u_n}\up
R_{s,u}$, hence $R_{s,u}$ is locally finite.\qed

\subsection{Reflected paths}\label{S:inter5}

In this section, we prove Lemma~\ref{L:extl}. We start with two preparatory
lemmas.
\bl\label{L:netint}{\bf[Forward and dual paths spend zero time together]}\\
Almost surely, for any $s<u$, one has
\be\label{piha}
\int_s^u1_{\{\xi^{(s)}_t\cap\hat\xi^{(u)}_t\neq\emptyset\}}\di t=0,
\ee
where, for $t\in(s,u)$,
\bc
\dis\xi^{(s)}_t&:=&\dis\{\pi(t):\pi\in\Ni,\ \sig_\pi=s\},\\[5pt]
\dis\hat\xi^{(u)}_t&:=&\dis\{\hat\pi(t):\hat\pi\in\hat\Ni,\ \hat\sig_{\hat\pi}=u\}.
\ec
In particular, one has $\int_s^u 1_{\{\pi(t)=\hat\pi(t)\}}\di t=0$ for
any $\pi\in\Ni$ and $\hat\pi\in\hat\Ni$ with $\sig_\pi=s$ and
$\hat\sig_{\hat\pi}=u$.
\el
{\bf Proof} It suffices to prove the statement for deterministic
times. In that case, the expectation of the quantity in (\ref{piha})
is given by
\be
\int_s^u\P[\xi^{(s)}_t\cap\hat\xi^{(u)}_t\neq\emptyset]\di t.
\ee
By Proposition~1.12 of \cite{SS06}, $\xi^{(s)}_t$ and $\hat\xi^{(u)}_t$ are
stationary point processes with finite intensity. By the independence
of $\Ni|_{(-\infty,t]}$ and $\hat\Ni|_{[t,\infty)}$, $\xi^{(s)}_t$ and
$\hat\xi^{(u)}_t$ are independent. It follows that
$\P[\xi^{(s)}_t\cap\hat\xi^{(u)}_t\neq\emptyset]=0$ for each $t\in(s,u)$.
Fubini's Theorem then implies (\ref{piha}) almost surely.
\qed

Part~(b) of the next lemma is similar to Lemma~\ref{L:cros}.
\bl\label{L:lpi}{\bf[Crossing of dual net paths]}
\begin{itemize}
\item[\rm(a)] Almost surely, for each $r\in\Wr$ and $\hat\pi\in\hat\Ni$ such
  that $\sig_r<\hat\sig_{\hat\pi}$ and $\hat\pi(\sig_r)<r(\sig_r)$, one has
  $\hat\pi\leq r$ on $[\sig_r,\hat\sig_{\hat\pi}]$.

\item[\rm(b)] Almost surely, for each $r\in\Wr$ and $\hat\pi\in\hat\Ni$ such
  that $\sig_r<\hat\sig_{\hat\pi}<\infty$ and $r(\sig_r)<\hat\pi(\sig_r)$, if
  $S:=\{s\in(\sig_r,\hat\sig_{\hat\pi}):r(s)=\hat\pi(s),\ (r(s),s)\mbox{ is of
    type }{\rm(p,pp)_s}\}=\emptyset$, then $r\leq\hat\pi$ on
  $[\sig_r,\hat\sig_{\hat\pi}]$, while otherwise, $r$ crosses $\hat\pi$ from
  left to right at time $\tau:=\inf(S)$, and one has $r\leq\hat\pi$ on
  $[\sig_r,\tau]$ and $\hat\pi\leq r$ on $[\tau,\hat\sig_{\hat\pi}]$.
\end{itemize}
Analogous statements hold for left-most paths crossing dual Brownian
net paths from right to left.
\el
{\bf Proof} To prove part~(a), imagine that $r(t)<\hat\pi(t)$ for some
$\sig_r<t\leq\hat\sig_{\hat\pi}$. Let $\hat r\in\hat\Wr$ be the left-most
element of $\hat\Ni(\hat\pi(t),t)$. Then, by \cite[Prop.~3.2~(c)]{SS06}, one
has $r\leq\hat r\leq\hat\pi$ on $[\sig_r,t]$, contradicting our assumption
that $\hat\pi(\sig_r)<r(\sig_r)$.

To prove part~(b), set
\bc
t&:=&\dis\sup\{s\in(\sig_r,\hat\sig_{\hat\pi}):r(s)<\hat\pi(s)\},\\[5pt]
t'&:=&\dis\inf\{s\in(\sig_r,\hat\sig_{\hat\pi}):\hat\pi(s)<r(s)\}.
\ec
Then $t<t'$ by what we have just proved and therefore $t=t'$ by
Lemma~\ref{L:netint}, hence $t$ is the unique crossing point of $r$ and
$\hat\pi$, $r\leq\hat\pi$ on $[\sig_r,t]$, and $\hat\pi\leq r$ on
$[t,\hat\sig_{\hat\pi}]$.

To see that $(r(t),t)$ is a separation point, by Proposition~\ref{P:sep}, it
suffices to show that there exists some $\hat l\in\hat\Wl$ that crosses $r$ at
time $t$. By \cite[Lemma~3.4~(b)]{SS06},
we can without loss of generality assume that $r$ starts from a deterministic
point. Choose $\hat t>t$ with $\hat\pi(\hat t)<r(\hat t)$ and $\hat t\in\Ti$,
where $\Ti\sub\R$ is some fixed, deterministic countable dense set. By
Lemma~\ref{L:nocro}, there exist unique $\hat t>\tau_1>\cdots>\tau_M>\sig_r$,
$\tau_0:=\hat t$, $\tau_{M+1}:=-\infty$, with $0\leq M<\infty$, and $\hat
l_0\in\hat\Wl(r(\tau_0),\tau_0),\ldots,\hat l_M\in\Wl(r(\tau_M),\tau_M)$, such
that $\hat l_i\leq r$ on $[\tau_{i+1}\vee\sig_r,\tau_i]$ for each $0\leq i\leq
M$, and $\tau_{i+1}$ is the crossing time of $\hat l_i$ and $r$ for each
$0\leq i<M$. We claim that $\tau_i=t$ for some $i=1,\ldots,M$.

By \cite[Prop.~1.8]{SS06}, $\hat\pi\leq\hat l_0$ on $(-\infty,\hat t]$, so
$M\geq 1$ and $\tau_1\geq t$. If $\tau_1=t$ we are done. Otherwise, we claim
that $\hat\pi\leq\hat l_1$ on $(-\infty,\tau_1]$. To see this, assume that
$\hat l_1(s)<\hat\pi(s)$ for some $s<\tau_1$. Then we can start a left-most
path $l$ between $\hat l_1$ and $\hat\pi$. By what we have proved in
part~(a) and \cite[Prop.~3.2~(c)]{SS06}, $l$ is contained by $\hat l_1$ and
$\hat\pi$, hence $l$ and $r$ form a wedge of $(\Wl,\Wr)$ which by the
characterization of $\hat\Ni$ using wedges (Theorem~\ref{T:netchar}~(b2))
cannot be entered by $\hat\pi$, which yields a contradiction. This shows
that $\hat\pi\leq\hat l_1$ on $(-\infty,\tau_1]$. It follows that $M\geq 2$
and $\tau_2\geq t$. Continuing this process, using the finiteness of $M$,
we see that $\tau_i=t$ for some $i=1,\ldots,M$.

Conversely, if $r(s)=\hat\pi(s)$ for some $s\in(\sig_r,\hat\sig_{\hat\pi})$
such that $z:=(r(s),s)$ is of type $\rm(p,pp)_s$, then by Lemma~\ref{L:NT}~(b),
the path $\hat\pi$ is contained between the left-most and right-most paths
starting at $z$ that are {\em not} continuations of incoming left-most and
right-most paths. It follows that any incoming right-most path $r$ at $z$ must
satisfy $\hat\pi\leq r$ on $[s,\hat\sig_{\hat\pi}]$. This shows that $t$ is
the first separation point that $r$ meets on $\hat\pi$.\qed

\noi
{\bf Proof of Lemma~\ref{L:extl}} We will prove the statement for reflected
right-most paths; the statement for left-most paths then follows by
symmetry. It follows from the image set property (see formula (\ref{improp})
or \cite[Prop.~1.13]{SS06}) and the local equicontinuity of the Brownian net
that $r_{z,\hat\pi}(s):=\sup\{\pi(s):\pi\in\Ni(z),\ \pi\leq\hat\pi\mbox{ on
}[t,\hat\sig_{\hat\pi}]\}$ $(s\geq t)$ defines a path $r_{z,\hat\pi}\in\Ni$. Put
\be\ba{r@{\,}l}\label{Fdef}
\dis\Fi:=\big\{s\in(t,\hat\sig_{\hat\pi}):&\dis
\exists r\in\Wr\mbox{ such that $r$ crosses $\hat\pi$ at }(\hat\pi(s),s)
\mbox{ and }\exists\pi\in\Ni(z)\\[5pt]
&\dis\mbox{such that }\pi(s)=\hat\pi(s),
\mbox{ and }\pi\leq\hat\pi\mbox{ on }[t,s]\big\}.
\ec
By Lemma~\ref{L:lpi}, $(\hat\pi(s),s)$ is a $(t,\hat\sig_{\hat\pi})$-relevant
separation point for each $s\in\Fi$, so by Lemma \ref{L:finrel}, $\Fi$ is a
locally finite subset of $(t,\hat\sig_{\hat\pi})$.

We claim that
\be\ba{l}\label{contin}
\dis\forall s\in\Fi\mbox{ and }\pi\in\Ni(z)\mbox{ s.t.\ }
\pi(s)=\hat\pi(s)\mbox{ and }\pi\leq\hat\pi\mbox{ on }[t,s]\\[5pt]
\dis\quad\exists\pi'\in\Ni(z)\mbox{ s.t.\ }\pi'=\pi\mbox{ on }[t,s]
\mbox{ and }\pi\leq\hat\pi\mbox{ on }[t,\hat\sig_{\hat\pi}].
\ec
To prove this, set $\pi_0:=\pi$ and $s_0:=s$ and observe that
$(\pi_0(s_0),s_0)$ is a separation point where some right-most path $r$
crosses $\hat\pi$. Let $r'$ be the outgoing right-most path at
$(\pi_0(s_0),s_0)$ that is not equivalent to $r$. By Lemma~\ref{L:NT}~(b),
$r'\leq\hat\pi$ on $[s_0,s_0+\eps]$ for some $\eps>0$. Let $\pi_1$ be the
concatenation of $\pi_0$ on $[t,s_0]$ with $r'$ on $[s_0,\infty]$. Since the
Brownian net is closed under hopping \cite[Prop.~1.4]{SS06}, using the
structure of separation points, it is not hard to see that $\pi_1\in\Ni$.
(Indeed, we may hop from $\pi_0$ onto the left-most path entering
$(\pi_0(s_0),s_0)$ at time $s$ and then hop onto $r$ at some time $s+\eps$ and
let $\eps\down 0$, using the closedness of $\Ni$.) The definition of $\Fi$ and
Lemma~\ref{L:lpi}~(b) imply that either $\pi_1\leq\hat\pi$ on
$[t,\hat\sig_{\hat\pi}]$ or there exists some $s_1\in\Fi$, $s_1>s_0$ such that
$r'$ crosses $\hat\pi$ at $s_1$. In that case, we can continue our
construction, leading to a sequence of paths $\pi_n$ and times $s_n\in\Fi$
such that $\pi_n\leq\hat\pi$ on $[t,s_n]$. Since $\Fi$ is locally finite,
either this process terminates after a finite number of steps, or
$s_n\up\hat\sig_{\hat\pi}$. In the latter case, using the compactness of
$\Ni$, any subsequential limit of the $\pi_n$ gives the desired path $\pi'$.

We next claim $\Fi\sub\Fi(r_{z,\hat\pi})$, where the latter is defined as in
(\ref{Ridef}). Indeed, if $s\in\Fi$, then $(\hat\pi(s),s)$ is of type
$\rm(p,pp)_s$ by Lemma~\ref{L:lpi}~(b). Moreover, we can find some
$\pi\in\Ni(z)$ such that $\pi(s)=\hat\pi(s)$ and $\pi\leq\hat\pi$ on
$[t,s]$. By (\ref{contin}), we can modify $\pi$ on $[s,\hat\sig_{\hat\pi}]$ so
that it stays on the left of $\hat\pi$, hence $\pi\leq r_{z,\hat\pi}$ by the
maximality of the latter, which implies that $r_{z,\hat\pi}(s)=\hat\pi(s)$
hence $s\in\Fi(r_{z,\hat\pi})$.

Set $\Fi':=\Fi\cup\{t,\hat\sig_{\hat\pi},\infty\}$ if $\hat\sig_{\hat\pi}$ is
a cluster point of $\Fi$ and $\Fi':=\Fi\cup\{t,\infty\}$ otherwise. Let
$s,u\in\Fi^-$ satisfy $s<u$ and $(s,u)\cap\Fi'=\emptyset$. Let $\Ti$ be some
fixed, deterministic countable dense subset of $\R$. By Lemma~\ref{L:netint}
we can choose $t_n\in\Ti$ such that $t_n\down s$ and
$r_{z,\hat\pi}(t_n)<\hat\pi(t_n)$. Choose $r_n\in\Wr(r_{z,\hat\pi}(t_n),t_n)$
such that $r_{z,\hat\pi}\leq r_n$ on $[t_n,\infty]$. By
\cite[Lemma~8.3]{SS06}, the concatenation of $r_{z,\hat\pi}$ on $[t,t_n]$ with
$r_n$ on $[t_n,\infty]$ is a path in $\Ni$, hence $r_n$ can cross $\hat\pi$
only at times in $\Fi$, and therefore $r_n\leq\hat\pi$ on $[t_n,u]$. Using the
compactness of $\Wr$, let $r\in\Wr$ be any subsequential limit of the
$r_n$. Then $r_{z,\hat\pi}\leq r$ on $[s,\infty]$, $r\leq\hat\pi$ on $[s,u]$,
and, since $\Ni$ is closed, the concatenation of $r_{z,\hat\pi}$ on $[t,s]$
and $r$ on $[s,\infty]$ is a path in $\Ni$.

If $u<\hat\sig_{\hat\pi}$, then $r_{z,\hat\pi}(u)\leq r(u)\leq\hat\pi(u)$
while $r_{z,\hat\pi}(u)=\hat\pi(u)$ by the fact that
$\Fi\sub\Fi(r_{z,\hat\pi})$. In this case, since $(\hat\pi(u),u)$ is a
separation point, by Lemma~\ref{L:lpi}~(b), the path $r$ crosses $\hat\pi$ at
$u$, $\hat\pi\leq r$ on $[u,\hat\sig_{\hat\pi}]$, and
$\inf\{s'>s:\hat\pi(s')<r(s')\}=u$. Let $\pi$ denote the concatenation of
$r_{z,\hat\pi}$ on $[t,s]$ and $r$ on $[s,\infty]$. By (\ref{contin}), if
$u<\hat\sig_{\hat\pi}$, then we can modify $\pi$ on $[u,\infty]$ such that it
stays on the left of $\hat\pi$. Therefore, whether $u<\hat\sig_{\hat\pi}$ or
not, by the maximality of $r_{z,\hat\pi}$, we have $r_{z,\hat\pi}=r$ on
$[s,u]$.

To complete our proof we must show that $\Fi\supset\Fi(r_{z,\hat\pi})$. To see
this, observe that if $s\in\Fi(r_{z,\hat\pi})$, then, since $r_{z,\hat\pi}$ is
a concatenation of right-most paths, by Lemma~\ref{L:lpi}~(b), some right-most
path crosses $\hat\pi$ at $s$, hence $s\in\Fi$.\qed

\section{Incoming paths and the image set}\label{S:image}

\subsection{Maximal $T$-meshes}\label{S:image1}

Let $\Ni_T:=\{\pi\in\Ni:\sig_\pi=T\}$ denote the set of paths in the
Brownian net starting at a given time $T\in[-\infty,\infty]$ and let
$N_T$ be its image set, defined in (\ref{NT}). We call $N_T$ the {\em
image set of the Brownian net started at time $T$}. In the present
section, we will identify the connected components of the complement
of $N_T$ relative to $\{(x,t)\in\Rc:t\geq T\}$.

The next lemma is just a simple observation.
\bl\label{L:bex}{\bf[Dual paths exit meshes through the bottom point]}\\
If $M(r,l)$ is a mesh with bottom point $z=(x,t)$ and $\hat\pi\in\hat\Ni$
starts in $M(r,l)$, then $r(s)\leq\hat\pi(s)\leq l(s)$ for all
$s\in[t,\hat\sig_{\hat\pi}]$.
\el
{\bf Proof} Immediate from Lemma~\ref{L:lpi}~(a).\qed

We will need a concept that is slightly stronger than that of a
mesh.
\bd\label{D:starmesh}{\bf[$\ast$-meshes]}\\
A mesh $M(r,l)$ with bottom point $z=(x,t)$ is called a $*$-mesh if
there exist $\hat r\in\hat\Wr$ and $\hat l\in\hat\Wl$ with $\hat
r\sim^z_{\rm in}\hat l$, such that $r$ is the right-most element of
$\Wr(z)$ that passes on the left of $\hat r$ and $l$ is the
left-most element of $\Wl(z)$ that passes on the right of $\hat l$.
\ed

\bl\label{L:star}{\bf[Characterization of $*$-meshes]}\\
Almost surely for all $z=(x,t)\in\R^2$ and $\hat r\in\hat\Wr,\hat l\in\hat\Wl$
such that $\hat r\sim^z_{\rm in}\hat l$, the set
\be\label{Mhatlr}
M_z(\hat r,\hat l):=\big\{z'\in\R^2:\forall\hat\pi\in\hat\Ni(z')\ \exists\eps>0
\mbox{ s.t.\ }\hat r\leq\hat\pi\leq\hat l\mbox{ on }[t,t+\eps]\big\}
\ee
is a $*$-mesh with bottom point $z$. Conversely, each $*$-mesh is of this form.
\el
{\bf Remark} A simpler characterization of $\ast$-meshes (but one that is
harder to prove) is given in Lemma~\ref{L:starchar} below. Once
Theorem~\ref{T:classnet} is proved, it will turn out that if a mesh is not a
$*$-mesh, then its bottom point must be of type $\rm(o, ppp)$. (See
Figure~\ref{fig:LRwebmult}.)\med

\noi
{\bf Proof of Lemma \ref{L:star}} Let $z=(x,t)\in\R^2$ and $\hat
r\in\hat\Wr,\hat l\in\hat\Wl$ satisfy $\hat r\sim^z_{\rm in}\hat l$. Let $r$
be the right-most element of $\Wr(z)$ that passes on the left of $\hat r$ and
let $l$ be the left-most element of $\Wl(z)$ that passes on the right of $\hat
l$. Then, obviously, $M(r,l)$ is a $*$-mesh, and each $*$-mesh is of this
form. We claim that $M(r,l)=M_z(\hat r,\hat l)$.

To see that $M(r,l)\supset M_z(\hat r,\hat l)$, note that if
$z'=(x',t')\not\in M(r,l)$ and $t'>t$, then either there exists an
$\hat r\in\hat\Wr(z')$ that stays on the left of $r$, or there exists
an $\hat l\in\hat\Wl(z')$ that stays on the right of $l$; in either
case, $z'\not\in M_z(\hat r,\hat l)$.

To see that $M(r,l)\sub M_z(\hat r,\hat l)$, assume that
$z'=(x',t')\in M(r,l)$. Since each path in $\hat\Ni(z')$ is contained
by the left-most and right-most dual paths starting in $z'$, it
suffices to show that each $\hat r'\in\hat\Wr(z')$ satisfies $\hat
r\leq\hat r'$ on $[t,t+\eps]$ for some $\eps>0$ and each $\hat
l'\in\hat\Wl(z')$ satisfies $\hat l'\leq\hat l$ on $[t,t+\eps']$ for
some $\eps'>0$. By symmetry, it suffices to prove the statement for
$\hat r'$. So imagine that $\hat r'<\hat r$ on $(t,t')$. Then there
exists an $r'\in\Wr(z)$ that stays between $\hat r'$ and $\hat r$,
contradicting the fact that $r$ is the right-most element of $\Wr(z)$
that passes on the left of $\hat r$.\qed

\bd\label{D:tmaxmesh}{\bf[Maximal $T$-meshes]}\\
For a given $T\in[-\infty,\infty)$, we call a mesh $M(r,l)$ with
bottom point $z=(x,t)$ a $T$-mesh if $M(r,l)$ is a $*$-mesh and $t\geq
T$, and a maximal $T$-mesh if it is not contained in any other
$T$-mesh. A maximal $T$-mesh with $T=-\infty$ is called a mesh of the
backbone of the Brownian net.
\ed
\bp\label{P:tmaxmesh}{\bf[Properties of maximal $T$-meshes]}\\
Almost surely for each $T\in[-\infty,\infty)$:
\begin{itemize}
\item[\rm(a)] A set of the form (\ref{Mhatlr}), with $z=(x,t)$, is a
maximal $T$-mesh if and only if $t=T$ or if $t>T$ and $\hat r<\hat l$
on $(T,t)$.
\item[\rm(b)] The maximal $T$-meshes are mutually disjoint, and their
union is the set $\{(x,t)\in\Rc:t\geq T\}\beh N_T$, where $N_T$ is
defined in (\ref{NT}).
\end{itemize}
\ep
{\bf Proof} Let $M(r,l)$ and $M(r',l')$ be $\ast$-meshes with bottom
points $z=(x,t)$ and $z'=(x',t')$, respectively, with $t\geq t'\geq
T$. Then either $M(r,l)$ and $M(r',l')$ are disjoint, or there exists
a $z''\in M(r,l)\cap M(r',l')$. In the latter case, by
Lemma~\ref{L:star}, any $\hat r\in\hat\Wr(z'')$ and $\hat
l\in\hat\Wl(z'')$ must pass through $z$ and $z'$ (in this order) and
we have $M(r,l)=M_z(\hat r,\hat l)$ and $M(r',l')=M_{z'}(\hat r,\hat
l)$. It follows that $M_z(\hat r,\hat l)\sub M_{z'}(\hat r,\hat l)$,
where the inclusion is strict if and only if $t>t'$. Moreover, a set
of the form (\ref{Mhatlr}), with $z=(x,t)$, is a maximal $T$-mesh if
and only if there exists no $z'=(x',t)$ with $t'\in[T,t)$ such that
$\hat r$ and $\hat l$ are equivalent incoming paths at $z'$. This
proves part~(a).

We have just proved that $T$-meshes are either disjoint or one is
contained in the other, so maximal $T$-meshes must be mutually
disjoint. Let $O_T:=\{(x,t)\in\Rc:t\geq T\}\beh N_T$. It is easy to
see that $O_T\sub\R\times(T,\infty)$. Consider a point
$z=(x,t)\in\R^2$ with $t>T$. If $z\not\in O_T$, then by
Lemma~\ref{L:NT} there exist $\hat r\in\hat\Wr(z)$ and $\hat
l\in\hat\Wl(z)$ with the property that there does not exist a
$z'=(x',t')$ with $t'\geq T$ such that $\hat r\sim^{z'}_{\rm in}\hat
l$, hence by Lemma~\ref{L:star}, $z$ is not contained in any
$T$-mesh. On the other hand, if $z\in O_T$, then by
Lemma~\ref{L:NT}~(b) and the nature of convergence of paths in the
Brownian web (see \cite[Lemma~3.4~(a)]{SS06}), there exist $\hat
r\in\hat\Wr$ and $\hat l\in\hat\Wl$ starting from points $(x_-,t)$ and
$(x_+,t)$, respectively, with $x_-<x<x_+$, such that $\hat r(s)=\hat
l(s)$ for some $s\in(T,t)$. Now, setting $u:=\inf\{s\in(T,t):\hat
r(s)=\hat l(s)\}$ and $z':=(\hat r(u),u)$, by Lemma~\ref{L:star},
$M_{z'}(\hat r,\hat l)$ is a maximal $T$-mesh that contains $z$.\qed

\subsection{Reversibility}\label{S:image2}

Recall from (\ref{NT}) the definition of the image set $N_T$ of the
Brownian net started at time $T$. It follows from
\cite[Prop.~1.15]{SS06} that the law of $N_{-\infty}$ is symmetric
with respect to time reversal. In the present section, we extend this
property to $T>-\infty$ by showing that locally on
$\R\times(T,\infty)$, the law of $N_T$ is absolutely continuous with
respect to its time-reversed counterpart. This is a useful property,
since it allows us to conclude that certain properties that hold a.s.\
in the forward picture also hold a.s.\ in the time-reversed
picture. For example, meeting and separation points have a similar
structure, related by time reversal. (Note that this form of
time-reversal is different from, and should not be confused with, the
dual Brownian net.)

We write $\mu\ll\nu$ when a measure $\mu$ is absolutely
continuous with respect to another measure $\nu$, and $\mu\sim\nu$ if $\mu$
and $\nu$ are equivalent, i.e., $\mu\ll\nu$ and $\nu\ll\mu$.
\bp\label{P:rev}{\bf[Local reversibility]}\\
Let $-\infty<S,T<\infty$ and let $N_T$ be the image set of the
Brownian net started at time $T$. Define $\Ri_S:\R^2\to\R^2$ by
$\Ri_S(x,t):=(x,S-t)$. Let $K\sub\R^2$ be a compact set such that
$K,\Ri_S(K)\sub\R\times(T,\infty)$. Then
\be\label{rev}
\P[N_T\cap\Ri_S(K)\in\cdot\,]\sim\P[\Ri_S(N_T\cap K)\in\cdot\,].
\ee
\ep
{\bf Proof} By the reversibility of the backbone, it suffices to prove that
\be\label{backcomp}
\P[N_T\cap K\in\cdot\,]\sim\P[N_{-\infty}\cap K\in\cdot\,].
\ee
Choose some $T<s<\min\{t:(x,t)\in K\}$. By \cite[Prop.~1.12]{SS06}, the set
\be
\hat\xi^K_s:=\{\hat\pi(s):\hat\pi\in\hat\Ni(K)\}
\ee
is a.s.\ a finite subset of $\R$, say
\be
\hat\xi^K_s=\{X_1,\ldots,X_M\}\quad\mbox{with}\quad X_1<\cdots<X_M.
\ee
For $U=T,-\infty$, let us write
\be
\xi^{(U)}_s:=\{\pi(s):\pi\in\Ni,\ \sig_\pi=U\}.
\ee
By Lemma~\ref{L:NT} and the fact that $s$ is deterministic, for any
$z=(x,t)\in K$, one has $z\in N_U$ if and only if there exist $\hat
r\in\hat\Wr(z)$ and $\hat l\in\hat\Wl(z)$ such that $\hat r<\hat l$ on
$[s,t)$ and $\xi^{(U)}_s\cap(\hat r(s),\hat l(s))\neq\emptyset$. Thus, we
can write
\be
N_U\cap K=\bigcup_{i\in I_U}N_i,
\ee
where
\be\ba{r@{\,}l}
\dis N_i:=\big\{z=(x,t)\in K:
&\dis\exists\hat r\in\hat\Wr(z),\ \hat l\in\hat\Wl(z)\\[5pt]
&\dis\mbox{ s.t.\  $\hat r<\hat l$ on $[s,t)$ and }
\hat r(s)\leq X_i,\ \hat l(s)\geq X_{i+1}\big\}
\ec
and
\be
I_U:=\big\{i:1\leq i\leq M-1,\ \xi^{(U)}_s\cap(X_i,X_{i+1})\neq\emptyset\big\}.
\ee
It follows that
\be\label{cond1}
\P[N_U\cap K\in\cdot\,]
=\E\Big[\P\Big[\bigcup_{i\in I_U}N_i\in\cdot\,\,\Big|\,
\{X_1,\ldots,X_M\}\Big]\Big],
\ee
where
\be\ba{l}\label{cond2}
\dis\P\Big[\bigcup_{i\in I_U}N_i\in\cdot\,\,\Big|\,
\{X_1,\ldots,X_M\}\Big]\\[5pt]
\dis\quad=\sum_{\Ii\sub\{1,\ldots,M-1\}}
\P\Big[\bigcup_{i\in\Ii}N_i\in\cdot\,\,\Big|\,\{X_1,\ldots,X_M\}\Big]
\P\Big[I_U=\Ii\,\,\Big|\,\{X_1,\ldots,X_M\}\Big].
\ec
Here the sum ranges over all subsets $\Ii$ of $\{1,\ldots,M-1\}$. The
statement of the proposition now follows from (\ref{cond1}) and
(\ref{cond2}) by observing that
\be
\P\Big[I_U=\Ii\,\,\Big|\,\{X_1,\ldots,X_M\}\Big]>0\quad{\rm a.s.}
\ee
for all $\Ii\sub\{1,\ldots,M-1\}$ and $U=T,-\infty$.\qed

\subsection{Classification according to incoming paths}\label{S:image4}

In this section, we give a preliminary classification of points
in the Brownianb
 net that is entirely based on incoming paths. Note
that if there is an incoming path $\pi\in\Ni$ at a point $z=(x,t)$,
then $z\in N_T$ for some $T<t$, where $N_T$ is the image set of the
Brownian net started at time $T$, defined in (\ref{NT}). Therefore in
this section, our main task is to classify the special points
of $N_T$.

For a given $T\in[-\infty,\infty)$, let us say that a point
$z=(x,t)\in N_T\cap(\R\times(T,\infty))$ is {\em isolated from the
left} if
\be
\sup\{x'\in\R:(x',t)\in N_T,\ x'<x\}<x.
\ee
Points that are isolated from the right are defined similarly, with
the supremum replaced by an infimum and both inequality signs reversed.
\bl\label{L:iso}{\bf[Isolated points]}\\
A point $z=(x,t)\in N_T\cap(\R\times(T,\infty))$ is isolated from
the left if and only if there exists an incoming path $l\in\Wl$ at~$z$.
An analogous statement holds if $z$ is isolated from the right.
\el
{\bf Proof} By Proposition~\ref{P:tmaxmesh}~(b), if $z=(x,t)$ is
isolated from the left, then there exists a maximal $T$-mesh $M(r,l)$
with left and right boundary $r$ and $l$, bottom time strictly smaller
than $t$, and top time strictly larger than $t$, such that
$l(t)=x$. Conversely, if there exists an incoming path $l\in\Wl$ at
$z$, then by \cite[Lemma~6.5]{SS06}, there exists a mesh $M(r,l)$ with
bottom time in $(T,t)$ and top time in $(t,\infty)$, such that
$l(t)=x$. Therefore, by the characterization of the Brownian net with
meshes (see Theorem~\ref{T:netchar}~(b3)), $z$ is isolated from the
left.\qed

We now give a classification of points in $\R^2$ based on incoming
paths in the Brownian net. Recall Definition~\ref{D:IMS} of
intersection, meeting, and separation points.
\bd\label{D:imgpts}{\bf[Classification by incoming paths]}\\
We say that a point $z=(x,t)\in\R^2$ is of type
\begin{itemize}
\item[${\rm(C_o)}$] if there is no incoming $\pi\in\Ni$ at $z$,
\item[${\rm(C_n)}$] if there is an incoming $\pi\in\Ni$ at $z$, but there
is no incoming $\pi\in\Wl\cup\Wr$ at $z$,
\item[${\rm(C_l)}$] if there is an incoming $l\in\Wl$ at $z$, but there is
no incoming $r\in\Wr$ at $z$,
\item[${\rm(C_r)}$] if there is an incoming $r\in\Wr$ at $z$, but there is
no incoming $l\in\Wl$ at $z$,
\item[${\rm(C_s)}$] if $z$ is a separation point of some $l\in\Wl$ and
$r\in\Wr$,
\item[${\rm(C_m)}$] if $z$ is a meeting point of some $l\in\Wl$ and $r\in\Wr$,
\item[${\rm(C_p)}$] if there is an incoming $l\in\Wl$ at $z$ and an
incoming $r\in\Wr$ at $z$, and $z$ is not of type ${\rm(C_s)}$ or~${\rm(C_m)}$.
\end{itemize}
\ed

\noi
Note that by Lemma~\ref{L:NT}, for any $T<t$ such that $z\in N_T$,
points of types $\rm(C_l),(C_r)$, and $\rm(C_n)$ are either not
isolated from the left, or not isolated from the right, or both. In
view of this, we call these points {\em cluster points}. Points of the
types $\rm(C_l)$ and $\rm(C_r)$ are called {\em one-sided cluster points}
and points of type $\rm(C_n)$ {\em two-sided cluster points}. Proposition
\ref{P:imgclas} below shows that, among other things, cluster points are the
limits of nested sequences of excursions between left-most and
right-most paths.

The main result of this section is the following.
\bl\label{L:imgpts}{\bf[Classification by incoming paths]}
\begin{itemize}
\item[\rm(a)] Almost surely, each point in $\R^2$ is of exactly one of
the types ${\rm(C_o)}$, ${\rm(C_n)}$, ${\rm(C_l)}$, ${\rm(C_r)}$,
${\rm(C_s)}$, ${\rm(C_m)}$, and ${\rm(C_p)}$, and each of these types
occurs in $\R^2$.
\item[\rm(b)] For deterministic $t\in\R$, a.s.\ each point in
$\R\times\{t\}$ is of type ${\rm(C_o)}$ or ${\rm(C_p)}$, and both
these types occur.
\item[\rm(c)] Each deterministic $z\in\R^2$ is a.s.\ of type ${\rm(C_o)}$.
\end{itemize}
\el
{\bf Proof of Lemma~\ref{L:imgpts}~(b) and (c)}
If $t\in\R$ is deterministic, and $T<t$, then by
\cite[Prop~1.12]{SS06}, the set $N_T\cap(\R\times\{t\})$ is locally
finite. In particular, if there is an incoming path $\pi\in\Ni$ at a
point $z\in\R\times\{t\}$, then $z$ is isolated from the left and
right. Since each meeting or separation point of some $l\in\Wl$ and
$r\in\Wr$ is the meeting or separation point of some left-most and
right-most path chosen from a fixed, deterministic countable dense
set, and since paths started at deterministic starting points a.s.\ do
not meet or separate at deterministic times, $z$ must be of type
${\rm(C_p)}$.

If $z=(x,t)\in\R^2$ is deterministic, then by \cite[Prop~1.12]{SS06},
for each $n\geq 1$, a.s.\ there is no path $\pi\in\Ni$ with
$\sig_\pi=t-1/n$ and $\pi(t)=x$, hence $z$ must be of type
${\rm(C_o)}$.\qed

Before proving Lemma~\ref{L:imgpts}~(a), we first establish some basic
properties for each type of points in Definition \ref{D:imgpts}. We start with
a definition and a lemma.

By definition, we say that two paths $\pi,\pi'\in\Pi$ make an {\em excursion}
{f}rom each other {\em on a time interval} $(s,u)$ if $\sig_\pi,\sig_{\pi'}<s$
(note the strict inequality), $\pi(s)=\pi'(s)$, $\pi\neq\pi'$ on $(s,u)$, and
$\pi(u)=\pi'(u)$. The next lemma says that excursions between left-most and
right-most paths are rather numerous.
\bl\label{L:exc}{\bf[Excursions along left-most paths]}\\
Almost surely, for each $l\in\Wl$ and for each open set $O$ such that
$l\cap O\neq\emptyset$, there exists an $r\in\Wr$ such that $r$ makes
an excursion from $l$ during a time interval $(s,u)$, and
$\{(x,t): t\in [s,u], x\in [l(t), r(t)]\} \sub O$.
\el
{\bf Proof} Choose $s<u$ such that $\{l(t):t\in[s,u]\}\sub O$ and choose some
$t$ from a fixed, deterministic countable dense subset of $\R$ such that
$t\in(s,u)$. By Lemma~\ref{L:imgpts}~(b), there exists a unique incoming path
$r\in\Wr$ at $(l(t),t)$, and $r$ does not separate from $l$ at time $t$. Since
by \cite[Prop.~3.6~(b)]{SS06}, the set $\{v:r(v)=l(v)\}$ is nowhere dense, we
can find $u_1>s_1>u_2>s_2>\cdots$ such that $u_n\down t$ and $r$ makes an
excursion from $l$ during the time interval $(s_n,u_n)$ for each $n\geq 1$. By
the local equicontinuity of the Brownian net, we can choose $n$ large enough
such that $\{(x,t): t\in [s_n,u_n], x\in [l(t), r(t)]\} \sub O$.\qed

\begin{figure}[tp] 
\centering
\plaat{0.5cm}{\includegraphics[width=10cm]{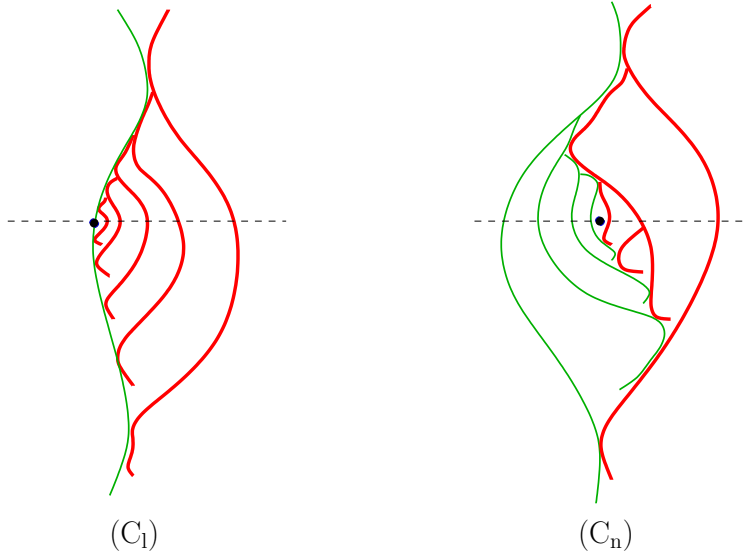}}
\caption{Nested excursions around cluster points.}
\label{fig:clusterpts}
\end{figure}

\bp\label{P:imgclas}{\bf[Structure of points with incoming paths]}
\begin{itemize}
\item[\rm(a)] If $z=(x,t)$ is of type ${\rm(C_l)}$, then there exist $l\in\Wl$
  and $r_n\in\Wr$ $(n\geq 1)$, such that $l(t)=x<r_n(t)$, each $r_n$ makes an
  excursion away from $l$ on a time interval $(s_n,u_n)\ni t$,
  $[s_n,u_n]\sub(s_{n-1},u_{n-1})$, $s_n\up t$, $u_n\down t$, and $r_n(t)\down
  x$. (See Figure \ref{fig:clusterpts}). By symmetry, an analogous statement
  holds for points of type ${\rm(C_r)}$.

\item[\rm(b)] If $z=(x,t)$ is of type ${\rm(C_n)}$, then there exist
  $l_1\in\Wl,r_2\in\Wr,l_3\in\Wl,\ldots$, such that $l_{2n+1}(t)<x<r_{2n}(t)$,
  each path ($l_n$ for $n$ odd, $r_n$ for $n$ even) in the sequence makes an
  excursion away from the previous path on a time interval $(s_n,u_n)\ni t$,
  $[s_n,u_n]\sub(s_{n-1},u_{n-1}]$, $s_n\up t$, $u_n\down t$, $l_{2n+1}(t)\up
  x$, and $r_{2n}(t)\down x$ (the monotonicity here need not be strict).

\item[\rm(c)] If $z=(x,t)$ is of type ${\rm(C_s)}$, then for each $T<t$ with
  $z\in N_T$, there exist maximal $T$-meshes $M(r,l)$ and $M(r',l')$ with
  bottom times strictly smaller than $t$ and top times strictly larger than
  $t$, and a maximal $T$-mesh $M(r'',l'')$ with bottom point $z$, such that
  $l\sim^z_{\rm in}r'$, $l\sim^z_{\rm out}r''$, and $l''\sim^z_{\rm out}r'$.

\item[\rm(d)] If $z=(x,t)$ is of type ${\rm(C_m)}$, then for each $T<t$ with
  $z\in N_T$, there exist maximal $T$-meshes $M(r,l)$ and $M(r',l')$ with
  bottom times strictly smaller than $t$ and top times strictly larger than
  $t$, and a maximal $T$-mesh $M(r'',l'')$ with top point $z$, such that
  $l\sim^z_{\rm in}r''$, $l''\sim^z_{\rm in}r'$, and $l\sim^z_{\rm out}r'$.

\item[\rm(e)] If $z=(x,t)$ is of type ${\rm(C_p)}$, then for each $T<t$ with
  $z\in N_T$, there exist maximal $T$-meshes $M(r,l)$ and $M(r',l')$ with
  bottom times strictly smaller than $t$ and top times strictly larger than
  $t$, such that $l\sim^z_{\rm in}r'$ and $l\sim^z_{\rm out}r'$.
\end{itemize}
\ep
{\bf Proof } (a): Let $l$ be an incoming left-most path at $z$ and choose
$T<t$ such that $l\sub N_T$. Choose $t_n$ from some fixed, deterministic
countable dense subset of $\R$ such that $t_n\up t$. By
Lemma~\ref{L:imgpts}~(b), for each $n$ there is a unique incoming path
$r_n\in\Wr$ at $(l(t_n),t_n)$. By assumption, $r_n$ does not pass through $z$,
hence $r_n$ makes an excursion away from $l$ on a time interval $(s_n,u_n)$
with $t_n\leq s_n<t<u_n\leq\infty$. Clearly $s_n\up t$ and $u_n\downarrow
u_\infty$ for some $u_\infty\geq t$. We claim that $u_\infty=t$. Indeed, if
$u_\infty>t$, then $(l(s_n),s_n)$ $(n\geq 1)$ are $(s_1,u_\infty)$-relevant
separation points, hence the latter are not locally finite on
$\R\times(s_1,u_\infty)$, contradicting Lemma~\ref{L:finrel}. By going to a
subsequence if necessary, we can asssure that $s_{n-1}<s_n$ and $u_n<u_{n-1}$.
The fact that $r_n(t)\down x$ follows from the local equicontinuity of the
Brownian net.

(b): Let $\pi\in\Ni$ be an incoming path at $z$ and choose $T<t$ such
that $\pi\sub N_T$. By Lemma~\ref{L:NT}~(a), there exist $\hat
r\in\hat\Wr(z)$ and $\hat l\in\hat\Wl(z)$ such that $\hat r<\hat l$ on
$(T,t)$ and $\hat r\leq\pi\leq\hat l$ on $[T,t]$. Choose an arbitrary
$l_1\in\Wl$ starting from $z_1=(x_1,s_1)$ with $s_1\in(T,t)$ and
$x_1\in(\hat r(s_1),\hat l(s_1))$. Since there is no incoming
left-most path at $z$, and since $l_1$ cannot cross $\hat l$, the path
$l_1$ must cross $\hat r$ at some time $s_2<t$. By
Proposition~\ref{P:sep}, there exists $r_2\in\Wr$ such that
$(l_1(s_2),s_2)$ is a separation point of $l_1$ and $r_2$, and $r_2$
lies on the right of $\hat r$. Since there is no incoming right-most
path at $z$, and since $r_2$ cannot cross $\hat r$, the path $r_2$
must cross $\hat l$ at some time $s_3<t$, at which a path
$l_3\in\Wl$ separates from $r_2$, and so on. Repeating this procedure
gives a sequence of paths $l_1\in\Wl,r_2\in\Wr,l_3\in\Wl,\ldots$ such
that the $n$-th path in the sequence separates from the $(n-1)$-th
path at a time $s_n\in(s_{n-1},t)$, and we have $l_1(t)\leq
l_3(t)\leq\cdots<x<\cdots\leq r_4(t)\leq r_2(t)$. By the a.s.\ local
equicontinuity of $\Wl\cup\Wr$ and $\hat\Wl\cup\hat\Wr$, it is clear
that $s_n\uparrow t$, $l_{2n+1}(t) \uparrow x$ and $r_{2n}(t)\downarrow x$.
Let $u_n:=\inf\{s\in(s_n,\infty):l_{n-1}=r_n\}$ if $n$ is even
and $u_n:=\inf\{s\in(s_n,\infty):r_{n-1}=l_n\}$ if $n$ is odd.
Then $u_n\down u_\infty\geq t$. The same argument
as in the proof of part~(a) shows that $u_\infty=t$.

To prepare for parts (c)--(e), we note that if there exist incoming
paths $l\in\Wl$ and $r\in\Wr$ at a point $z\in\R^2$ and $T<t$ is such
that $z\in N_T$, then by Lemma~\ref{L:iso}, $z$ is isolated from the
left and from the right, hence by Proposition~\ref{P:tmaxmesh}~(b),
there exist maximal $T$-meshes $M(r,l)$ and $M(r',l')$ with bottom
times strictly smaller than $t$ and top times strictly larger than
$t$, such that $l(t)=x=r'(t)$.  We now prove (c)--(e).

(c): If $z$ is a separation point of some left-most and right-most
path, then by Proposition~\ref{P:sep}, $z$ must be a separation point
of the paths $l$ and $r'$ mentioned above, and there exist paths
$r''\in\Wr(z)$ and $l''\in\Wl(z)$ such that $l\sim^z_{\rm out}r''$,
and $l''\sim^z_{\rm out}r'$, and $z$ is the bottom point of the mesh
$M(r'',l'')$. Note that $M(r'',l'')$ is a maximal $T$-mesh by
Definition~\ref{D:tmaxmesh}.

(d): If $z$ is a meeting point of some left-most and right-most path,
then $z$ must be a meeting point of the paths $l$ and $r'$ mentioned
above. This means that $z$ is a point where two maximal $T$-meshes
meet. Recall from Proposition~\ref{P:tmaxmesh} that the maximal
$T$-meshes are the connected components of the open set
$(\R\times(T,\infty))\beh N_T$. Thus, reversing time for $N_T$ in some
compact environment of $z$ makes $z$ into a separation point of two
maximal $T$-meshes. Therefore, by local reversibility
(Proposition~\ref{P:rev}) and what we have just proved about points of
type $\rm(C_s)$, there exists a maximal $T$-mesh $M(r'',l'')$ with top
point $z$, such that $l\sim^z_{\rm in}r''$, $l''\sim^z_{\rm in}r'$,
and $l\sim^z_{\rm out}r'$.

(e): If $z$ is not a separation or meeting point of any left-most and
right-most path, then the paths $l$ and $r'$ mentioned above must
satisfy $l\sim^z_{\rm in}r'$ and $l\sim^z_{\rm out}r'$.\qed
\bigskip

\noi
{\bf Proof of Lemma~\ref{L:imgpts}~(a)} The statement that each point in
$\R^2$ belongs to exactly one of the types ${\rm(C_o)}$, ${\rm(C_n)}$,
${\rm(C_l)}$, ${\rm(C_r)}$, ${\rm(C_s)}$, ${\rm(C_m)}$, and
${\rm(C_p)}$ is entirely self-evident, except that we have to show
that a point cannot at the same time be a separation point of some
$l\in\Wl$ and $r\in\Wr$, and a meeting point of some (possibly
different) $l'\in\Wl$ and $r'\in\Wr$. This however follows from
Proposition~\ref{P:sep}.

It follows from parts~(b) and (c) of the lemma (which have already been
proved) that points of the types ${\rm(C_o)}$ and ${\rm(C_p)}$
occur. Obviously, points of types ${\rm(C_s)}$ and ${\rm(C_m)}$ occur as
well. To prove the existence of cluster points, it suffices to establish the
existence of nested sequences of excursions, which follows from
Lemma~\ref{L:exc}.\qed

Proposition~\ref{P:imgclas} yields a useful consequence.
\bl\label{L:meetsep}{\bf[Separation and meeting points]}\\
If $z\in\R^2$ is a separation (resp.\ meeting) point of two paths
$\pi,\pi'\in\Ni$, then $z$ is of type $\rm(p,pp)_s$ (resp.\
$\rm(pp,p)$).
\el
{\bf Proof} We start by showing that if $z\in\R^2$ is a separation
(resp.\ meeting) point of two paths $\pi,\pi'\in\Ni$, then $z$ is of
type $\rm(C_s)$ (resp.\ $\rm(C_m)$). Thus, we need to show that a.s.\
for any $z\in\R^2$, if $\pi,\pi'\in\Ni_T$ are incoming paths at $z$,
one has $\pi\sim^z_{\rm in}\pi'$ if $z$ is not a meeting point, and
$\pi\sim^z_{\rm in}\pi'$ if $z$ is not a separation point.

If $z$ is not a cluster point, these statements follow from the
configuration of maximal $T$-meshes around $z$ as described in
Proposition~\ref{P:imgclas}~(c)--(e).

If $z$ is a two-sided cluster point, then any path $\pi\in\Ni$ must
pass through the top points of the nested excursions around $z$, hence
all paths in $\Ni(z)$ are equivalent as outgoing paths at $z$. If $z$
is a one-sided cluster point, then by \cite[Prop.~1.8]{SS06},
$l\leq\pi$ on $[t,\infty)$ for all incoming net paths at $z$, hence
all incoming net paths must pass through the top points of the nested
excursions and therefore be equivalent as outgoing paths.

Our previous argument shows that at a cluster point $z=(x,t)$, for any
$T<t$ such that $z\in N_T$, all paths $\pi\in\Pi$ such that $\pi(t)=x$
and $\pi\sub N_T$ are equivalent as outgoing paths at $z$. By local
reversibility (Proposition~\ref{P:rev}), it follows that all paths
$\pi\in\Pi$ such that $\pi(t)=x$ and $\pi\sub N_T$ are equivalent as
ingoing paths at $z$. (Note that reversing time in $N_T$ does not
change the fact that $z$ is a cluster point.)

This completes the proof that if $z\in\R^2$ is a separation (resp.\
meeting) point of two paths $\pi,\pi'\in\Ni$, then $z$ is of type
$\rm(C_s)$ (resp.\ $\rm(C_m)$). If $z$ is of type $\rm(C_s)$, then by
Proposition~\ref{P:sep}, $z$ is of type $\rm(p,pp)_s$. If $z$ is of
type $\rm(C_m)$, then by Proposition~\ref{P:imgclas}~(d), there are
exactly two incoming left-right pairs at $z$, and there is at least
one outgoing left-right pair at $z$. Therefore, $z$ must be of type
$(2,1)$ in both $\Wl$ and $\Wr$, hence there are no other outgoing
paths in $\Wl\cup\Wr$ at $z$, so $z$ is of type $\rm(pp,p)$.\qed

\noi
{\bf Remark} Lemma~\ref{L:meetsep} shows in particular that any
meeting point of two paths $\pi,\pi'\in\Wl\cup\Wr$ is of type
$\rm(pp,p)$. This fact can be proved by more elementary methods as
well. Consider the Markov process $(L,R,L')$ given by the
unique weak solutions to the SDE
\bc\label{LRL}
\dis\di L_t&=&\dis 1_{\txt\{L_t<R_t\}}\di B^{\rm l}_t
+1_{\txt\{L_t=R_t\}}\di B^{\rm s}_t-\di t,\\[5pt]
\dis\di R_t&=&\dis 1_{\txt\{L_t<R_t\}}\di B^{\rm r}_t
+1_{\txt\{L_t=R_t\}}\di B^{\rm s}_t+\di t,\\[5pt]
\dis\di L'_t&=&\dis\di B^{{\rm l}'}_t-\di t,
\ec
where $B^{\rm l},B^{\rm r},B^{\rm s}$, and $B^{{\rm l}'}$ are
independent Brownian motions, and we require that $L_t\leq R_t$ for
all $t\geq 0$. Set
\be\label{lrltau}
\tau=\inf\{t\geq 0:R_t=L'_t\}.
\ee
Then the claim follows from the fact that the process started in
$(L_0,R_0,L'_0)=(0,0,\eps)$ satisfies
\be\label{aim2}
\lim_{\eps\to 0}\P^{(0,0,\eps)}\big[L_\tau=R_\tau\big]=1,
\ee
which can be shown by a submartingale argument. Since this proof is of interest
on its own, we give it in Appendix~\ref{A:extra2}.\med

The next lemma is a simple consequence of Lemma~\ref{L:meetsep}.
\bl\label{L:starchar}{\bf[Characterization of $\ast$-meshes]}\\
A mesh $M(r,l)$ with bottom point $z=(x,t)$ is a $*$-mesh if and only if there
exists no $\pi\in\Wl(z)\cup\Wr(z)$, $\pi\neq l,r$, such that $r\leq\pi\leq l$
on $[t,t+\eps]$ for some $\eps>0$.
\el
{\bf Proof} If $M(r,l)$ is a mesh with bottom point $z=(x,t)$ and there exists
some $l\neq l'\in\Wl(z)$ such that $r\leq l'\leq l$ on $[t,t+\eps]$ for some
$\eps>0$, then $r<l'$ on $[t,t+\eps]$ by \cite[Prop.~3.6~(a)]{SS06}, hence by
Lemma~\ref{L:bex} we can find $\hat r\in\hat\Wr$ and $\hat l\in\hat\Wl$ such
that $r\leq\hat r\leq\hat l\leq l'$ on $[t,t+\eps]$, which implies that
$M(r,l)$ is not a $\ast$-mesh. By symmetry, the same is true if there exists
some $r\neq r'\in\Wr(z)$ such that $r\leq r'\leq l$ on $[t,t+\eps]$.

For any mesh $M(r,l)$, by Lemma~\ref{L:bex}, we can find $\hat r\in\hat\Wr$
and $\hat l\in\hat\Wl$ such that $r\leq\hat r\leq\hat l\leq l'$ on
$[t,t+\eps]$ for some $\eps>0$. We claim that $\hat r\sim^z_{\rm in}\hat l$ if
$M(r,l)$ satisfies the assumptions of Lemma~\ref{L:starchar}; from this, it
then follows that $M(r,l)$ is a $\ast$-mesh. To prove our claim, assume that
$\hat r\not\sim^z_{\rm in}\hat l$. Then, by Lemma~\ref{L:meetsep}, there exist
$\hat l'\in\hat\Wl$ and $\hat r'\in\hat\Wr$ such that $\hat r\leq\hat l'<\hat
r'\leq\hat l$ on $[t,t+\eps']$ for some $\eps'>0$. By
\cite[Prop~3.6~(d)]{SS06}, this implies that there exist $l\in\Wl(z)$ such
that $\hat l'\leq l\leq\hat r'$ on $[t,t+\eps']$, hence $M(r,l)$ does not
satisfy the assumptions of Lemma~\ref{L:starchar}.\qed

\subsection{Special times}\label{S:image5}

For any closed set $K\sub\R^2$, set
\be\label{xiCt}
\xi^K_t:=\{\pi(t):\pi\in\Ni(K),\ \sig_\pi\leq t\}.
\ee
It has been proved in \cite[Theorem~1.11]{SS06} that for any closed
$A\sub\R$, the process $(\xi^{A\times\{0\}}_t)_{t\geq 0}$ is a Markov
process taking values in the space of closed subsets of $\R$. It was shown in
\cite[Prop~1.12]{SS06} that $\xi^{A\times\{0\}}_t$ is a.s.\ a locally
finite point set for each deterministic $t>0$. It was claimed without
proof there that there exists a dense set of times $t>0$ such that
$\xi^{A\times\{0\}}_t$ is not locally finite. Indeed, with the help of
Lemma~\ref{L:exc}, we can prove the following result.
\bp\label{P:notiso}{\bf[No isolated points]}\\
Almost surely, there exists a dense set $\Ti\sub(0,\infty)$ such that
for each $t\in\Ti$ and for each closed $A\sub\R$, the set
$\xi^{A\times\{0\}}_t$ contains no isolated points.
\ep
{\bf Proof} We claim that it suffices to prove the statement for
$A=\R$. To see this, suppose that $x\in\xi^{A\times\{0\}}_t$ is not isolated
{f}rom the left in $\xi^{\R\times\{0\}}_t$ for some $t>0$. It follows from the
characterization of the Brownian net using meshes (see
Theorem~\ref{T:netchar}~(b3)) that the pointwise infimum
$\pi:=\inf\{\pi'\in\Ni(A\times\{0\}):\pi'(t)=x\}$ defines a path
$\pi\in\Ni(A\times\{0\})$. Let $l$ be the left-most element of
$\Wl(\pi(0),0)$. By Lemma~\ref{L:iso}, there is no incoming left-most
path at $(x,t)$, hence we must have $l(t)<x$. Since $x$ is not
isolated from the left in $\xi^{\R\times\{0\}}_t$, there are
$\pi_n\in\Ni$ starting at time $0$ such that $\pi_n(t)\in(l(t),x)$ for
each $n$, and $\pi_n(t)\up x$. Now each $\pi_n$ must cross either $l$
or $\pi$, so by the fact that the Brownian net is closed under hopping
(\cite[Prop.~1.4]{SS06}), $x$ is not isolated from the left in
$\xi^{A\times\{0\}}_t$.

To prove the proposition for $A=\R$, we claim that for each $0<s<u$
and for each $n\geq 1$, we can find $s\leq s'<u'\leq u$ such that
\be\label{notis}
(x-\ffrac{1}{n},x+\ffrac{1}{n})\cap\xi^{\R\times\{0\}}_t\neq\{x\}\quad
\forall t\in(s',u'),\ x\in(-n,n)\cap\xi^{\R\times\{0\}}_t.
\ee
To show this, we proceed as follows. If (\ref{notis}) holds for $s'=s$
and $u'=u$ we are done. Otherwise, we can find some $t_1\in(s,u)$ and
$x_1\in(-n,n)\cap\xi^{\R\times\{0\}}_{t_1}$ such that
$(x_1-\ffrac{1}{n},x_1+\ffrac{1}{n})\cap\xi^{\R\times\{0\}}_{t_1}=\{x_1\}$. In
particular, $x_1$ is an isolated point so there is an incoming
$l_1\in\Wl$ at $(x_1,t_1)$, hence by Lemma~\ref{L:exc} we can find an
$r_1\in\Wr$ such that $r_1$ makes an excursion from $l_1$ during a time
interval $(s_1,u_1)$ with $s<s_1<u_1<u$, with the additional property
that $x_1-\ffrac{1}{2n}\leq l_1<r_1\leq x_1+\ffrac{1}{2n}$ on
$(s_1,u_1)$. Now either we are done, or there exists some
$t_2\in(s_1,u_1)$ and $x_2\in(-n,n)\cap\xi^{\R\times\{0\}}_{t_2}$ such that
$(x_2-\ffrac{1}{n},x_2+\ffrac{1}{n})\cap\xi^{\R\times\{0\}}_{t_2}=\{x_2\}$.
In this case, we can find $l_2\in\Wl$ and $r_2\in\Wr$ making an excursion
during a time interval $(s_2,u_2)$, with the property that $l_2$ and $r_2$
stay in $[x_2-\ffrac{1}{2n},x_2+\ffrac{1}{2n}]$. We iterate this process if
necessary. Since $x_{m+1}$ must be at least a distance $\ffrac{1}{2n}$ from
each of the points $x_1,\ldots,x_m$, this process terminates after a finite
number of steps, proving our claim.

By what we have just proved, for any $0<s<u$, we can find $s\leq
s_1\leq s_2\leq\cdots\leq u_2\leq u_1\leq u$ such that
\be
(x-\ffrac{1}{n},x+\ffrac{1}{n})\cap\xi^{\R\times\{0\}}_t\neq\{x\}\quad
\forall n\geq 1,\ t\in(s_n,u_n),\ x\in(-n,n)\cap\xi^{\R\times\{0\}}_t.
\ee
Necessarily $\bigcap_n (s_n, u_n) =\{t\}$ for some $t\in\R$, and we
conclude that $\xi^{\R\times\{0\}}_t$ contains no isolated points.\qed

\section{Excursions}\label{S:excur}

\subsection{Excursions between forward and dual paths}\label{S:excur1}

\begin{figure}[tp] 
\centering
\plaat{5cm}
{\includegraphics[width=3.8cm]{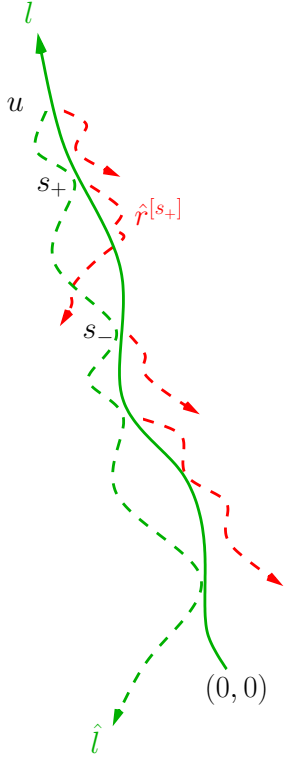}}
\caption{Excursions between a left-most and dual left-most path, with dual
  right-most paths starting at the top of each excursion.}
\label{fig:excur}
\end{figure}

Excursions between left-most and right-most paths have already been
studied briefly in Section~\ref{S:image4}. In this section, we
study them in more detail. In particular, in order to prove the
existence of points from groups (4) and (5) of
Theorem~\ref{T:classnet}, we will need to prove that, for a given
left-most path $l$ and a dual left-most path $\hat l$ that hits $l$ from
the left, there exist nested sequences of excursions of right-most
paths away from $l$, such that each excursion interval contains an
intersection point of $l$ and $\hat l$. As a first step towards
proving this, we will study excursions between $l$ and $\hat l$.

By definition, we say that a forward path $\pi\in\Pi$ and a backward
path $\hat\pi\in\hat\Pi$ make an {\em excursion} from each other {\em
on a time interval} $(s_-,s_+)$ if $\sig_\pi<s_-<s_+<\hat\sig_{\hat\pi}$
(note the strict inequalities), $\pi(s_-)=\hat\pi(s_-)$, $\pi\neq\pi'$ on
$(s_-,s_+)$, and $\pi(s_+)=\hat\pi(s_+)$. We write
\bc\label{Iidef}
\dis\Ii(\pi,\hat\pi)&:=&\dis\big\{t\in(\sig_\pi,\hat\sig_{\hat\pi}):
\pi(t)=\hat\pi(t)\big\},\\[5pt]
\dis\Ii_+(\pi,\hat\pi)&:=&\dis\big\{s_+:\mbox{$\pi$ and
$\hat\pi$ make an excursion}\\
&&\dis\phantom{\big\{s_+:}
\mbox{from each other on $(s_-,s_+)$ for some
$\sig_\pi<s_-<s_+<\hat\sig_{\hat\pi}$}\big\}.
\ec
Our next proposition is the main result of this section. Below,
for given $s_+\in\Ii_+(l,\hat l)$, it is understood that $s_-$ denotes
the unique time such that $l$ and $\hat l$ make an excursion from each
other on $(s_-,s_+)$. See Figure~\ref{fig:excur}.
\bp\label{P:exent}{\bf[Excursions entered by a dual path]}\\
Almost surely, for each $l\in\Wl$ and $\hat l\in\hat\Wl$ such that
$\sig_l<\hat\sig_{\hat l}$ and $\hat l\leq l$ on
$[\sig_l,\hat\sig_{\hat l}]$, there starts at each
$s_+\in\Ii_+(l,\hat l)$ a unique $\hat r^{[s_+]}\in\hat\Wr(l(s_+),s_+)$
such that $l\leq\hat r^{[s_+]}$ on $[s_+-\eps,s_+]$ for some
$\eps>0$. Moreover, the set
\be\label{Iprim}
\Ii'_+(l,\hat l):=\big\{s_+\in\Ii_+(l,\hat l):
\hat r^{[s_+]}\mbox{ crosses $l$ at some time in }(s_-,s_+)\big\}
\ee
is a dense subset of $\Ii(l,\hat l)$. The same statements hold when
$\hat l$ is replaced by a path $\hat r\in\hat\Wr$.
\ep
The proof of Proposition~\ref{P:exent} is somewhat long and depends on
excursion theory. We start by studying the a.s.\ unique left-most path
started at the origin. Let $l\in\Wl(0,0)$ be that path and fix some
deterministic $u>0$. By the structure of special points of the
Brownian web (see, e.g., \cite[Lemma~3.3~(b)]{SS06}),
$\hat\Wl(l(u),u)$ contains a.s.\ two paths, one on each side of
$l$. Let $\hat l$ be the one on the left of $l$. Set
\be
X_t:=\ffrac{1}{\sqrt 2}\big(l(u-t)-\hat l(u-t)\big)\qquad(t\in[0,u]).
\ee
Our first lemma, the proof of which can be found below, says that $X$
is standard Brownian motion reflected at the origin.
\bl\label{L:refBM}{\bf[Reflected Brownian motion]}\\
There exists a standard Brownian motion $B=(B_t)_{t\geq 0}$ such that
\be\label{XPsi}
X_t=B_t+\Psi_t\qquad(t\in[0,u])\quad\mbox{where}\quad
\Psi_t:=-\inf_{s\in[0,t]}B_s\qquad(t\geq 0).
\ee
\el
Extend $X_t$ to all $t\geq 0$ by (\ref{XPsi}), and put
\bc\label{tauincr}
\dis S_\tau&:=&\dis\inf\{t>0:\Psi_t>\tau\}\qquad(\tau\geq 0),\\[5pt]
\dis\Ti&:=&\dis\{\tau\geq 0:S_{\tau-}<S_\tau\}.
\ec
Then the intervals of the form $(S_{\tau-},S_\tau)$ with $\tau\in\Ti$
are precisely the intervals during which $X$ makes an excursion away
{f}rom $0$. Define a random point set $N$ on $(0,\infty)^2$ by
\be\label{Ndef}
N:=\{(h_\tau,\tau):\tau\in\Ti\}\quad\mbox{where}
\quad h_\tau:=S_\tau-S_{\tau-}.
\ee
The following fact is well-known.
\bl\label{L:PoisN}{\bf[Poisson set of excursions]}\\
The set $N$ is a Poisson point process with intensity measure $\nu(\di
x)\di\tau$, where
\be\label{nudef}
\nu(\di h)=\frac{\di h}{\sqrt{2\pi h^3}}.
\ee
\el
{\bf Proof} It follows from Brownian scaling that $(S_\tau)_{\tau\geq
0}$ is a stable subordinator with exponent $1/2$, and this
implies that $\nu(\di h)=ch^{-3/2}\di h$ for some $c>0$. The precise
formula (\ref{nudef}) can be found in \cite[Sect.~6.2.D]{KS91}.\qed

To explain the main idea of the proof of Proposition~\ref{P:exent}, we
formulate one more lemma, which will be proved later.
\bl\label{L:exent}{\bf[Dual paths at top points of excursions]}\\
Almost surely, at each $s_+\in\Ii_+(l,\hat l)$ there starts a unique
$\hat r^{[s_+]}\in\hat\Wr(l(s_+),s_+)$ such that $l\leq\hat r^{[s_+]}$
on $[s_+-\eps,s_+]$ for some $\eps>0$.
\el

Set
\be
\Ti_u:=\{\tau\in\Ti:S_\tau<u\}=\Ti\cap(0,\Psi_u).
\ee
and observe that
\be\label{Utau}
\Ii_+(l,\hat l)=\{U_\tau:\tau\in\Ti_u\}\quad\mbox{where}
\quad U_\tau:=u-S_{\tau-}\qquad(\tau\in\Ti).
\ee
For $d\geq 1$ and $h>0$, let $\Ei^d_h$ be the space of all continuous
functions $f:[0,h]\to\R^d$ such that $f(0)=0$, and set
$\Ei^d:=\{(f,h):h>0,\ f\in\Ei^d_h\}$. Using Lemma~\ref{L:exent},
for each $\tau\in\Ti_u$, we define a random function
$(\ti R^\tau,L^\tau,\hat L^\tau)\in\Ei^3_{h_\tau}$ by
\be
\label{RLL}
\begin{aligned}
\dis\ti R^\tau_t&:=\dis l(U_\tau)
-\big(\hat r^{[U_\tau]}(U_\tau-t)\vee\hat l(U_\tau-t)\big),\\[5pt]
\dis L^\tau_t&:=\dis l(U_\tau)-l(U_\tau-t),\\[5pt]
\dis\hat L^\tau_t&:=\dis l(U_\tau)-\hat l(U_\tau-t),
\end{aligned}
\qquad \qquad 0\leq t\leq h_\tau,
\ee
where $U_\tau$ is as in (\ref{Utau}). Note that modulo translation and
time reversal, the triple $(\ti R^\tau, L^\tau, \hat L^\tau)$
is just $(\hat r^{[U_\tau]}, l, \hat l)$ during the time interval
$[U_\tau-h_\tau, U_\tau]$ when $\hat l$ and $l$ make an excursion
away from each other, and $\hat r^{[U_\tau]}$ coalesces with $\hat l$
upon first hitting $\hat l$. Let $N^3_u$ be the random subset of
$\Ei^3\times(0,\Psi_u)$ defined by
\be\label{Nu}
N^3_u:=\big\{(\ti R^\tau,L^\tau,\hat L^\tau,h_\tau,\tau):\tau\in\Ti_u\big\}.
\ee
We will show that $N^3_u$ can be extended to a Poisson point process $N^3$ on
$\Ei^3\times(0,\infty)$. Proposition~\ref{P:exent} will then be established by
showing that $N^3_u$ contains infinitely many points $(\ti R^\tau,L^\tau,\hat
L^\tau,h_\tau,\tau)$ with the property that $\ti R^\tau$ crosses $L^\tau$
before $h_\tau$. (Note that since we are only interested in the time when
$\hat r^{[U_\tau]}$ crosses $l$, there is no need to follow $\hat
r^{[U_\tau]}$ after it meets $\hat l$. This is why we have defined $\ti
R^\tau_t$ in such a way that it coalesces with $\hat L^\tau_t$.)\med

\noi
{\bf Proof of Lemma~\ref{L:refBM}} Set
\be\ba{r@{\,}c@{\,}l}\label{LLdef}
\dis L_t&:=&\dis l(u)-l(u-t),\\[5pt]
\dis\hat L_t&:=&\dis l(u)-\hat l(u-t)
\ea\qquad(t\in[0,u]).
\ee
We know from \cite{STW00} (see also Lemma~2.1 and
\cite[formula~(6.17)]{SS06}) that conditioned on $l$, the dual path
$\hat l$ is distributed as a Brownian motion with drift $-1$, Skorohod
reflected off $l$. Therefore, on $[0,u]$, the paths $L$ and $\hat L$
are distributed as solutions to the SDE
\bc\label{LhatL}
\dis\di L_t&=&\dis\di B^{\rm l}_t-\di t,\\[5pt]
\dis\di\hat L_t&=&\dis\di B^{\rm\hat l}_t-\di t+\di\Phi_t,
\ec
where $B^{\rm l}_t$ and $B^{\rm\hat l}_t$ are independent, standard
Brownian motions, $\Phi_t$ is a nondecreasing process, increasing only
when $L_t=\hat L_t$, and one has $L_t\leq\hat L_t$ for all
$t\in[0,u]$. Extending our probability space if necessary, we may
extend solutions of (\ref{LhatL}) so that they are defined for all
$t\geq 0$. Set
\be\ba{r@{\,}c@{\,}lr@{\,}c@{\,}l}\label{LLXY}
\dis B^-_t&:=&\dis
\ffrac{1}{\sqrt 2}\big(B^{\rm\hat l}_t-B^{\rm l}_t\big),\qquad&
\dis B^+_t&:=&\dis
\ffrac{1}{\sqrt 2}\big(B^{\rm\hat l}_t+B^{\rm l}_t\big),\\[5pt]
\dis X_t&:=&\dis\ffrac{1}{\sqrt 2}\big(\hat L_t-L_t\big),\qquad&
\dis Y_t&:=&\dis\ffrac{1}{\sqrt 2}\big(\hat L_t+L_t\big).
\ec
Then $B^-_t$ and $B^+_t$ are independent standard Brownian motions,
\bc\label{XYsde}
\dis X_t&=&\dis B^-_t+\Psi_t,\\[5pt]
\dis Y_t&=&\dis B^+_t-\sqrt2t+\Psi_t,
\ec
where $X_t\geq 0$ and $\Psi_t:=\ffrac{1}{\sqrt 2}\Phi_t$ increases
only when $X_t=0$. In particular, setting $B:=B^-$ and noting that
$X$ in (\ref{XYsde}) solves a Skorohod equation, the claims in
Lemma~\ref{L:refBM} then follow.\qed

For each $\tau\in\Ti$, we define $X^\tau\in\Ei_{h_\tau}$ by
\be
X^\tau(t):=X_{S_{\tau-}+t}\quad \quad (t\in[0,h_\tau]),
\ee
and we define a point process $N^1$ on $\Ei^1\times\half$ by
\be
N^1:=\big\{(X^\tau,h_\tau,\tau):\tau\in\Ti\big\}.
\ee
The following facts are well-known.
\bl\label{L:exBM}{\bf[Excursions of reflected Brownian motion]}\\
There exists a \si-finite measure $\mu^1$ on $\Ei^1$ such that $N^1$ is a
Poisson point process on $\Ei^1\times\half$ with intensity measure
$\mu^1(\di(f,h))\di\tau$. The measure $\mu^1$ may be written as
\be
\mu^1(\di(f,h))=\mu^1_h(\di f)\nu(\di h),
\ee
where $\nu$ is the measure in (\ref{nudef}) and the $\mu^1_h$ are
probability measures on $\Ei^1_h$ $(h>0)$. There exists a
random function $F:[0,1]\to\half$ with $F(0)=0=F(1)$ and $F>0$ on $(0,1)$,
such that
\be\label{muh}
\mu^1_h=\P[F_h\in\cdot\,]\quad\mbox{where}\quad
F_h(t):=\sqrt{h}F(t/h)\qquad(t\in[0,h]).
\ee
\el
{\bf Proof} The existence of the excursion measure $\mu^1$ follows from
general excursion theory, see \cite[Chapter~IV]{Ber96} or
\cite[Chapter~VI.8]{RW94}; a precise description of $\mu^1$ for
reflected Brownian motion can be found in \cite[Sect.~VI.55]{RW94}.
Since $\nu$ is the marginal of the measure $\mu^1$, it has to be the
measure in (\ref{nudef}). Formula (\ref{muh}) is a result of Brownian
scaling.\qed

Let $(L,\hat L)$ be the solution to (\ref{LhatL}), extended to all
$t\geq 0$. For each $\tau\in\Ti$, we define $(L^\tau,\hat
L^\tau)\in\Ei^2_{h_\tau}$ by
\be
\begin{aligned}
\dis L^\tau(t)&:=\dis L_{S_{\tau-}+t}-L_{S_{\tau-}},\\[5pt]
\dis\hat L^\tau(t)&:=\dis\hat L_{S_{\tau-}+t}-L_{S_{\tau-}},
\end{aligned}
\qquad \qquad \qquad t\in[0,h_\tau],
\ee
and we define a point process $N^2$ on $\Ei^2\times\half$ by
\be
N^2:=\big\{(L^\tau,\hat L^\tau,h_\tau,\tau):\tau\in\Ti\big\}.
\ee
\bl\label{L:exleft}{\bf[Excursions between a left and dual left path]}\\
The set $N^2$ is a Poisson point process on $\Ei^2\times\half$ with
intensity $\mu^2_h(\di f)\nu(\di h)$, where $\nu$ is the measure in
(\ref{nudef}) and the $\mu^2_h$ are probability measures on $\Ei^2_h$
$(h>0)$ given by
\be
\mu^2_h=\P[(F^-_h,F^+_h)\in\cdot\,]\quad\mbox{with}\quad
F^\pm_h(t):=B_t-t\pm \ffrac{1}{\sqrt 2}F_h(t)\quad(t\in[0,h]),
\ee
where $F_h$ is a random variable as in (\ref{muh}) and $B$ a Brownian
motion independent of $F_h$.
\el
{\bf Proof} This follows from the fact that, by (\ref{LLXY}) and
(\ref{XYsde}),
\bc\label{LLhat}
\dis L_t&=&\dis\ffrac{1}{\sqrt 2}\big(Y_t-X_t\big)=B^+_t-t
+\ffrac{1}{\sqrt 2}\Psi_t-\ffrac{1}{\sqrt 2}X_t,\\[5pt]
\hat L_t&=&\dis\ffrac{1}{\sqrt 2}\big(Y_t+X_t\big)=B^+_t-t
+\ffrac{1}{\sqrt 2}\Psi_t+\ffrac{1}{\sqrt 2}X_t,
\ec
where $B^+$ is a standard Brownian motion independent of $X$ and
$\Psi$. Note that restrictions of $B^+$ to disjoint excursion
intervals are independent and that, since $\Psi$ increases only at
times $t$ when $L_t=\hat L_t$, it drops out of the formulas for
$L^\tau$ and $\hat L^\tau$.\qed

\noi
{\bf Remark} General excursion theory tells us how a strong Markov process can
be constructed by piecing together its excursions from a singleton. In our
situation, however, we are interested in excursions of the process $(L_t,\hat
L_t)$ from the set $\{(x,x): x\in\R\}$, which is not a singleton. Formula
(\ref{LLhat}) shows that apart from motion during the excursions, the process
$(L_t, \hat L_t)$ also moves along the diagonal at times when $L_t=\hat L_t$,
even though such times have zero Lebesgue measure. (Indeed, it is possible to
reconstruct $(L_t, \hat L_t)$ from its excursions and local time in
$\{(x,x): x\in\R\}$, but we do not pursue this here.)\med

Together with Lemma~\ref{L:exleft}, the next lemma implies that the point
process $N^3_u$ on $\Ei^3\times(0,\Psi_u)$ defined in (\ref{Nu}) is a Poisson
point process, as claimed.
\bl\label{L:crosdis}{\bf[Distribution of crossing times]}\\
The paths $(\ti R^\tau)_{\tau\in\Ti_u}$ are conditionally independent given
$l$ and $\hat l$, and their conditional law up to coalescence with $\hat
L^\tau$ is given by the solution to the
Skorohod equation
\be
\begin{aligned}
\di\ti R^\tau_t&=\di B_t+\di t-\di\De_t,\quad\quad&\De_t<T,\\
\di\ti R^\tau_t&=\di B_t+\di t+\di\De_t,\quad\quad &T\leq\De_t,
\end{aligned}
\qquad
0\leq t\leq h_\tau\wedge\inf\{s\geq 0:\ti R^\tau_s=\hat L^\tau_s\},
\ee
where $B$ is a standard Brownian motion, $\De$ is a nondecreasing
process increasing only when $\ti R^\tau_t=L^\tau_t$, $T$ is an independent
mean $1/2$ exponential random variable, and $\ti R^\tau$ is subject to the
constraints that $\ti R^\tau_t\leq L^\tau_t$ resp.\ $L^\tau_t\leq \ti R^\tau_t$
when $\De_t<T$ resp.\ $T\leq\De_t$.
\el
We will prove Lemmas~\ref{L:exent} and \ref{L:crosdis} in one stroke.\med

\noi
{\bf Proof of Lemmas~\ref{L:exent} and \ref{L:crosdis}} We start by showing
that almost surely, at each $s_+\in\Ii_+(l,\hat l)$ there starts at most one
$\hat r^{[s_+]}\in\hat\Wr(l(s_+),s_+)$ such that $l\leq\hat r^{[s_+]}$
on $[s_+-\eps,s_+]$ for some $\eps>0$. Indeed, since dual right-most paths
cannot cross $l$ from left to right \cite[Prop.~3.6~(d)]{SS06}, if there is
more than one dual right-most path starting at $(l(s_+),s_+)$ on the right of
$l$, then there start at least three dual right-most paths at this point. It
follows that $(l(s_+),s_+)$ is a meeting point of $\Wr$ and hence, by
Lemma~\ref{L:meetsep}, also a meeting point of $\Wl$. This contradicts the
existence of an incoming dual left-most path (see Theorem~\ref{T:classweb}).

To prove the other statements we use discrete approximation (compare the proof
of Lemma~\ref{L:forback}). As in \cite{SS06}, we consider systems of
branching-coalescing random walks on $\Z^2_{\rm even}$ with branching
probabilities $\eps_n\to 0$. Diffusively rescaling space and time as
$(x,t)\mapsto(\eps_n x,\eps_n^2 t)$ then yields the Brownian net in the
limit. In the discrete system, we consider the left-most path $l_n$ starting
at the origin, we choose $u_n\in\N$ such that $\eps_nu_n\to u$ and consider
the dual left-most path $\hat l_n$ started at time $u_n$ at distance one to
the left of $l_n$. For each $0\leq i\leq u_n$, we let $i^+:=\inf\{j\geq
i:l_n(j)-\hat l_n(j)=1\}$ and we let $\hat r_n(i)$ denote the position at time
$i$ of the dual right-most path started at time $i^+$ at distance one on the
right of $l_n$. Then $\hat r_n$ is the concatenation of dual right-most paths,
started anew immediately on the right of $l_n$ each time $\hat l_n$ is at
distance one from $l_n$. In analogy with the definition of $\ti R^\tau$ in
(\ref{RLL}), we set $\ti r_n(i):=\hat l_n(i)\vee\hat r_n(i)$.

Now $(l_n,\hat l_n)$, diffusively rescaled, converges in distribution to
$(l, \hat l)$ where $l$ is the left-most path in the Brownian net starting at
the origin and $\hat l$ is the a.s.\ unique dual left-most path starting at
$(l(u),u)$ that lies on the left of $l$. Moreover, the set of times when $\hat
l_n$ is at distance one from $l_n$, diffusively rescaled, converges to the set
$\{s\in(0,u]:\hat l(s)=l(s)\}$. (This follows from the fact that the
reflection local time of $\hat l_n$ off $l_n$ converges to its continuum
analogue, and the latter increases whenever $\hat l(s)=l(s)$.)

In the diffusive scaling limit, the path $\ti r_n$ converges to a path $\ti r$
such that $t\mapsto\ti r(u-t)$ is right-continuous and is set back to $l(u-t)$
each time $\hat l(u-t)$ meets $l(u-t)$. Between these times, in the same way
as in the proof of Lemma~\ref{L:forback}, we see that the conditional law of
$\ti r$ given $l$ and $\hat l$ is as described in Lemma~\ref{L:crosdis}.

Since $\ti r_n$ is the concatenation of dual right-most paths, we see that at
each time $s_+\in\Ii_+(l,\hat l)$ there starts at least one dual right-most
path that lies on the right of $l$ on $[s_+-\eps,s_+]$ for some $\eps>0$,
completing the proof of Lemma~\ref{L:exent}.\qed

\noi
{\bf Proof of Proposition~\ref{P:exent}} We first prove the claims for
the a.s.\ unique paths $l\in\Wl(0,0)$ and $\hat l\in\hat\Wl(l(u),u)$
such that $\hat l$ lies on the left of $l$. For each $\tau\in\Ti_u$, set
\be\label{Ctau}
C_\tau:=h_\tau\wedge\inf\{t\in[0,h_\tau]:L^\tau(t)<\ti R^\tau(t)\}.
\ee
Let $N$ be the Poisson point process in (\ref{Ndef}), let $N_u$ denote the
restriction of $N$ to $(0,\infty)\times(0,\Psi_u)$, and set
\be
N'_u:=\{(h_\tau,\tau):\tau\in\Ti_u,\ C_\tau<h_\tau\}.
\ee
Then $N'_u$ is a thinning of $N_u$, obtained by independently keeping
a point $(h_\tau,\tau)\in N_u$ with probability $\rho(h_\tau)$, where
$\rho:(0,\infty)\to(0,\infty)$ is some function. Indeed, by
Lemmas~\ref{L:exleft} and \ref{L:crosdis}, $\rho(h)$ has
the following description. Pick a random variable $F_h$ as in
Lemma~\ref{L:exBM}, two standard Brownian motions $B,B'$, and a mean
$1/2$ exponential random variable $T$, independent of each other.
Set
\be
L(t):=B_t-t-\ffrac{1}{\sqrt 2}F_h\qquad(t\in[0,h]).
\ee
Let $(R',\De)$ be the solution to the Skorohod equation
\be
\di R'_t=\di B'_t+\di t-\di\De_t\qquad(0\leq t\leq h),
\ee
reflected to the left off $L$, and set
\be
C:=h\wedge\inf\{t\in[0,h]:\De_t\geq T\}.
\ee
Then
\be
\rho(h)=\P[C<h]=\E[1-e^{-2\De_h}]\geq\E[1-e^{-2\De'_h}]=:\rho'(h),
\ee
where
\bc
\dis\De_h&=&\dis\sup_{t\in[0,h]}\big(B'_t+t-L_t\big)
=\sup_{t\in[0,h]}\big(B'_t-B_t+2t+\ffrac{1}{\sqrt 2}F_h(t)\big)\\[5pt]
&\geq&\dis\sup_{t\in[0,h]}\big(B'_t-B_t+\ffrac{1}{\sqrt 2}F_h(t)\big)
=:\De'_h.
\ec
It follows from Brownian scaling (see (\ref{muh})) that
\be
\rho'(h)=\E[1-e^{-2\sqrt{h}\De'_1}],
\ee
hence, since $h^{-1/2}(1-e^{-2\sqrt{h}\De'_1})\up 2\De'_h$ as $h\down 0$,
\be
\lim_{h\to 0}h^{-1/2}\rho'(h)=2\E[\De'_1]>0.
\ee
By (\ref{nudef}), it follows that, for some $c>0$,
\be
\int_{0+}\nu(\di h)\rho(h)\geq c\int_{0+}h^{-3/2}h^{1/2}\di h=\infty,
\ee
i.e., the intensity measure of the thinned Poisson point process is not
integrable, hence the set $\Ti'_u:=\{\tau:(h_\tau,\tau)\in N'_u\}$ is a dense
subset of $(0,\Psi_u)$, hence $\Ii'_+(l,\hat l)$ is dense in $\Ii(l,\hat
l)$. This completes the proof for the special paths $l\in\Wl(0,0)$ and $\hat
l\in\hat\Wl(l(u),u)$ such that $\hat l$ lies on the left of $l$.

To prove the same statement for $l\in\Wl(0,0)$ and $\hat r\in\hat\Wr(l(u),u)$
where $\hat r$ lies on the left of $l$, first note that because $u$ is
deterministic, Lemma \ref{L:forback} implies the existence of such an $\hat r$.
Set
\be
\hat R_t:=l(u)-\hat r(u-t)\qquad(t\in[0,u]),
\ee
and let $L$ be as in (\ref{LLdef}). Then, by Lemma~\ref{L:forback},
$L$ and $\hat R$ are distributed as solutions to the SDE (compare
(\ref{LhatL}))
\bc\label{LhatR}
\dis\di L_t&=&\dis\di B^{\rm l}_t-\di t,\\[5pt]
\dis\di\hat R_t&=&\dis\di B^{\rm\hat r}_t+\di t+\di\Phi_t,
\ec
where $B^{\rm l}_t$ and $B^{\rm\hat r}_t$ are independent, standard
Brownian motions, $\Phi_t$ is a nondecreasing process, increasing only
when $L_t=\hat R_t$, and $L_t\leq\hat R_t$ for all $t\in[0,u]$. By
Girsanov, solutions of (\ref{LhatR}) are equivalent in law to
solutions of (\ref{LhatL}), hence we can reduce this case to the case
of a dual left-most path. In particular, by what we have already
proved, almost surely for each $s_+\in\Ii_+(l,\hat r)$ there exists a
unique $\hat r^{[s_+]}\in\hat\Wr(l(s_+),s_+)$ that lies on the right
of $l$, and the set
\be
\Ii'_+(l,\hat r):=\big\{s_+\in\Ii_+(l,\hat r):
\hat r^{[s_+]}\mbox{ crosses $l$ at some time in }(s_-,s_+)\big\}
\ee
is a dense subset of $\Ii(l,\hat r)$.

By translation invariance, $\Ii'_+(l,\hat l)$ is dense in $\Ii(l,\hat
l)$ for each $l\in\Wl$ started from a point $z=(x,t)\in\Q^2$
and $\hat l\in\hat\Wl(l(u),u)$ that lies on the left of $l$. By
\cite[Lemma~3.4~(b)]{SS06}, we can generalize this to arbitrary
$l\in\Wl$. Since any dual left-most path that hits $l$ from the left
at a time $s$ must have coalesced with some left-most path started in
$(l(u),u)$ for some $u\in\Q$ with $u>s$ and lying on the left of $l$,
we can generalize our statement to arbitrary $l\in\Wl$ and $\hat
l\in\hat\Wl$. The argument for dual right-most paths is the same.\qed

\subsection{Excursions around hitting points}\label{S:excur2}

With the help of Proposition~\ref{P:exent}, we can prove the following result.
Recall the definition of $\Ii(\pi,\hat\pi)$ from (\ref{Iidef}).
\bp\label{P:exint}{\bf[Excursions around intersection points]}\\
Almost surely, for each $l\in\Wl$ and $\hat l\in\hat\Wl$ such that
$\sig_l<\hat\sig_{\hat l}$ and $\hat l\leq l$ on
$[\sig_l,\hat\sig_{\hat l}]$, the set
\be\ba{r@{\,}l}
\dis\Ii''(l,\hat l):=\big\{t\in\Ii(l,\hat l):
&\dis\exists r\in\Wr\mbox{ s.t.\ $r$ makes an excursion from $l$}\\
&\dis\mbox{during a time interval $(s,u)$ with }
\sig_l<s<t<u<\hat\sig_{\hat l}\big\}
\ec
is a dense subset of $\Ii(l,\hat l)$. The same statements hold when
$\hat l$ is replaced by a path $\hat r\in\hat\Wr$.
\ep
{\bf Proof} Since the set $\Ii'_+(l,\hat l)$ from (\ref{Iprim}) is
dense in $\Ii(l,\hat l)$, for each $t\in\Ii(l,\hat l)\beh\Ii_+(l,\hat
l)$ we can find $s^{(n)}_+\in\Ii'_+(l,\hat l)$ such that $s^{(n)}_+\up
t$. Our claim will follow provided we show that infinitely many of the
$s^{(n)}_+$ are in $\Ii''(l,\hat l)$. It suffices to show that at
least one $s^{(n)}_+$ is in $\Ii''(l,\hat l)$; then the same argument
applied to the sequence started after $s^{(n)}_+$ gives the existence
of another such point, and so on, ad infinitum. So imagine that
$s^{(n)}_+\not\in\Ii''(l,\hat l)$ for all $n$. Let $z_n=(l(c_n),c_n)$,
with $c_n\in(s^{(n)}_-,s^{(n)}_+)$, be the point where $\hat
r^{[t_n]}$ crosses $l$. By Proposition~\ref{P:sep}, there exists an
incoming path $r_n\in\Wr$ at $z_n$ such that $r_n$ separates from $l$
in $z_n$. Set $\tau_n:=\inf\{s>c_n:l(s)=r_n(s)\}$. If
$\tau_n<s^{(n)}_+$ then $l$ and $r_n$ form a wedge of $(\Wl,\Wr)$
which cannot be entered by $\hat r^{[t_n]}$, leading to a
contradiction. If $\tau_n=s^{(n)}_+$, then $s^{(n)}_+$ is a meeting
point of $l$ and $r_n$, hence by Lemma~\ref{L:meetsep},
$(l(s^{(n)}_+),s^{(n)}_+)$ is of type $(2,1)/(0,3)$ in $\Wl$, which
contradicts the existence of the dual incoming path $\hat l$. Finally,
we cannot have $\tau_n\in(s^{(n)}_+,\hat\sig_{\hat l})$ because of our
assumption that $s^{(n)}_+\not\in\Ii''(l,\hat l)$, so we conclude that
$\tau_n\geq\sig_{\hat l}$ for all $n$. It follows that the $z_n$ are
$(\sig_l,\hat\sig_{\hat l})$-relevant separation points, contradicting
Lemma~\ref{L:finrel}.\qed

\section{Structure of special points}\label{S:struct}

\subsection{Classification of special points}\label{S:struct1}

Recall the preliminary classification of points in the Brownian net
given in Section~\ref{S:image4}, which is based only on the structure
of incoming paths. In this section, we turn our attention to
the more detailed classification from Theorem~\ref{T:classnet}, which
also uses information about outgoing paths that are not continuations
of incoming paths. We start with a preliminary lemma.
\bl\label{L:noin}{\bf[No incoming paths]}\\
Almost surely for each $z=(x,t)\in\R^2$, there is no incoming path
$\pi\in\Ni$ at $z$ if and only if $\hat\Wr(z)$ and $\hat\Wl(z)$ each
contain a single path, $\hat r$ and $\hat l$, say, and $\hat r
\sim^z_{\rm out}\hat l$.
\el
{\bf Proof} Let $\hat r$ be the left-most element of $\hat\Wr(z)$ and
let $\hat l$ be the righ-most element of $\hat\Wl(z)$. Then $\hat
r\leq\hat l$, and by Lemma~\ref{L:NT}, there is an incoming path
$\pi\in\Ni$ at $z$ if and only if $\hat r<\hat l$ on $(t-\eps,t)$ for
some $\eps>0$.\qed

Theorem~\ref{T:classnet} follows from the following result.
\bt\label{T:classnet2}{\bf[Classification of points in the Brownian net]}\\
Let $\Ni$ be the standard Brownian net and let $\hat\Ni$ be its
dual. Then, using the classification of Definition~\ref{D:imgpts}, almost
surely, each point in $\R^2$ is of one of the following 19 types in
$\Ni/\hat\Ni$:
\begin{itemize}
\item[{\rm (1)}] ${\rm(C_o)/(C_o)}$, ${\rm(C_o)/(C_p)}$, ${\rm(C_p)/(C_o)}$,
${\rm(C_o)/(C_m)}$, ${\rm(C_m)/(C_o)}$, ${\rm(C_p)/(C_p)}$;
\item[{\rm (2)}] ${\rm(C_s)/(C_s)}$;
\item[{\rm (3)}] ${\rm(C_l)/(C_o)}$, ${\rm(C_o)/(C_l)}$, ${\rm(C_r)/(C_o)}$,
${\rm(C_o)/(C_r)}$;
\item[{\rm (4)}] ${\rm(C_l)/(C_p)}$, ${\rm(C_p)/(C_l)}$, ${\rm(C_r)/(C_p)}$,
${\rm(C_p)/(C_r)}$;
\item[{\rm (5)}] ${\rm(C_l)/(C_l)}$, ${\rm(C_r)/(C_r)}$;
\item[{\rm (6)}] ${\rm(C_n)/(C_o)}$, ${\rm(C_o)/(C_n)}$;
\end{itemize}
and all of these types occur. Moreover, these points correspond to the
types listed in Theorem~\ref{T:classnet} (in the same order), where
points of type ${\rm(C_p)/(C_p)}$ are either of type
${\rm(p,pp)_l/(p,pp)_l}$ or ${\rm(p,pp)_r/(p,pp)_r}$. For each
deterministic time $t\in\R$, almost surely, each point in
$\R\times\{t\}$ is of either type ${\rm(C_o)/(C_o)}$,
${\rm(C_o)/(C_p)}$, or ${\rm(C_p)/(C_o)}$, and all of these types
occur. A deterministic point $(x,t)\in\R^2$ is almost surely of type
${\rm(C_o)/(C_o)}$.
\et
For clarity, we split the proof into three lemmas.
\bl\label{L:fordu}{\bf[Forward and dual types]}\\
Almost surely, for all $z\in\R^2$:
\begin{itemize}
\item[\rm(a)] If $z$ is of type ${\rm(C_s)}$ in $\Ni$, then $z$ is of
type ${\rm(C_s)}$ in $\hat\Ni$.
\item[\rm(b)] If $z$ is of type ${\rm(C_m)}$ in $\Ni$, then $z$ is of
type ${\rm(C_o)}$ in $\hat\Ni$.
\item[\rm(c)] If $z$ is of type ${\rm(C_n)}$ in $\Ni$, then $z$ is of
type ${\rm(C_o)}$ in $\hat\Ni$.
\item[\rm(d)] If $z$ is of type ${\rm(C_l)}$ in $\Ni$, then $z$ is not of
type ${\rm(C_r)}$ in $\hat\Ni$.
\end{itemize}
In particular, each point in $\R^2$ is of one of the 19 types listed in
Theorem~\ref{T:classnet2}.
\el
\bl\label{L:exis}{\bf[Existence of types]}\\
Almost surely, all types of points listed in Theorem~\ref{T:classnet2}
occur. For each deterministic time $t\in\R$, almost surely, each point
in $\R\times\{t\}$ is of either type ${\rm(C_o)/(C_o)}$,
${\rm(C_o)/(C_p)}$, or ${\rm(C_p)/(C_o)}$, and all of these types
occur. A deterministic point $(x,t)\in\R^2$ is almost surely of type
${\rm(C_o)/(C_o)}$.
\el
\bl\label{L:struc}{\bf[Structure of points]}\\
Almost surely, with respect to $\Ni/\hat\Ni$ and for all $z\in\R^2$:
\begin{itemize}
\item[\rm(a)] If $z$ is of type ${\rm(C_o)/(C_o)}$, then $z$ is of
type ${\rm(o,p)/(o,p)}$.
\item[\rm(b)] If $z$ is of type ${\rm(C_p)/(C_o)}$, then $z$ is of
type ${\rm(p,p)/(o,pp)}$.
\item[\rm(c)] If $z$ is of type ${\rm(C_m)/(C_o)}$, then $z$ is of
type ${\rm(pp,p)/(o,ppp)}$.
\item[\rm(d)] If $z$ is of type ${\rm(C_p)/(C_p)}$, then $z$ is of
type ${\rm(p,pp)_l/(p,pp)_l}$ or ${\rm(p,pp)_r/(p,pp)_r}$.
\item[\rm(e)] If $z$ is of type ${\rm(C_s)/(C_s)}$, then $z$ is of
type ${\rm(p,pp)_s/(p,pp)_s}$.
\item[\rm(f)] If $z$ is of type ${\rm(C_l)/(C_o)}$, then $z$ is of
type ${\rm(l,p)/(o,lp)}$.
\item[\rm(g)] If $z$ is of type ${\rm(C_l)/(C_p)}$, then $z$ is of
type ${\rm(l,pp)_r/(p,lp)_r}$.
\item[\rm(h)] If $z$ is of type ${\rm(C_l)/(C_l)}$, then $z$ is of
type ${\rm(l,lp)_r/(l,lp)_r}$.
\item[\rm(i)] If $z$ is of type ${\rm(C_n)/(C_o)}$, then $z$ is of
type ${\rm(n,p)/(o,lr)}$.
\end{itemize}
Analogous statements hold for the remaining types in Theorem \ref{T:classnet2}
by left-right and forward-backward symmetry.
\el
{\bf Proof of Lemma~\ref{L:fordu}} Part~(a) follows from
Proposition~\ref{P:sep}. If $z$ is of type ${\rm(C_m)}$ in $\Ni$, then
by Lemma~\ref{L:meetsep}, there is a single outgoing left-right pair
at $z$, hence by Lemma~\ref{L:noin}, $z$ is of type ${\rm(C_o)}$ in
$\hat\Ni$, proving part~(b). If $z$ is of type ${\rm(C_n)}$ in $\Ni$,
then each outgoing Brownian net path at $z$ must pass through the top
points of the nested excursions around $z$ as described in
Proposition~\ref{P:imgclas}~(a). Hence also in this case, there is a
single outgoing left-right pair at $z$, so by Lemma~\ref{L:noin}, $z$
is of type ${\rm(C_o)}$ in $\hat\Ni$, proving part~(c). To prove
part~(d), suppose that $z$ is of type ${\rm(C_l)}$ and that $\hat
r\in\hat\Wr$ enters $z$. By Proposition~\ref{P:imgclas}~(a), there
exist right-most paths $r_n$ making a sequence of nested excursions
away from $l$. Since each $r_n$ forms with $l$ a wedge of $(\Wl,\Wr)$
that cannot be entered by $\hat r$, the latter must satisfy $\hat
r\leq l$ on $(t,t+\eps)$ for some $\eps>0$. Since $\hat r$ reflects
off $l$, we can find some $u>t$ such that $\hat r(u)<l(u)$. Now any
path in $\hat\Wl$ started in $(\hat r(u),l(u))\times\{u\}$ must enter
$z$, hence $z$ is not of type ${\rm(C_r)}$ in $\hat\Ni$. It is easy
to check that (a)--(d) rules out all but the 19 types listed in
Theorem~\ref{T:classnet2}.
\qed

\noi
{\bf Proof of Lemma~\ref{L:exis}} It follows from Lemma~\ref{L:imgpts}
that each deterministic $z\in\R^2$ is of type ${\rm(C_o)/(C_o)}$, and
that for each deterministic $t\in\R$, a.s.\ each point in
$\R\times\{t\}$ is of one of the types ${\rm(C_o)/(C_o)}$,
${\rm(C_p)/(C_o)}$, ${\rm(C_o)/(C_p)}$, and ${\rm(C_p)/(C_p)}$. By the
structure of the Brownian web (see, e.g., \cite[Lemma~3.3]{SS06}),
points of type ${\rm(C_p)/(C_p)}$ do not occur at deterministic
times. Since points of type ${\rm(C_p)}$ do occur at deterministic
times both in $\Ni$ and in $\Ni$, we conclude that a.s.\ each point in
$\R\times\{t\}$ is of one of the types ${\rm(C_o)/(C_o)}$,
${\rm(C_p)/(C_o)}$, and ${\rm(C_o)/(C_p)}$, and all of these types
occur.

To prove the existence of all 19 types of points listed in
Theorem~\ref{T:classnet2}, by symmetry between left and right and
between forward and dual paths, it suffices to prove the existence of
the 9 types of points listed in Lemma~\ref{L:struc}. We have just
established the existence of points of types ${\rm(C_o)/(C_o)}$ and
${\rm(C_p)/(C_o)}$. The existence of points of types
${\rm(C_s)/(C_s)}$, ${\rm(C_m)/(C_o)}$, and ${\rm(C_n)/(C_o)}$ follows
{f}rom Lemma~\ref{L:fordu} and the existence of points of types
${\rm(C_s)}$, ${\rm(C_m)}$, and ${\rm(C_n)}$ in $\Ni$. Hence, we are
left with the task of establishing the existence of points of types
${\rm(C_p)/(C_p)}$, ${\rm(C_l)/(C_o)}$, ${\rm(C_l)/(C_p)}$, and
${\rm(C_l)/(C_l)}$.

By the structure of the Brownian web, there exist $l\in\Wl$ and $\hat
l\in\hat\Wl$ such that $l\leq\hat l$ on $(\sig_l,\hat\sig_{\hat l})$
and the set $\Ii(l,\hat l)$ defined in (\ref{Iidef}) is not
empty. Since $\Ii(l,\hat l)$ is uncountable and since the set of all
crossing points in $\R^2$ is countable, there exist lots of points
$z=(x,t)\in\Ii(l,\hat l)$ such that no path in $\hat\Wr$ crosses $l$
and no path in $\Wr$ crosses $\hat l$ at $z$. For such points, we can
find $r\in\Wr$ and $\hat r\in\hat\Wr$ such that $l\leq r\leq\hat
r\leq\hat l$ on $[t-\eps,t+\eps]$ for some $\eps>0$. Of all types of
points listed in Theorem~\ref{T:classnet2}, only points of type
${\rm(C_p)/(C_p)}$ have incoming paths in $\Wl$, $\Wr$, $\hat\Wr$, and
$\hat\Wl$, hence, by Lemma~\ref{L:fordu}, $z$ must be of this type.

We are left with the task of establishing the existence of points of
types ${\rm(C_l)/(C_o)}$, ${\rm(C_l)/(C_p)}$, and ${\rm(C_l)/(C_l)}$.
By Lemma~\ref{L:iso} and Proposition~\ref{P:imgclas}~(b), a point $z$
is of type ${\rm(C_l)}$ in $\Ni$ if and only if there is an incoming
path $l\in\Wl$ at $z$ and there are right-most paths $r_n$ making a
nested sequence of excursions from $l$, as described in
Proposition~\ref{P:imgclas}~(b). We need to show that we can choose
these nested excursions in such a way that $z$ is of type
${\rm(C_o)}$, ${\rm(C_p)}$, or ${\rm(C_l)}$ in $\hat\Ni$.

Fix a path $l\in\Wl$ and let $\{\hat l_n\}_{n\in\N}$ be paths in
$\hat\Wl$ starting from a deterministic countable dense subset of
$\R^2$. Since $\hat l_1$ is reflected off $l$, using
Lemma~\ref{L:exc}, we can find a path $r_1\in\Wr$ such that $r_1$
makes an excursion from $l$ during an interval $(s_1,u_1)$ with
$l_1\neq l$ on $[s_1,u_1]$. By the same arguments, we can inductively
find paths $r_n\in\Wr$ such that $r_n$ makes an excursion from $l$
during an interval $(s_n,u_n)$ with $l_n\neq l$ on
$[s_n,u_n]\sub(s_{n-1},u_{n-1})$. We can choose the $r_n$ such that
$\bigcap_n[s_n,u_n]=\{t\}$ for some $t\in\R$. Then setting
$z:=(l(t),t)$ yields a point such that $z$ is of type ${\rm(C_l)}$ in
$\Ni$ and no path in $\hat\Wl$ enters $z$, hence (by
Lemma~\ref{L:fordu}~(d)), $z$ must be of type ${\rm(C_l)/(C_o)}$.

To construct points of type ${\rm(C_l)/(C_p)}$, we fix paths $l\in\Wl$
and $\hat r\in\hat\Wr$ with $\hat r\leq l$ on $[\sigma_l, \hat\sigma_{\hat r}]$
such that the set $\Ii(l,\hat r)$ defined in
(\ref{Iidef}) is not empty. By Proposition~\ref{P:exint} and
Lemma~\ref{L:exc}, we can inductively find paths $r_n\in\Wr$ such that
$r_n$ makes an excursion from $l$ during an interval $(s_n,u_n)$ with
$[s_n,u_n]\sub(s_{n-1},u_{n-1})$ and $(s_n,u_n)\cap\Ii(l,\hat
r)\neq\emptyset$. Choosing $\{t\}=\bigcap_n[s_n,u_n]$ and setting
$z:=(l(t),t)$ then yields a point of type ${\rm(C_l)}$ in $\Ni$ such
that $\hat r$ enters $z$, hence (by Lemma~\ref{L:fordu}~(d)), $z$ must
be of type ${\rm(C_l)/(C_p)}$.

Finally, to construct points of type ${\rm(C_l)/(C_l)}$, we fix paths
$l\in\Wl$ and $\hat l\in\hat\Wr$ with $\hat l\leq l$ on $[\sigma_l,
  \hat\sigma_{\hat l}]$ such that the set $\Ii(l,\hat l)$ defined in
(\ref{Iidef}) is not empty. By Proposition~\ref{P:exint} and
Lemma~\ref{L:exc}, we can inductively find intervals $(s_n,u_n)$ with
$[s_n,u_n]\sub(s_{n-1},u_{n-1})$ and $(s_n,u_n)\cap\Ii(l,\hat
l)\neq\emptyset$, and paths $r_{2n}\in\Wr$ and $\hat r_{2n+1}\in\hat\Wr$ such
that $r_{2n}$ makes an excursion from $l$ during $(s_{2n},u_{2n})$ and $\hat
r_{2n+1}$ makes an excursion from $\hat l$ during
$(s_{2n+1},u_{2n+1})$. Choosing $\{t\}=\bigcap_n[s_n,u_n]$ and setting
$z:=(l(t),t)$ then yields a point of type ${\rm(C_l)/(C_l)}$.\qed

\noi
{\bf Proof of Lemma~\ref{L:struc}} Part~(e) about separation points
has already been proved in Proposition~\ref{P:sep}. It follows from
Lemma~\ref{L:meetsep} that at points of all types except ${\rm(C_m)}$,
all incoming paths in $\Ni$ are equivalent. Therefore, points of type
${\rm(C_o)}$ are of type ${\rm(o,\cdot)}$, points of type ${\rm(C_n)}$
are of type ${\rm(n,\cdot)}$, points of type ${\rm(C_l)}$ are of type
${\rm(l,\cdot)}$, and points of types ${\rm(C_p)}$ are of type
${\rm(p,\cdot)}$, while by Lemma~\ref{L:meetsep}, points of type
${\rm(C_m)}$ are of type ${\rm(pp,p)}$. We next show that, the configuration
of incoming paths in $\Ni$ at $z$ determines the configuration of outgoing
paths in $\hat\Ni$ at $z$:
\begin{itemize}
\item[\rm(o)] If $z$ is of type ${\rm(C_o)}$ in $\Ni$, then $z$ is of type
${\rm(\cdot,p)}$ in $\hat\Ni$.
\item[\rm(p)] If $z$ is of type ${\rm(C_p)}$ in $\Ni$, then $z$ is of type
${\rm(\cdot,pp)}$ in $\hat\Ni$.
\item[\rm(m)] If $z$ is of type ${\rm(C_m)}$ in $\Ni$, then $z$ is of type
${\rm(\cdot,ppp)}$ in $\hat\Ni$.
\item[\rm(l)] If $z$ is of type ${\rm(C_l)}$ in $\Ni$, then $z$ is of type
${\rm(\cdot,lp)}$ in $\hat\Ni$.
\item[\rm(n)] If $z$ is of type ${\rm(C_n)}$ in $\Ni$, then $z$ is of type
${\rm(\cdot,lr)}$ in $\hat\Ni$.
\end{itemize}

Statement~(o) follows from Lemma~\ref{L:noin}.

To prove statement~(p), we observe that if $z=(x,t)$ is of type ${\rm(C_p)}$,
then by Proposition~\ref{P:imgclas}~(e), for each $T<t$ with $z\in N_T$, there
exist maximal $T$-meshes $M^-_T=M(r',l)$ and $M^+_T=M(r,l')$ with bottom times
strictly smaller than $t$ and top times strictly larger than $t$, such that
$l\sim^z_{\rm in}r$ and $l\sim^z_{\rm out}r$. Therefore, by the structure of
the Brownian web and ordering of paths, there exist $\hat l^-,\hat
l^+\in\hat\Wl(z)$ and $\hat r^-,\hat r^+\in\hat\Wr(z)$ such that $r'\leq\hat
r^-\leq\hat l^-\leq l\leq r\leq\hat r^+\leq\hat l^+\leq l'$ on $[t-\eps,t]$,
for some $\eps>0$. By Lemma~\ref{L:bex}, the paths $\hat l^-$ and $\hat r^-$
both pass through the bottom point of $M^-_T$. Since $T$ can be chosen
arbitrarily close to $t$, we must have $\hat l^-\sim^z_{\rm out}\hat r^-$, and
similarly $\hat l^+\sim^z_{\rm out}\hat r^+$.

The proof of statement~(m) is similar, where in this case we use the
three maximal $T$-meshes from Proposition~\ref{P:imgclas}~(d).

If $z=(x,t)$ is of type ${\rm(C_l)}$, then by the structure of the
Brownian web and ordering of paths, there exist $\hat l^-,\hat
l^+\in\hat\Wl(z)$ such that $\hat l^-\leq l\leq\hat l^+$, and there
exists a unique $\hat r\in\hat\Wr(z)$ with $\hat r\leq \hat l^-$.
By Lemma~\ref{L:iso}, $z$ is
isolated from the left in $N_T$ for any $T<t$ with $z\in N_T$, hence by
Proposition~\ref{P:tmaxmesh}~(b), there exists a maximal $T$-mesh
$M(r,l)$ with bottom time strictly smaller than $t$ and top time
strictly larger than $t$, such that $l(t)=x$. Now any path in
$\hat\Ni$ on the left of $l$ must exit $M(r,l)$ through its bottom
point, hence the same argument as before shows that $\hat
l^+\sim^z_{\rm out}\hat r$, proving statement~(l).

Finally, if $z=(x,t)$ is of type ${\rm(C_n)}$, then by the
structure of the Brownian web, $\hat\Wr(z)$ and $\hat\Wl(z)$ each
contain a unique path, say $\hat r$ and $\hat l$. By
Lemma~\ref{L:noin}, $\hat r\not\sim^z_{\rm out}\hat l$, so by ordering
of paths, $\hat r<\hat l$ on $(t-\eps,t)$ for some $\eps>0$.
This proves statement~(n).

To complete our proof, we must determine how incoming paths continue
for points of the types ${\rm(C_p)/(C_p)}$, ${\rm(C_l)/(C_p)}$, and
${\rm(C_l)/(C_l)}$.

Points of type ${\rm(C_p)/(C_p)}$ must be either of type
${\rm(p,pp)_l/(p,pp)_l}$ or ${\rm(p,pp)_r/(p,pp)_r}$. Hence, points of
at least one of these types must occur, so by symmetry between left and
right, both types must occur.

If $l\in\Wl$ and $\hat l\in\hat\Wl$ are incoming paths at a point $z$
and $l\leq\hat l$, then it has been shown in the proof of
Lemma~\ref{L:exis} that $z$ must either be of type ${\rm(C_p)/(C_p)}$
or of type ${\rm(C_s)/(C_s)}$. It follows that at points of type
${\rm(C_l)/(C_p)}$ and ${\rm(C_l)/(C_l)}$, the incoming path $l\in\Wl$
continues on the right of the incoming path $\hat l\in\hat\Wl$.\qed

\subsection{Structure of special points}\label{S:struct2}

In this section, we prove Theorems~\ref{T:locspec} and
\ref{T:netpt}. We start with a preparatory lemma.
\bl\label{L:refexc}{\bf[Excursions between reflected paths]}\\
Let $\hat\pi_1,\hat\pi_2\in\hat\Ni$, $z_i=(x_i,t_i)\in\R^2$ with
$t_i<\hat\sig_{\hat\pi_i}$ and $x_i\leq\hat\pi_i(t_i)$, and let
$r_i=r_{z_i,\hat\pi_i}$ $(i=1,2)$ be reflected right-most paths off
$\hat\pi_1$ and $\hat\pi_2$. Assume that $r_1$ makes an excursion from $r_2$
during a time interval $(s,u)$ with
$u<\hat\sig_{\hat\pi_1},\hat\sig_{\hat\pi_2}$. Then there exists a $t\in(s,u]$
such that $\hat\pi_1$ separates from $\hat\pi_2$ at time $t$.
\el
{\bf Proof} Without loss of generality we may assume that $r_1<r_2$ on
$(s,u)$. Since $r_1$ separates from $r_2$ at time $s$, by the structure of
reflected paths (Lemma~\ref{L:extl}), we must have $s\in\Fi(r_1)$ but
$s\not\in\Fi(r_2)$, and hence $\hat\pi_1(s)=r_1(s)=r_2(s)<\hat\pi_2(s)$. By
the structure of separation points and Lemma~\ref{L:NT}~(a), we have
$r_1\leq\hat\pi_1\leq r_2$ on $[s,s+\eps]$ for some $\eps>0$. Set
$\tau:=\inf\{t>u:r_2(t)<\hat\pi_1(t)\}$. Then $\tau\leq u$ since $r_1$ does
not cross $\hat\pi_1$ and spends zero Lebesgue time with $\hat\pi_1$ (see
Lemma~\ref{L:netint}). By Lemma~\ref{L:lpi}~(a) and the structure of reflected
paths, we have $\tau\in\Fi(r_2)$, hence $\hat\pi_1(\tau)=\hat\pi_2(\tau)$. It
follows that $\hat\pi_1$ and $\hat\pi_2$ separate at some time
$t\in(s,\tau]$.\qed
\bigskip

\noi
{\bf Proof of Theorem~\ref{T:locspec}~(a)} To prove part~(a), assume that
$l_i=l_{z_i,\hat r_i}\in\El_{\rm in}(z)$ $(i=1,2)$ satisfy $l_1\sim^z_{\rm
  in}l_2$ but $l_1\neq^z_{\rm in}l_2$. Then $l_1$ and $l_2$ make excursions
away from each other on a sequence of intervals $(s_k,u_k)$ with $u_k\up t$.
By Lemma~\ref{L:refexc}, it follows that $\hat r_i\in\hat\Wr_{\rm out}(z)$
$(i=1,2)$ separate at a sequence of times $t_k\up t$, which contradicts
(\ref{wequiv}). If $l_i=l_{z,\hat r_i}\in\El_{\rm out}(z)$ $(i=1,2)$ satisfy
$l_1\sim^z_{\rm out}l_2$ but $l_1\neq^z_{\rm out}l_2$, then the same argument
shows that the paths $\hat r_i\in\hat\Er_{\rm in}(z'_i)$ with $z'_i=(x'_i,t)$
and $x'_i\leq x$ separate at a sequence of times $t_k\down t$. It follows that
$z'_1=z'_2=:z'$, $r_1\sim^{z'}_{\rm in}r_2$, but $r_1\neq^{z'}_{\rm in}r_2$,
contradicting what we have just proved.\qed

Before we continue we prove one more lemma.
\bl\label{L:lrcont}{\bf[Containment between extremal paths]}\\
Almost surely, for each $z=(x,t)\in\R^2$, if $l$ is the left-most element of
$\Wl(z)$ and $r$ is the right-most element of $\Wr$, then every $\pi\in\Ni(z)$
satisfies $l\leq\pi\leq r$ on $[t,\infty)$ and every
$\hat\pi\in\hat\Ni_{\rm in}(z)$ satisfies $l\leq\hat\pi\leq r$ on
$[t,\hat\sig_{\hat\pi}]$.
\el
{\bf Proof} By symmetry, it suffices to prove the statements for $l$. The
first statement then follows by approximation of $l$ with left-most paths
starting on the left of $x$, using the fact that paths in $\Ni$ cannot
paths in $\Wl$ from right to left \cite[Prop.~1.8]{SS06}. The second statement
follows from Lemma~\ref{L:NT}~(a).\qed

\noi
{\bf Proof of Theorems~\ref{T:locspec}~(b)--(c) and \ref{T:netpt}~(a)--(c)}.
Since these are statements about incoming paths only, or about outgoing paths
only, it suffices to consider the following three cases: {\bf Case I }points
of type $\rm(\cdot,p)/(o,\cdot)$, {\bf Case II }points of types
$\rm(\cdot,pp)/(p,\cdot)$ and $\rm(\cdot,ppp)/(pp,\cdot)$, and {\bf Case III}
points of types $\rm(\cdot,lp)/(l,\cdot)$, $\rm(\cdot,pr)/(r,\cdot)$, or
$\rm(\cdot,lr)/(n,\cdot)$. In each of these cases, when we prove
Theorem~\ref{T:netpt}~(a), we actually prove the analogue  for
$\hat\Ni_{\rm in}(z)$.

{\bf Case I }In this case, by Lemma~\ref{L:lrcont}, all paths in $\Ni_{\rm
  out}(z)$ are contained between the equivalent paths $l$ and $r$. Note that
this applies in particular to paths in $\El_{\rm out}(z)$ and $\Er_{\rm
  out}(z)$, so by Theorem~\ref{T:locspec}~(a), $\El_{\rm out}(z)=\Wl_{\rm
  out}(z)$ and $\El_{\rm out}(z)=\Wl_{\rm out}(z)$ up to strong
equivalence. By Lemma~\ref{L:NT}~(a), $\hat\Ni_{\rm in}(z)=\emptyset$.

{\bf Case II }We start with points of type $\rm(\cdot,pp)/(p,\cdot)$. Write
$\Wl(z)=\{l,l'\}$ and $\Wr=\{r,r'\}$ where $l\sim^z_{\rm out}r'$ and
$l'\sim^z_{\rm out}r$. By Lemma~\ref{L:starchar}, $r'$ and $l'$ form a maximal
$t$-mesh, hence, by Proposition~\ref{P:tmaxmesh}~(b) and Lemma~\ref{L:lrcont},
all paths in $\Ni_{\rm out}(z)$ are either contained between $l$ and $r'$ or
between $l'$ and $r$. By Theorem~\ref{T:locspec}~(a), this implies that
$\El_{\rm out}(z)=\{l,l'\}$ and $\Er_{\rm out}(z)=\{r,r'\}$ up to strong
equivalence. To prove the statements about dual incoming paths, we note that
by Lemma~\ref{L:struc}, points of type $\rm(\cdot,pp)/(p,\cdot)$ are of type
$\rm(C_p)$ or $\rm(C_s)$ in $\hat\Ni$. Therefore, by
Proposition~\ref{P:imgclas}~(c) and (e), there are unique paths $\hat
r\in\hat\Wr_{\rm in}(z)$ and $\hat l\in\hat\Wl_{\rm in}(z)$ and all paths in
$\hat\Ni_{\rm in}(z)$ are eventually contained between the equivalent paths
$\hat r$ and $\hat l$. By Theorem~\ref{T:locspec}~(a), $\hat\El_{\rm
  in}(z)=\{\hat l\}$ and $\hat\Er_{\rm in}(z)=\{\hat r\}$ up to strong
equivalence.

For points of type $\rm(\cdot,ppp)/(pp,\cdot)$ the argument is similar,
except that there are now two maximal $t$-meshes with bottom point $z$, and
for the dual incoming paths we must use Proposition~\ref{P:imgclas}~(d).

{\bf Case III }Our proof in this case actually also works for points of type
$\rm(\cdot,pp)/(p,\cdot)$, although for these points, the argument given in
Case~II is simpler.

Let $l$ denote the left-most element of $\Wl(z)$ and let $r$ denote the
right-most element of $\Wr(z)$. By Lemma~\ref{L:NT}~(a), any path
$\pi\in\hat\Ni_{\rm in}(z)$ must stay in $V:=\{(x_+,t_+):t_+>t,\ l(t_+)\leq
x_+\leq r(t_+),\ l<r\mbox{ on }(t,t_+)\}$. Choose any $z_+\in V$ and define
(by Lemma~\ref{L:extl}) reflected paths by $\hat r:=\hat
r_{z_+,l}=\min\{\hat\pi'\in\Ni(z_+):l\leq\hat\pi'\mbox{ on }[t,t_+]\}$ and
$\hat l:=\hat l_{z_+,r}$. By Lemma~\ref{L:meetsep}, all paths in
$\hat\Ni_{\rm in}(z)$ are equivalent, hence $\hat r\sim^z_{\rm in}\hat l$ and,
by Theorem~\ref{T:locspec}~(a), $\hat\El_{\rm in}(z)=\{\hat l\}$ and
$\hat\Er_{\rm in}(z)=\{\hat r\}$ up to strong equivalence. Since our
construction does not depend on the point $z_+$, each path
$\hat\pi\in\hat\Ni_{\rm in}(z)$ satisfies $\hat r\leq\hat\pi\leq\hat l$ on
$[t,t+\eps]$ for some $\eps>0$.

Now fix some $z_+=(x_+,t_+)\in V$ and put $\hat r:=\hat r_{z_+,l}$ and $\hat
l:=\hat l_{z_+,r}$. By what we have just proved we can choose $z_+$ in such a
way that $\hat l\in\hat\Wl_{\rm in}(z)$ if the latter set is nonempty and
$\hat r\in\hat\Wl_{\rm in}(z)$ if that set is nonempty. In any case
$\hat\El_{\rm in}(z)=\{\hat l\}$ and $\hat\Er_{\rm in}(z)=\{\hat r\}$ up to
strong equivalence. Define (by Lemma~\ref{L:extl}) reflected paths
$r':=r_{z,\hat l}$ and $l':=l_{z,\hat r}$. Since any path in $\hat
l'\in\hat\El_{\rm in}$ is strongly equivalent to $\hat l$, each reflected path
of the form $r_{z,\hat l'}$ is strongly equivalent to $r'$. On the other hand,
it is easy to see that each reflected path of the form $r_{z,\hat l'}$ with
$\hat l'\in\El_{\rm in}(z')$, $z'=(x',t')$, $x<x'$ is strongly equivalent to
$r$. It follows that, up to strong equivalence, $\Er_{\rm out}(z)=\{r,r'\}$,
and by a symmetric argument $\El_{\rm out}(z)=\{l,l'\}$.

Since $\hat r\sim^z_{\rm in}\hat l$, since $\hat r$ and $\hat l$ spend
positive Lebesgue time together whenever they meet while $r'$ and $\hat l$
spend zero Lebesgue time together, we must have $r'(t_n)<\hat r(t_n)$ for a
sequence of times $t_n\down t$. Since $r'$ can cross $\hat r'$ only at times
in $\Fi(\hat r)$, we must have $l(t'_n)=r'(t'_n)$ for a sequence of times
$t'_n\down t$, i.e., $l\sim^z_{\rm out}r'$. A symmetric argument shows that
$l'\sim^z_{\rm out}r$. It follows from the definition of reflected paths that
any $\pi\in\Ni(z)$ such that $\pi\leq\hat l$ on $[t,t_+]$ must satisfy
$l\leq\pi\leq r'$ on $[t,t_+]$. Since all paths in $\hat\El_{\rm in}(z)$ are
strongly equivalent, this shows that if $\pi\in\Ni(z)$ satisfies
$\pi\leq\hat l'$ on $[t,t+\eps]$ for some $\hat l'\in\El_{\rm in}(z)$ and
$\eps>0$, then $l\leq\pi\leq r'$ on $[t,t+\eps']$ for some $\eps'>0$. A
similar conclusion holds if $\hat r'\leq\pi$ on $[t,t+\eps]$ for some $\hat
r'\in\Er_{\rm in}(z)$ and $\eps>0$.

By Lemma~\ref{L:lpi}~(a), a path $\pi\in\Ni(z)$ can cross $\hat l$ from
right to left only at times in $\Fi(\hat l)$ and $\pi$ can cross $\hat r$
{f}rom left to right only at times in $\Fi(\hat r)$. In particular, if $\hat
l\in\hat\Wl$, then any path $\pi\in\Ni(z)$ that enters the region to the right
of $\hat l$ must stay there till time $t_+$. It follows that for points of
type $\rm(\cdot,lp)/(l,\cdot)$ or $\rm(\cdot,pp)/(p,\cdot)$, any
$\pi\in\Ni(z)$ must either satisfy $\hat l\leq\pi$ on $[t,t_+]$ or
$\pi\leq\hat l$ on $[t,t+\eps]$ for some $\eps>0$. By what we have just
proved, this implies that either $l'\leq\pi\leq r$ on $[t,t_+]$ or
$l\leq\pi\leq r'$ on $[t,t+\eps']$ for some $\eps'>0$. By symmetry, a similar
argument applies to points of type $\rm(\cdot,pr)/(r,\cdot)$.

This completes the proof of Theorems~\ref{T:locspec}~(b)--(c) and
\ref{T:netpt}~(a)--(b). To prove Theorem~\ref{T:netpt}~(c), we must show that
at points of type $\rm(o,lr)/(n,p)$ there exist paths $\pi$ that infinitely
often cross over between $l$ and $r$ in any time interval $[t,t+\eps]$ with
$\eps>0$. Since $r'\not\in\Wr_{\rm out}(z)$, we have $\Fi(r')=\{t_n:n\geq 0\}$
for some $t_n\down t$. Let $r_n\in\Wr$ be the right-most path such that
$r'=r_n$ on $[t_{n+1},t_n]$ and set $\tau_n:=\inf\{u>t_n:r_n(u)=r(u)\}$.  By
the compactness of $\Wr$ we have $r_n\up r'$ for some $r'\in\Wr(z)$ and hence,
since the latter set has only one element, $r_n\up r$. This shows that
$\tau_n\down t$. By a symmetry, there exist times $t'_n,\tau'_n\down t$ and
$l_n\in\Wl$ such that $l_n=l'$ on $[t'_{n+1},t'_n]$ and
$l_n(\tau'_n)=l(\tau'_n)$. Using these facts and the hopping construction of
the Brownian net, it is easy to construct the desired path $\pi\in\Ni(z)$ such
that $l\sim^z_{\rm out}\pi\sim^z_{\rm out}r$.\qed
\bigskip

\noi {\bf Proof of Theorem~\ref{T:netpt}~(d)} At all points of types with the
subscript ${\rm l}$, there is a dual path $\hat r$ that, by
Lemma~\ref{L:lpi}~(a), cannot be crossed by any path $\pi\in\Ni$ entering $z$
(as we have have just seen it must) between the unique incoming paths
$l\in\Wl_{\rm in}(z)$ and $r\in\Wr_{\rm in}(z)$, and hence prevents $\pi$ from
leaving $z$ in the right outgoing equivalence class. By symmetry, an analogue
statement holds for points of types with the subscript ${\rm r}$.

In all other cases where $\Ni_{\rm in}(z)$ is not empty, either the point is
of type $\rm(p,pp)_s$ or there is a single outgoing equivalence class. In all
these cases, using the fact that the Brownian net is closed under hopping, it
is easy to see that any concatenation of a path in $\Ni_{\rm in}(z)$ up to
time $t$ with a path in $\Ni_{\rm out}(z)$ after time $t$ is again a path in
$\Ni$.\qed

\appendix

\section{Extra material}\label{A:extra}

\subsection{Some simple lemmas}\label{A:extra1}

\bl\label{L:lrincome}{\bf[Incoming paths]}\\
Let $(\Wl,\Wr,\hat\Wl,\hat\Wr)$ be the standard left-right Brownian
web and its dual. For $z\in\R^2$, Let $m^{\rm lr}_{\rm in}(z)$ and
$m^{\rm lr}_{\rm out}(z)$ denote respectively the number of
equivalence classes of paths in $(\Wl, \Wr)$ entering and leaving
$z$. Let $\hat m^{\rm lr}_{\rm in}(z)$ and $\hat m^{\rm lr}_{\rm
out}(z)$ be defined similarly for $(\hat\Wl, \hat\Wr)$. Then
\begin{itemize}
\item[{\rm(a)}] For each deterministic $z\in\R^2$, almost surely
$m^{\rm lr}_{\rm in}(z)+\hat m^{\rm lr}_{\rm in}(z)=0$.

\item[{\rm(b)}] For each deterministic $t\in\R$, almost surely $m^{\rm
lr}_{\rm in}(x,t)+\hat m^{\rm lr}_{\rm in}(x,t)\leq 1$ for all
$x\in\R$.

\item[{\rm(c)}] Almost surely, $m^{\rm lr}_{\rm in}(z)+\hat m^{\rm
lr}_{\rm in}(z)\leq 2$ for all $z\in\R^2$.
\end{itemize}
\el
{\bf Proof} By Lemma~3.4~(b) of~\cite{SS06}, it suffices to consider
only paths in $(\Wl, \Wr, \hat\Wl, \hat\Wr)$ starting from a
deterministic countable dense set $\Di \subset \R^2$, which we denote
by
\[
(\Wl(\Di), \Wr(\Di),\hat\Wl(\Di),\hat\Wr(\Di)).
\]
Part~(a) follows from the fact that a Brownian motion (with constant
drift) almost surely does not hit a deterministic space-time point.

Part (b) follows from the fact, paths in $\Wl(\Di)\cup\Wr(\Di)$,
resp.\ $\hat\Wl(\Di)\cup\hat\Wr(\Di)$, evolve independently before
they meet, hence almost surely, no two paths in
$\Wl(\Di)\cup\Wr(\Di)$, resp.\ $\hat\Wl(\Di)\cup\hat\Wr(\Di)$, can
first meet at a deterministic time. Also note that,
$(\Wl(\Di),\Wr(\Di))|_{(-\infty, t)}$, resp.\ $(\Wl(\Di),
\Wr(\Di))|_{(t,\infty)}$, the restriction of paths in $(\Wl(\Di),
\Wr(\Di))$ to the time interval $(-\infty, t)$, resp.\ $(t, \infty)$,
are independent. Hence it follows that $(\Wl(\Di),
\Wr(\Di))|_{(-\infty,t)}$ and $(\hat\Wl(\Di),
\hat\Wr(\Di))|_{(t,\infty)}$ are also independent, and almost surely,
no path in $\Wl(\Di)\cup\Wr(\Di)$ starting before time $t$ can be at
the same position at time $t$ as some path in
$\hat\Wl(\Di)\cup\hat\Wr(\Di)$ starting after time $t$.

The proof of (c) is similar. By symmetry, it suffices to show that:
(1) the probability of $\pi_1, \pi_2,\pi_3 \in\Wl(\Di)\cup\Wr(\Di)$
first meeting at the same space-time point is zero; (2) the
probability of $\pi_1, \pi_2 \in\Wl(\Di)\cup\Wr(\Di)$ first meeting at
a point $(x,t)$, and $\hat \pi(t)=x$ for some $\hat\pi \in\hat\Wl(\Di)
\cup\hat\Wr(\Di)$ starting after time $t$ is zero. Claim (1) holds
because $\pi_1, \pi_2$ and $\pi_3$ evolve independently before they
meet. For claim (2), note that conditioned on $T=\inf\{t>
\sig_{\pi_1}\vee\sig_{\pi_2}:\pi_1(t)=\pi_2(t)\}$, the
distribution of
$(\Wl(\Di\cap((T,\infty)\times\R)),\Wr(\Di\cap((T,\infty)\times\R)))$
is independent of $(\pi_1, \pi_2)|_{(-\infty, T)}$, which is a
consequence of the strong Markov property of a collection of
left-right coalescing Brownian motions. Hence $(\hat\Wl(\Di),
\hat\Wr(\Di))|_{(T, \infty)}$ is also independent of $(\pi_1,
\pi_2)|_{(-\infty, T)}$ conditioned on $T$. The claim then follows.
\qed

\bl\label{L:intpt}{\bf[Intersection points of $\Wl$ and $\Wr$]}\\
For $l\in\Wl$ and $r\in\Wr$ with deterministic starting points $z_{\rm
l}$, resp.\ $z_{\rm r}$, let $T = \inf\{ s\geq \sig_{l} \vee
\sig_r : l(s) = r(s)\}$, $I=\{ s > T: l(s) =r(s)\}$, and let $\mu_I$
denote the measure on $\R$ defined by $\mu_I(A) = \ell(I\cap A)$,
where $\ell$ is the Lebesgue measure. Then
\begin{itemize}
\item[{\rm (a)}] Almost surely on the event $T<\infty$,
$\lim_{t\downarrow 0} t^{-1} \mu_I([T, T+t]) = 1$.

\item[{\rm (b)}] $\P(l(T+t) = r(T+t)\, |\, T<\infty) \uparrow 1$ as
$t\downarrow 0$.

\end{itemize}
\el
{\bf Proof} By the strong Markov property of the unique weak solution
of (\ref{lrsde}) (see Proposition~2.1 of~\cite{SS06}), it suffices to
verify (a) and (b) for the solution of (\ref{lrsde}) with initial
condition $L_0=R_0=0$. Let $W_\tau = \sqrt{2} B_\tau + 2\tau$, where
$B_\tau$ is a standard Brownian motion. Let
\be
X_{\tau} = W_\tau + R_\tau \quad \mbox{with} \quad R_\tau
= -\inf_{0\leq s\leq\tau} W_s
\ee
denote $W_\tau$ Skorohod reflected at 0. Then recall from the proof of
Proposition~2.1 and Lemma~2.2 in \cite{SS06} that, the process
$D_t=R_t-L_t$ is a time change of $X_\tau$, converting the local time
of $X$ at 0 into real time. More precisely, $D_t$ is equally
distributed with $X_{T_t}$, where the inverse of $T$ is defined by
$T^{-1}_\tau = \tau + \frac{1}{2}R_\tau$. In particular, $\int_0^t
1_{\{X_{T_s}=0\}}ds = \frac{1}{2}R_{T_t}$. Therefore to show part~(a),
it suffices to show that, almost surely
\be
\lim_{t\down 0}\frac{R_{T_t}}{2t} =\lim_{\tau\down 0}
\frac{1}{1+ \frac{2\tau}{R_\tau}} = 1.
\ee
By Girsanov's theorem, it suffices to show that
\be\label{subor}
\lim_{\tau\down 0}\frac{\tau}{R'_\tau}=0\quad{\rm a.s.},
\ee
where $R'_\tau = -\sqrt{2} \inf_{0\leq s\leq \tau} B_s$ for a standard
Brownian motion $B_s$. This can be proved by a straightforward
Borel-Cantelli argument, which we leave to the reader.

To prove part~(b), it suffices to show the monotonicity of
$\P[L_s=R_s]$ in $s$, because by (a), we have $\lim_{t\down 0}
\frac{1}{t}\int_0^t\P[L_s=R_s]ds=\lim_{t\downarrow
0}\E\big[\frac{1}{t}\int_0^t 1_{\{L_s=R_s\}}ds\big]=1$. Note that by
the same arguments as in the proof of Proposition~16 in \cite{SS06},
$D_t=R_t-L_t$ is the unique nonnegative weak solution of the SDE
\be\label{Dsde}
dD_t=1_{\{D_t>0\}}\sqrt{2}\,dB_t+2dt,
\ee
which defines a time homogeneous strong Markov process with continuous
paths in $[0,\infty)$. By coupling, it is clear that any such Markov
process preserves the stochastic order, i.e., if $D_s$ and $D'_s$ are
two solutions of (\ref{Dsde}) with initial laws
$\Li(D_0)\leq\Li(D'_0)$, then $\Li(D_t)\leq\Li(D'_t)$ for all $t\geq 0$,
where $\leq$ denotes stochastic order. In particular, if $D_0=0$, then
$\Li(D_s)\leq\Li(D_t)$ for all $s\leq t$, since
$\de_0\leq\Li(D_{t-s})$. Therefore $\P(D_s=0) \geq\P(D_t=0)$ for all
$s\leq t$, which is the desired monotonicity.\qed

\subsection{Meeting points}\label{A:extra2}

In this appendix, we give an alternative proof of the fact that any
meeting point of two paths $\pi,\pi'\in\Wl\cup\Wr$ is of type
$\rm(pp,p)$. Recall that any pair $L\in\Wl$ and $R\in\Wr$ with
deterministic starting points solve the SDE (\ref{lrsde}). Consider
the SDE (\ref{LRL}). Let us change variables and denote $X_t=
L'_t-L_t$, $Y_t=R_t-L_t$, so
\bc\label{XY}
\dis\di X_t&=&\dis 1_{\txt\{Y_t>0\}}(\di B^{\rm l'}_t -\di B^{\rm l}_t)
+1_{\txt\{Y_t=0\}}(\di B^{\rm l'}_t-\di B^{\rm s}_t), \\[5pt]
\dis\di Y_t&=&\dis 1_{\txt\{Y_t>0\}}(\di B^{\rm r}_t - \di B^{\rm l}_t)
+ 2\di t,
\ec
with the constraint that $Y_t\geq 0$. In terms of the process
$(X,Y)$, the stopping time $\tau$ from (\ref{lrltau}) becomes
$\tau=\inf\{t\geq 0: X_t=Y_t\}$, and (\ref{aim2}) becomes
\be\label{aim}
\lim_{\eps\to 0} \P^{\eps,0}[Y_\tau =0] = 1.
\ee
The generator of the process $(X,Y)$ is
\bc
G = \diff{x} + 1_{\{y>0\}}\diff{y} + 1_{\{y>0\}}\difif{x}{y} + 2\dif{y}.
\ec
Define a function $g(x,y)$ on the domain $0\leq y<x$ by
\be
g(x,y)=\frac{y}{x},
\ee
so that
\be\label{Gg} Gg(x,y) = \frac{2y}{x^3} - 1_{\{y>0\}}\frac{1}{x^2}
+ \frac{2}{x} = 1_{\{y>0\}}\frac{2y-x}{x^3} +\frac{2}{x}.
\ee
Now define the domain $\De=\{(x,y): 0\leq y < x/2,\ 0<x< 1\}$. Then
the first term on the right-hand side of (\ref{Gg}) is nonpositive on
$\De$. We compensate the second term on the right-hand side of
(\ref{Gg}) by adding another function to $g$. Namely, let
$f=8\sqrt{x}$, then $Gf(x,y) = -2x^{-3/2}$. Hence
\be
G(g+f)\leq 0\quad\mbox{on}\quad\De,
\ee
which means that $(g+f)(X_t,Y_t)$ is a local supermartingale before it
exits $\De$. Let $\partial\De$ denote the boundary of $\De$. Observe
that
\be\ba{l}
g+f\geq 0\quad\mbox{on}\quad\De, \\[5pt]
g+f\geq\ffrac{1}{2}\quad\mbox{on}\quad \partial
\De\backslash\{(0,0)\}, \\[5pt] \dis\lim_{\eps\to 0}\,(g+f)(\eps,0)=0.
\ec
It follows that as $\eps\downarrow 0$, the probability that the
process $(X,Y)$ started at $(\eps,0)$ exits $\De$ from
$\partial\De\backslash\{(0,0)\}$ tends to zero. Consequently, the
process must exit from $(0,0)$, and (\ref{aim}) and (\ref{aim2}) then
follow.

Now we prove that any meeting point of two paths
$\pi,\pi'\in\Wl\cup\Wr$ is of type $\rm(pp,p)$. By Lemma~3.4~(b) of
\cite{SS06}, it suffices to show this for paths in $(\Wl, \Wr)$
starting from a deterministic countable dense subset $\Di\subset\R^2$.
Without loss of generality, assume $l_1, l_2\in\Wl$ start respectively
{f}rom $z_1=(x_1,t_1), z_2=(x_2, t_2)\in\Di$ with $t_1\leq t_2$ and
$l_1(t_2) < l_2(t_2)$. Let $\tau_\eps = \inf\{s\geq t_2 :
l_2(s)-l_1(s) \leq \eps\}$. Then for any $\delta>0$, we can find
$z'_1\in\Di$ and $r_1\in\Wr(z'_1)$ with $z'_1$ sufficiently close to
$(l_1(\tau_\eps), \tau_\eps)$, such that with probability at least
$1-\delta$, $l_1$ and $r_1$ will meet, at which time $l_1$ and $l_2$
still have not met and are at most $2\eps$ distance apart. Since
$\eps, \delta>0$ can be arbitrary, (\ref{aim2}) implies that almost
surely, the first meeting point $z$ of $l_1$ and $l_2$ is also a
meeting point of $l_2$ and some $r_1\in\Wr(\Di)$ bounded between $l_1$
and $l_2$. In fact, by Lemma~\ref{L:lrincome}~(c), we have
$l_1\sim^z_{\rm in} r_1$. Repeating the argument for $r_1$ and $l_2$,
there must exist $r_2\in\Wr(\Di)$ which enters $z$ and $r_2\sim^z_{\rm
in} l_2$. The case of meeting points of $l\in\Wl(\Di)$ and
$r\in\Wr(\Di)$ is similar.\qed

\bigskip
\noi
{\bf\large Acknowledgement} We thank Jon Warren for useful email conversations.
We also thank the referees for constructive criticism.

\newcommand{\noopsort}[1]{}

\parbox[t]{4.7cm}{\small
Emmanuel Schertzer\\
Mathematics Department\\
Columbia University\\
Room 509, MC 4406\\
2990 Broadway\\
New York, NY 10027\\
e-mail: schertzer@math.columbia.edu}
\hspace{.3cm}
\parbox[t]{4.7cm}{\small
Rongfeng~Sun\\
Department of Mathematics\\
National University\\
of Singapore\\
2 Science Drive 2\\
117543 Singapore\\
e-mail: matsr@nus.edu.sg}
\hspace{.3cm}
\parbox[t]{4.7cm}{\small
Jan M.~Swart\\
Institute of Information\\ Theory and Automation\\ of
the ASCR (\' UTIA)\\
Pod vod\'arenskou v\v e\v z\' i 4\\
18208 Praha 8\\
Czech Republic\\
e-mail: swart@utia.cas.cz}

\end{document}